\documentclass[11 pt, reqno]{article}

\usepackage{microtype}
\usepackage{amsthm}
\usepackage{amsmath}
\usepackage{amssymb}
\usepackage{amsfonts}
\usepackage{latexsym}
\usepackage{hyperref}
\usepackage{bbm}
\usepackage{comment} 
\usepackage{enumitem}
\usepackage{cite}
\usepackage{braket}
\usepackage{graphicx}
\usepackage{chngcntr}
\counterwithin*{equation}{section}
\usepackage{tikz-cd}
 \usepackage{wrapfig}
 \usepackage{bm}

\usepackage[a4paper, left=2.5cm, right=2.5cm, top=2.5 cm, bottom= 2.5 cm]{geometry}
\usepackage{color}
\usepackage{fancyhdr}
\usepackage{latexsym}
\newtheorem{ccounter}{ccounter}[section]
\newtheorem{thm}[ccounter]{Theorem}
\newtheorem{lem}[ccounter]{Lemma}
\newtheorem{cor}[ccounter]{Corollary}
\newtheorem{defn}[ccounter]{Definition}
\newtheorem{prop}[ccounter]{Proposition}
\newtheorem{ass}[ccounter]{Assumption}
\newtheorem{ex}[ccounter]{Example}

\def\bet{\begin{thm}}
\def\eet{\end{thm}}
\def\bel{\begin{lem}}
\def\eel{\end{lem}}
\def\bas{\begin{ass}}
\def\eas{\end{ass}}
\def\bec{\begin{cor}}
\def\eec{\end{cor}}
\def\bed{\begin{defn}}
\def\eed{\end{defn}}
\def\bep{\begin{prop}}
\def\eep{\end{prop}}
\def\beq{\begin{equation}}
\def\eeq{\end{equation}}
\def\proof{\noindent {\bf Proof.}\ \ }
\def\bea{\begin{equation*}}
\def\eea{\end{equation*}}
\def\tr{\mathrm{tr}}
\def\bex{\begin{ex}}
\def\eex{\end{ex}}
\def\bp{\proof}
\def\ep{\qed \\}
\def\be{\beq}
\def\ee{\eeq}

\def\remark{\noindent{\bf Remark.}\ \ }

\def\benr{\begin{enumerate}[label=(\roman*)]}
\def\eenr{\end{enumerate}}

\def\N{\mathbb{N}}
\def\C{\mathbb{C}}
\def\T{\mathbb{T}}
\def\Z{\mathbb{Z}}
\def\R{\mathbb{R}}
\def\P{\mathbb{P}}
\def\E{\mathbb{E}}
\def\S{\mathbb{S}}
\def\B{\mathcal{B}}
\def\mfc{m_{\mathrm{fc}, t } }

\def\sc{\rho_{\mathrm{sc}} }
\def\msc{m_{\mathrm{sc}} }
\def\gsc{\gamma^{\mathrm sc}}
\def\one{{\mathbbm 1}}
\def\eps{\varepsilon}
\def\dsc{\rho_{\mathrm{fc}, t }}
\def\Im{\operatorname{Im}}
\def\Re{\operatorname{Re}}

\newcommand{\dd}{\mathrm{d}}
\newcommand{\bma}{\begin{bmatrix}}
\newcommand{\ema}{\end{bmatrix}}
\def\trans{\dagger}
\def\tr{\operatorname{Tr}}

\newcommand{\abs}[1]{\lvert #1 \rvert}
\newcommand{\absa}[1]{\left\lvert #1 \right\rvert}

\newcommand{\norm}[1]{\lVert #1 \rVert}
\newcommand{\norma}[1]{\left\lVert #1 \right\rVert}
\def\U{{\mathcal U}}
\def\B{{\mathcal B}}
\def\Fa{{\mathcal F_\alpha}}
\def\l{\ell}
\def\hc{\hat{\mathcal C}}
\def\flat{\gamma^{\mathrm f}}
\def\hz{\hat{z}}

\begin{document}

%%\maketitle

%\begin{comment}
\begin{table}
\centering
\begin{tabular}{c}
\multicolumn{1}{c}{\Large{\bf Universality of the least singular value for sparse random matrices}}\\
\\
\\
\end{tabular}
\begin{tabular}{c c c}
Ziliang Che & & Patrick Lopatto \\
\\
\multicolumn{3}{c}{ \small{Department of Mathematics} } \\
 \multicolumn{3}{c}{ \small{Harvard University} } \\
 \\
 
\small{zche@math.harvard.edu} & & \small{lopatto@math.harvard.edu} \\
\\
\end{tabular}
\\
\begin{tabular}{c}
\multicolumn{1}{c}{\today}\\
\\
\end{tabular}

\begin{tabular}{p{15 cm}}
\small{{\bf Abstract:} We study the distribution of the least singular value associated to an ensemble of sparse random matrices. Our motivating example is the ensemble of $N\times N$ matrices whose entries are chosen independently from a Bernoulli distribution with parameter $p$. These matrices represent the adjacency matrices of random Erd\H{o}s--R{\'e}nyi digraphs and are sparse when $p\ll 1$. We prove that in the regime $pN\gg 1$, the distribution of the least singular value is universal in the sense that it is independent of $p$ and equal to the distribution of the least singular value of a Gaussian matrix ensemble. We also prove the universality of the joint distribution of multiple small singular values. Our methods extend to matrix ensembles whose entries are chosen from arbitrary distributions that may be correlated, complex valued, and have unequal variances.}
\end{tabular}
\end{table}
%\end{comment} %\begin{comment}
{
\tableofcontents
}

\section{Introduction} \label{sec:int}
Random real symmetric and complex Hermitian matrices have been intensely studied since Wigner's discovery that, in the large $N$ limit, their eigenvalue densities are universal and follow the semicircle distribution. Recent investigations have culminated in a proof of the Wigner--Dyson--Mehta conjecture, which asserts the universality of the local eigenvalue statistics in the limit \cite{EPR10, wigfixed, EYbook}. {\let\thefootnote\relax\footnote{Z.C. is partially supported by NSF grant DMS-1607871. P.L. is partially supported by the NSF Graduate Research Fellowship Program under Grant DGE-1144152.}

In the case of non-Hermitian matrices, there has been a similar study of the singular values. Given a $N\times N$ matrix $M$, its singular values are the eigenvalues of $\sqrt{M^\trans M}$, which we label
\beq 0 \le \lambda_1 \le \dots \le \lambda_N. \eeq
Traditionally, one studies the squares of the singular values, the eigenvalues of $M^\trans M$. In the bulk, the limiting distribution of the squares of the singular values is universal under fairly general hypotheses and follows the Marchenko--Pastur law \cite{MP, BSbook, sparsecovar}. Averaged-energy universality for the local correlation functions around a fixed energy in the bulk was shown in \cite{PY,localrelax1}, and \cite{PY} also showed universality for the largest singular value $\lambda_N$. However, the methods in \cite{PY} do not suffice to prove universality of the least singular value. The analysis of this case is more subtle because the Marchenko--Pastur distribution has a density with a singularity at the origin. 

Early work on the least singular value of random matrices was motivated by the analysis of algorithms in computer science. An important recent development in this area is the method of smoothed analysis, which estimates the practical performance of algorithms \cite{ST02}. In these applications, it is important to estimate the probability that $\lambda_1$ is small for various random matrix models. Since the inverse of the least singular value of a matrix is equal to the operator norm of its inverse, these estimates control the probability that the inverse has large norm.  We refer the reader to \cite{RV} for an introduction to this line of research. Recent results include \cite{TV09a} on Bernoulli matrices, \cite{C16} on structured random matrices, \cite{TV10a} on matrices with i.i.d. entries shifted by a deterministic matrix, and \cite{LR12, BR15} on sparse matrices. 

The universality of the least singular value distribution was considered in \cite{leastsv} from a viewpoint inspired by the method of property testing in computer science and combinatorics. The main result is that if $\xi$ is a real random variable with $\E \xi =0$ and $\E \xi^2 =1$, and such that $\E |\xi|^C < \infty$ for some sufficiently large absolute constant $C$, then the distribution of the least singular value of the ensemble of $N\times N$ random matrices $M_N$ with entries chosen i.i.d. with distribution $\xi/\sqrt{N}$ satisfies 
\beq \P\left( N\lambda_1(M_N) \le  r \right) = 1 - e^{-r^2/2 -r} + O(N^{-c}),\eeq
for some $c>0$.
Also given in \cite{leastsv} are analogous results for complex matrices and for the joint distribution of multiple smallest singular values. However, these results require that the entries are independent and have equal variances. 

 This paper studies the universality of the least singular value from the same dynamical viewpoint that was used to prove the Wigner--Dyson--Mehta conjecture. Our motivating example is the ensemble of $N\times N$ matrices whose entries are chosen independently from a Bernoulli distribution with parameter $p$. These matrices represent the adjacency matrices of random Erd\H{o}s--R{\'e}nyi digraphs, which are directed graphs on $N$ vertices where each possible directed edge is present with probability $p$. Such matrices are sparse when $p\ll 1$, and our result implies that the distribution of the least singular value is universal in the regime $pN\gg 1$. We also apply our method to prove universality of the least singular value for matrices whose entries have unequal variances and weak correlations.
 
An important feature of our proof is that we consider an analogue of Dyson Brownian motion where the particles move in the Weyl chamber corresponding to the hyperoctahedral group. This is in contrast to the literature on the universality of eigenvalue statistics, which studies the traditional Dyson Brownian motion with particles restricted to the Weyl chamber corresponding to the symmetric group. For background on Brownian motion in a Weyl chamber, we refer the reader to \cite{G99}. An essential technical input is showing that these dynamics, which govern the evolution of the singular value distribution, reach local equilibrium after short times $t \gg N^{-1}$.\\

\noindent{\bf Acknowledgments.} The authors thank H.-T. Yau for suggesting the problem, useful discussions, and helpful comments on our preliminary draft. P.L. thanks A. Aggarwal, B. Landon, and \mbox{P. Sosoe} for useful discussions. Finally, the authors thank the referee for their comments and suggestions, which significantly improved the paper.

\section{Overview and main result}\label{sec:main}

In this section we state and prove our main theorem on the universality of the least singular value for sparse matrices, invoking several preliminary results proved in the forthcoming sections. We adapt the three-step method \cite{localrelax1, localrelax2, survey,EYbook} and take advantage of recent technical improvements \cite{fixed, HLY15, LY}. 
\mbox{Step 1} is to obtain local control of the singular values of a deformed ensemble, and in particular to prove their rigidity, which is necessary for the stochastic analysis in the second step. {Step 2} is to obtain short time universality. We prove that for a sparse matrix, the least singular value is universal after time $t=N^{-1+\eps}$ when the singular values are evolved according to the singular value analogue of Dyson Brownian motion. {Step 3} is to remove the time evolution and prove universality for the original ensemble.

Section \ref{sec:deterministic} contains the main estimate of this paper, the analysis of the time evolution of the singular values, which is necessary for Step 2. In Section \ref{sec:dqclaw} we carry out Step 1. Step 3 is accomplished in Section \ref{sec:removal} by a Green function comparison theorem. In Section \ref{a:correlated} we discuss how to extend our methods to matrices with correlated entries and unequal variances, and in Appendix \ref{a:sde} we verify the SDE governing the evolution of the singular values is well-posed.

We now define the primary matrix model we study in this work. To motivate this definition, consider an $N\times N$ matrix $X$ of independent Bernoulli random variables $\{x_{ij}\}$ that take the value $1$ with probability $p$ and $0$ with probability $1-p$ for some parameter $p$, which may depend on $N$. Heuristically, in order to place the bulk eigenvalues on an interval of constant order, we must normalize the matrix by the average $\ell^2$ sum of a row, which in this case is $\sqrt{pN}$. Setting $q = \sqrt{pN}$, $M=X/q$ and extracting the mean $f$ of the resulting entries, we may decompose the normalized matrix as 
\beq\nonumber
	{M} = B + f \ket{w}\bra{w}
\eeq
where $f = \sqrt{p}$. The parameter $q$ also plays a role in $B$, where the $k$-th moment of each entry is $O(N^{-1}q^{2-k})$. 

We generalize this setting and come to the following definition.

\bed\label{def:sparse} We say a sequence of matrices $(M_N)_{N=1}^\infty$ is a sparse random matrix ensemble with sparsity parameter $q$ and mean $f$ if, for every $N$, $M_N$ is a $N\times N$ matrix of the form
\beq M_N = B_N + f \ket{w}\bra{w}\eeq
where $w= N^{-1/2}(1,\dots, 1)^\trans$, $f$ is a parameter (which may depend on $N$) such that $0\le f \le N^{1/2}$, and $B_N$ is a real matrix with independent entries $b_{ij}$ such that \beq\label{growth} \E[b_{ij}] =0,\quad E[b_{ij}^2] = s^{(N)}_{ij}, \quad \E[|b_{ij}|^k] \le \frac{C^k}{Nq^{k-2}}\eeq
for all $k$. We assume there exists $\alpha>0$ such that $q$ satisfies 
\beq N^\alpha \le q \le N^{1/2}. \eeq 
We also assume the variance matrices $S_N$ are doubly stochastic with elements $s_{ij}^{(N)}$ satisfying
\beq \frac{c}{N} \le s^{(N)}_{ij} \le \frac{C}{N}\eeq
for some universal constants $c$ and $C$.
\eed

\remark  We have chosen to state our results for the model in Definition \ref{def:sparse} to simplify the exposition. However, several generalizations are possible.

First, the assumption that the variance matrices $S_N$ are doubly stochastic is made for convenience, so that we are in the usual case where the limiting global distribution of the singular values is a quarter circle. We prove universality in the technically more involved case of matrices with correlated entries in Section \ref{a:correlated}, and our work there subsumes the case of matrices with independent entries and a non-stochastic variance matrix. 

Second, the condition \eqref{growth} may be relaxed considerably to require bounds on only a finite number of moments. These moment bounds are used to prove certain stochastic domination estimates in Section \ref{sec:removal} by applying Markov's inequality to large moments of the entries (see Definition \ref{d:sd}). However, while Definition \ref{d:sd} requires a certain estimate to hold for all  $D$ and $\eps > 0$, our proof requires this estimate only for a fixed large $D$ and fixed small $\eps >0$ (independent of the size of the matrix $N$). Hence, the bound  \eqref{growth} is needed only for some large but fixed number of moments, and the condition can be weakened to requiring that
\beq  \E[|b_{ij}|^k] \le \frac{C}{Nq^{k-2}} \eeq
for some constant $C$ for all $k\le K_0$, where $K_0$ is some large constant. In particular, our result extends to prove universality for random matrices with entries of the form $b_{ij} \xi_{ij}$, where $\{ b_{ij} \}$ are independent Bernoulli random variables with parameter $p = N^{-1 + \delta}$ for any $\delta >0$  and and $\{\xi_{ij}\}$ are i.i.d. random variables with all moments finite. 

Finally, the results of \cite{EKYY13} suggest universality of the least singular value should hold for matrices with sparsity parameter as small as $q = (\log N)^C$ for some large constant $C$ (and even this is probably not optimal, in light of some recent results on sparse matrices \cite{dumitriu2018sparse, he2018local} and the relaxation time of Dyson Brownian motion in the Hermitian case \cite{claeys2017boundaries}). Our approach would generalize to this case if we could show short time universality for $t= (\log N)^C/N$. However, we show this only for $t = N^{-1+\eps}$. It is plausible that a careful examination of our proof would yield the stronger result, but we do not take this up here. \\

Our main result shows that the least singular value of a sparse random matrix ensemble is universal in the large $N$ limit. The Wigner ensemble case of this result was proved in Theorem 1.3 of \cite{leastsv}.

\bet\label{thm:mainresult} Let $(M_N)_{N=1}^\infty$ be a sparse matrix ensemble with least singular values $\lambda_1(M_N)$. For all $r\ge 0$, we have 
\beq \P\left( N\lambda_1(M_N) \le  r \right) = 1 - e^{-r^2/2 -r} + O(N^{-c}) \eeq
where $c>0$ is an absolute constant uniform in $r$.
\eet

\begin{comment}
\begin{figure}
  \centering
    \includegraphics[width=0.72\textwidth]{Figure1}
    \caption{Plot of the limiting probability density function $f(x) = (x+1)\exp( - x^2/2 - x)$,} which corresponds to the cumulative distribution function $\P\left( N\lambda_1(M_N) \le  x \right) = 1 - e^{-x^2/2 -x} + O(N^{-c})$ in Theorem \ref{thm:mainresult}.
\end{figure}
\end{comment}

Since our short time universality estimate Theorem \ref{thm:homomain} holds not just for the smallest singular value, but also for the smallest $k$ singular values when $k$ is at most a small power of $N$, we also obtain universality of the joint distribution of the smallest $k$ singular values for fixed $k$. This is the content of the following theorem. While we do not state the universal distribution explicitly, an exact expression can be derived, and we refer the reader to Section 6 of \cite{leastsv} for details.

\bet\label{t:k} Fix a positive integer $k$. Let $(M_N)^\infty_{N=1}$ be a sparse ensemble and $(G_N)_{N=1}^\infty$ be an ensemble with independent entries of distribution $\mathcal N(0, N^{-1})$. For any matrix $A$, define
\beq \Lambda_k(A) = (N \lambda_1(A), \dots , N \lambda_k(A) ), \eeq
and for any choice of energies $\widehat{E} = (E_i)\in \mathbb R^k$, define 
\beq \Omega{(\widehat {E} )} = \left\{  x\in \R^k \colon x_i \le E_i \text{ for all } i \le k \right\}.\eeq
Let $\widehat E \pm N^{-c}$ denote the vector $(E_i \pm N^{-c}) \in \R^k$. Then
$$ \P\left( \Lambda_k(G_N) \in \Omega(\widehat E - N^{-c} ) \right)- N^{-c}  \le \P\left(\Lambda_k( M_N) \in \Omega(\widehat E)\right)\le \P\left(\Lambda_k(G_N)\in \Omega( \widehat E  + N^{-c}) \right)+ N^{-c},$$
uniformly in all choices of $\widehat{E} \in \mathbb R^k$, for some $c>0$ and large enough $N$.

\eet

Finally, we mention that analogues of our results hold when the matrix entries are complex valued. The proofs are essentially identical. 

\subsection{Proof of Theorem \ref{thm:mainresult}}

Here we present the proof of our main result, Theorem \ref{thm:mainresult}. The proof of Theorem \ref{t:k} is analogous.

We suppose for the remainder of this paper that the entry distributions for all random matrices considered are absolutely continuous with respect to Lebesgue measure, so that their singular values are distinct almost surely and can be strictly ordered. The case of a general matrix $H$ is dealt with by considering $H(\eps)  =  H + \eps V$, where $V$ is a Gaussian matrix. The entires of $H(\eps)$ have distributions that are absolutely continuous for $\eps>0$, since their laws are convolutions with a Gaussian distribution. All results may be extended to $H$ by taking the limit as $\eps$ goes to $0$ and using Weyl's inequality for singular values.

 We recall the following definitions used in the proof.

\bed We say that an event $\mathcal F$ holds with overwhelming probability if for any $D>0$ we have $\P(\mathcal F^c) \le N^{-D}$ for large enough $N$. For a family of events $\{\mathcal F(u) \}$, we say $\{\mathcal F(u) \}$ holds with overwhelming probability if $\sup_u \P(\mathcal F(u)^c)\le N^{-D}$ for large enough $N$.
\eed

\bed\label{def:counting} Given a random matrix with eigenvalues $\{\lambda_i\}$ and $E_1\le E_2$, define the eigenvalue counting function \beq{n(E_1, E_2)=|\{ E_1 < \lambda_i < E_2 \}|}.\eeq
\eed

	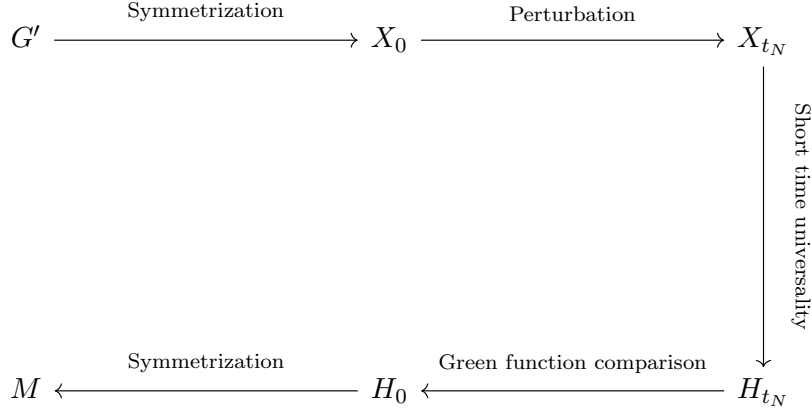
\begin{figure}[h]
	\vspace{0.2cm}
	\begin{center}
	\begin{tikzcd}[row sep = 4cm, column sep = 4cm]

	{ G' } \arrow{r }[yshift=1ex]{\text{Symmetrization}}&	{X_0} \arrow{r}[yshift=1ex] {\text{Perturbation} }& {X_{t_N}} \arrow{d}[anchor=center,rotate=-90,yshift=3ex]{\text{Short time universality}}   \\ 
		
		{M} & {H_0}  \arrow{l}[yshift=3.5ex] {\text{Symmetrization} }& { H_{t_N} } \arrow{l}[yshift=3.5ex]{ \text{Green function comparison} }
			
	\end{tikzcd} 
	\caption{
       		Diagram of the proof of Theorem \ref{thm:mainresult}. The sparse ensemble $M$ and a Gaussian ensemble $G'$ are both symmetrized in order to access their singular values through the eigenvalues of the symmetric matrices $H_0$ and $X_0$. The distribution of the eigenvalue counting function near zero for the perturbed Gaussian matrix $X_t$ is transferred to $H_t$ via a short time universality result. Finally, this distribution is pulled back to $M$ using a Green function comparison argument. 
    } \label{implications}
	\end{center}
	\end{figure}

\bp For any $N$, let $G_N$ be a $N\times N$ matrix with independent entries of distribution $\mathcal N(0, N^{-1})$. We set $M_t = M_N + \sqrt{t}G_N$ and let $\lambda^\circ_1(t)$ be the least singular value of $M_t$, suppressing dependence on $N$. Let $G'_N$ and $G''_N$ be matrices independent from $G_N$ with the same entry distributions, and let $\tilde \lambda^\circ_1(t)$ be the least singular value of $G_t = G'_N + \sqrt{t} G''_N$.
Define the block matrices 
\beq\label{eqn:ht} H_t = \bma  0  & M_t \\ M_t^\trans & 0 \ema  ,\quad X_t = \bma  0  & G_t \\ G_t^\trans & 0 \ema.\eeq
Note that the eigenvalues of $H_t$ are precisely the singular values of $M_t$ and their negatives, and same relation holds between $X_t$ and $G_t$.

For any $\eps>0$, Lemma \ref{lem:shorttime} on short time universality shows that there exists a coupling $(\lambda_i(t), \tilde \lambda_i(t) )$ of the $\lambda^\circ_i(t)$ and $\tilde \lambda^\circ_i(t)$ such that with overwhelming probability we have
\beq\label{eqn:st1} |\lambda_1(t_N) - \tilde \lambda_1(t_N) | \le \frac{1}{N^{1+\sigma}} \eeq
for some $N$-dependent parameter $t_N$ satisfying $t_N \le N^{-1+\eps}$ and some $\sigma >0$. 

We conclude by removing the time evolution. This is accomplished by comparing smoothed versions of the eigenvalue counting functions using results in Section \ref{sec:removal}. Let $n_t$ be the eigenvalue counting function for the $\lambda_i(t)$ and $\tilde n_t$ be the counting function for the $\tilde\lambda_i(t)$, as in Definition \ref{def:counting}. Set \beq E=r/N, \quad y=N^{-1-\sigma},\quad \eta = y N^{-32\sigma}.\eeq Recasting (\ref{eqn:st1}) in terms of these counting functions gives that, for any $D>0$,
 \beq \P(\tilde n_{t_N}(-E- y,E +y)=0)-N^{-D} \le \P(n_{t_N}(-E,E)=0) \le  \P(\tilde n_{t_N}(-E+y ,E-y)=0) +N^{-D}.\eeq
The conclusion of Lemma \ref{lem:qbound} is that 
\beq \E q(\tr \chi_{E+y}\star \theta_{\eta}(H_{t_N}))  - N^{-D} \le\P(n_{t_N}(-E,E)=0)\le \E q(\tr \chi_{E-y}\star \theta_{\eta}(H_{t_N})) + N^{-D},\eeq
where $q$ is a smooth cutoff function defined in Section \ref{sec:removal}, and 
\beq\theta_\eta = \frac{\eta}{\pi(x^2 + \eta^2)} = \frac{1}{\pi} \Im \frac{1}{x - i\eta}.\eeq

Further,  Lemma \ref{lem:gcompare} shows that there exists $\eps>0$ and $c>0$ such that
\beq \E q(\tr \chi_{E-y}\star \theta_{\eta}(H_t)) \le \E q(\tr \chi_{E-y}\star \theta_{\eta}(H_0)) + N^{-c}\le  \P(n_0(-E+2y,E-2y)=0)+ N^{-c}\eeq
for all $t\le N^\eps/N$, where the last inequality holds by another application of Lemma \ref{lem:qbound}. 

Fix this $\eps$ and corresponding $t_N \le N^\eps/N$ for the rest of the proof. After adjusting $c$ downward and combining this display with the previous one, we have for large enough $N$ that \beq \P(n_{t_N}(-E,E)=0) \le  \P(n_0(-E+2y,E-2y)=0)+ N^{-c}.\eeq
Similar reasoning, interchanging the roles of $n_0$ and $n_t$, gives 
\beq \P(n_0(-E-2y,E+2y)=0) \le  \P(n_{t_N}(-E,E)=0)+ N^{-c}.\eeq
We conclude, after setting $\tilde y =2y$, that there exists $c>0$ such that, for large enough $N$,
\beq\P(\tilde n_0(-E- \tilde y,E+ \tilde y)=0) - N^{-c} \le \P(n_0(-E,E)=0) \le \P(\tilde n_0(-E+ \tilde y,E- \tilde y)=0) + N^{-c}.\eeq
Rephrased in the language of cumulative distribution functions, this is 
\beq\label{eqn:pbounds}  \P\left(\tilde \lambda_1(0)\le E -  \tilde y \right) - N^{-c} \le \P\left( \lambda_1(0)\le E \right) \le  \P\left(\tilde \lambda_1(0)\le E +  \tilde y \right) + N^{-c}.\eeq

Since $\tilde \lambda_1$ is the least singular value of a matrix with i.i.d. entries, we may use Theorem 1.3 in \cite{leastsv} to control its distribution. After changing the variable of integration, we have for some small $c>0$ that
\beq\label{eqn:gexplicit} \P\left(\tilde \lambda_1(0)\le \frac{r}{N} \right) =  1 - e^{-r^2/2 -r} + O(N^{-c}). \eeq
Note that (\ref{eqn:gexplicit}) shows 
\beq \P\left(\tilde \lambda_1(0)\in (E- \tilde y, E+ \tilde y) \right) \le N^{-c}.\eeq
Together with (\ref{eqn:pbounds}), this completes the proof.
\ep

\section{Short time universality}\label{sec:deterministic}

In order to state the main result of this section, we introduce the following notation. We fix $\delta_1>0$ and let $g$ and $G$ be $N$-dependent parameters such that 

\beq N^{-1+\delta_1} \le g \le N^{-\delta_1},\quad G\le N^{-\delta_1}.\eeq

We consider a deterministic matrix $V$ of initial data and let $B_t=\{B_{ij}(t)\}_{1\leq i, j \leq N}$ be a matrix of i.i.d. real Brownian motions. We define
\beq M_t = V + \frac{1}{\sqrt N } B_t,\quad H_t =  \bma  0 & M_t \\ M_t^\trans & 0 \ema . \eeq 
Let $\{s_i(t) \}_{i = -N}^N$ (omitting the zero index) be the eigenvalues of $H_t$, which are the singular values of $M_t$ along with their negatives. We set
\beq m_V(z) = \sum_{i= - N}^N \frac{1}{s_i(0) - z}.\eeq 

\bed\label{def:regular}
With $g$ and $G$ defined as above, we say $V$ is $(g,G)$-regular if 
\beq c \le \Im m_V(E+i\eta) \le C\eeq
for $|E| \le G$ and $\eta \in [ g , 10]$ for large enough $N$, and if there exists a constant $C_V$ such that $|v_i| \le N^{C_V}$ for all $v_i$. \eed

We now state the main result of this section. Let $W$ be a random matrix whose entries are i.i.d. $\mathcal N(0, N^{-1})$ variables and let $\tilde B_t=\{\tilde B_{ij}(t)\}_{1\leq i, j\leq N}$ be a matrix of i.i.d. real Brownian motions. Define $W_t = W + N^{-1/2} \tilde{B}_t$.  Recall  $\{s_i(t)\}_{i=1}^N$ are the singular values of $M_t$, and let $\{r_i(t)\}_{i=1}^N$ be the singular values of $W_t$. %The choice of $\tilde{B}_t$ will be addressed in the next subsection.

\bet\label{thm:homomain} Fix $\sigma>0$, and let $V$ be a $(g,G)$-regular deterministic matrix. Let $M_t$, $W_t$, $\{s_i(t)\}$, and $\{r_i(t)\}$ be defined as above. Then there exists a coupling of the processes $\{s_i(t)\}$ and $\{r_i(t)\}$ such that the following holds. Given parameters $0< \omega_1 < \omega_0$ and times $t_0 = N^{-1+\omega_0}$, $t_1 = N^{-1+\omega_1}$, with the restrictions that 
\beq gN^{\sigma} \le t_0 \le N^{-\sigma} G^2,\quad 2\omega_1 < \omega_0,\eeq
there exist $\omega, \delta>0$ such that
\beq |s_i(t_a)) - r_i(t_a) | < N^{-1-\delta}\eeq
for $i < N^\omega$ and $t_a=t_0+t_1$.
\eet

\subsection{Preliminaries}

We prove Theorem \ref{thm:homomain} by analyzing the SDE that governs the time evolution of the singular values of $M_t$.
Given initial data $(s_i(0))_{i=1}^N$, the SDE for the perturbed singular values is
\beq\label{eqn:svsde} d s_i(t) = \frac{dB_i}{\sqrt{N}}  + \frac{1}{2N}\sum_{ j\neq i}\left( \frac{1}{s_i(t) - s_j(t)} + \frac{1}{s_i(t) + s_j(t)}  \right) dt. \eeq
To be precise, we mean that the stochastic process $(s_i(t))_{i=1}^N$ and the stochastic process given by the singular values of $M_t$ are equal in distribution.

\begin{comment}
Above, $dB_i = \sum_{j,k}T^\trans_{ij} dB_{jk}(t)R^\trans_{ki} $, where $T$ and $R$ are the left and right orthogonal matrices in the singular decomposition $V = T D  R$ with $D$ diagonal. Similarly, the SDE for $r_i(t)$ is nearly identical to \eqref{eqn:svsde}, with the only difference being that $dB_i$ is replaced by its counterpart, $d\tilde{B}_i = \sum_{j,k} \tilde T^\trans_{ij} d \tilde B_{jk}(t) \tilde R^\trans_{ki} $. Here $\tilde T$ and $\tilde R$ represent the (random) left and right orthogonal matrices in the singular value decomposition $W = \tilde T \tilde D \tilde R$.  

For the coupling, we pick the Brownian motion matrix $B_t$ arbitrarily and set 
\begin{equation}
\tilde B_t = \tilde T T^\trans B_t R^\trans \tilde R.
\end{equation}
This ensures that the SDEs are coupled so that $d B_i = d\tilde{B}_i$ for all $i$.
\end{comment}

Technical information about this SDE is contained in Appendix \ref{a:sde}, where we show existence and uniqueness of strong solutions and verify that it represents the claimed evolution of the singular values (in distribution). The arguments used in the appendix may also be used to show existence and uniqueness of strong solutions for the other SDEs considered in this section.

%%% Our proof of Theorem \ref{thm:homomain} contains two main ideas. First, we note that the study of the singular values of a $N\times N$ non-symmetric matrix is equivalent to the study of the eigenvalues of a $2N\times 2N$ symmetric matrix, as mentioned in the line following \eqref{eqn:ht}. Proving rigidity for small singular values after a short time then reduces to proving eigenvalue rigidity in the bulk around $E=0$ for the symmetric matrix model. This is handled by adapting the arguments of \cite{LY}. Rigidity of the singular values is crucial for our analysis of \eqref{eqn:svsde}, and we state the necessary result as Lemma \ref{lem:uniformlocal}. 

Our main idea is to analyze a symmetrized version of (\ref{eqn:svsde}). Define $s_{-k}(t)$ and $B_{-k}(t)$ for $1\le k \le N$ by
\[
s_k(t) = -s_{-k}(t),\quad B_{k}(t)=-B_{-k}(t),\quad t\geq 0.
\]
Then the SDE (\ref{eqn:svsde}) is equivalent to
\beq\label{eqn:sym}
\, d s_i = \frac{dB_i}{\sqrt{N}}+ \frac{1}{2N}\sum_{\substack{-N \le j \le N \\ j\neq \pm i, 0}} \frac{\, d t}{s_i - s_j }.
\eeq

Note we have labeled the particles from $-1$ to $-N$ and $1$ to $N$, so that the zero index is omitted. Unless otherwise stated, this will be our convention throughout this section, and we will no longer note this omission explicitly. 

Note also that (\ref{eqn:sym}) is the same as a standard DBM, except it is missing the repulsion term for $j=-i$. In particular, for the least singular value $i=1$, the force deflecting this particle from the origin is weaker than in the usual DBM. This complicates the analysis of \eqref{eqn:sym}, since we are now missing an important regularizing effect, and level repulsion estimates such as the one used in \cite{LY} to study the short time behavior of DBM seem out of reach. The key point of this section is that the smallest singular values reach equilibrium on short time scales even without this regularizing term.

We begin by studying the differences $s_i(t)-r_i(t)$. By demanding that the SDEs for the $r_i$ and $s_i$ are driven by the same Brownian motions, we force the Brownian motion term to vanish in the equation for the differences, which renders it considerably easier to work with. This coupling approach was developed in \cite{wigfixed, LY}. Here we use a refinement of this idea, introduced in \cite{fixed}, and construct a continuous interpolation. Next, we obtain precise control over this difference equation by comparing its semigroup elements to the solution kernel of an appropriate integral equation, which represents the homogenized version of the DBM. The solution kernel is approximately the Poisson kernel, which corresponds to the fact that the DBM difference equation is a discrete version of the square root of the Laplacian (see Section 7.11 of \cite{lieb2001analysis} for details in the continuous case). Theorem \ref{thm:homomain} then follows from a short calculation in Subsection \ref{sec:conclusion} exploiting cancellation in the kernel coming from the symmetrization of the SDE. 

Except for the final calculation, this plan follows closely \cite{fixed}, whose methods are fundamental to our approach. The work \cite{fixed} proved short time universality, at level of particle gaps, for DBM with deterministic initial data. This result, called ``fixed energy universality,'' is of great importance because it enables the proof of bulk universality for a vast array of matrix models. Our result is similarly flexible and we have stated it in more generality than needed for this work in order to facilitate future applications. \\

\remark In this section we deal only with the real variable case, where we apply a real Gaussian perturbation. If we consider complex matrices, the appropriate SDE is
	\beq\label{eqn:cmplx}
	\,d s_k = \sqrt{\frac{1}{2N}}\, d B_k + \frac{1}{2N}\sum_{ j\neq \pm k} \frac{\,d t}{s_k - s_j }.
	\eeq
The only difference here is that the diffusion term loses a factor of $\sqrt{2}$. The argument in the complex case is the same as the real case, and the main result of this section Theorem \ref{thm:homomain} holds with no difference.

\subsubsection{Definitions}

Given a probability measure $\rho(E)\, dE$ and some even number $2N$, we define the classical particle locations ${\gamma}_i$ as follows, suppressing dependence on $N$. For $i\ge 1$, we set 
\beq \gamma_i = \inf\left\{ x \colon \int_{-\infty}^x \rho(E)\, dE \ge \frac{N+i-1}{2N}  \right\}.\eeq
For $i\le -1$, we set 
\beq \gamma_i = \inf\left\{ x \colon \int_{-\infty}^x \rho(E)\, dE \ge \frac{N+i}{2N}  \right\}.\eeq
In accordance with our labeling convention, we do not define ${\gamma}_0$. Our labeling is chosen so that $\gamma_1=0$ for a symmetric distribution centered at $0$. 

We recall the Stieltjes transform of a $N\times N$ matrix $M$ with eigenvalues $\lambda_i$ is given by
\beq m_N(z) = \frac{1}{N} \sum_i \frac{1}{\lambda_i - z}.\eeq

\subsection{Interpolation}

Our goal is to compare solutions to the SDE \eqref{eqn:sym} with initial data $s_i(0)$, a $(g,G)$-regular deterministic set of singular values, and $r_i(0)$, the singular values of the Gaussian ensemble $G$. We introduce times $t_0=N^{-1+\omega_0}, t_1=N^{-1+\omega_1}$ with parameters $0< \omega_1 < \omega_0$ and such that $gN^\sigma \le t_0 \le N^{-\sigma}G^2$ and $2\omega_1 \le \omega_0$. 

For $\alpha \in [0,1]$ we introduce a continuous interpolation between $s_i(t)$ and $r_i(t)$:
\beq\label{eqn:interp}
\, d z_i(t,\alpha) = \frac{dB_i}{\sqrt{N}} + \frac{1}{2N}\sum_{j\neq \pm i, 0} \frac{\, d t}{z_i(t,\alpha) - z_j(t,\alpha) }.
\eeq
The initial values are given by
\beq z_i(0,\alpha)=(1-\alpha)r_i(t_0) + \alpha s_i(t_0),\eeq
and the processes $z(t,0)$ and $z(t,1)$ are shown below to provide the desired coupling of the singular values of $M_t$ and $W_t$.

More precisely, we construct solutions of \eqref{eqn:interp} for a countable dense set of $\alpha\in[0,1]$ and use the continuity of the map from initial data to solution paths to construct solutions for the remaining $\alpha$. This approach avoids the possibility that uncountably many sets of measure zero accumulate. See Appendix \ref{sec:interpolationdetails} for details.

In particular, $z_i(t,0)=r_i(t_0+t)$ and $z_i(t,1)=s_i(t_0 +t)$. The time shift is enforced to simplify the notation. Using the deformed law Theorem \ref{thm:deformedlaw}, we will establish rigidity of the particles after the short time $t_0$. Hence, the time shift means the $z_i(t)$ are rigid for times $t\ge 0$, which is notationally convenient. See Lemma \ref{lem:uniformlocal} for the precise statement.

We define the $\gamma_i(t)$ to be the classical locations associated to the eigenvalue density of the symmetrization of the deformed matrix $V + \sqrt{t}W$, which has eigenvalues $\{s_i(t)\}_{i=-N}^N$. We let $\gsc_i$ denote the classical locations for the semicircle law. 
\subsubsection{Interpolating Measures}

We would like to introduce a family of measures $\rho(E,t,\alpha)\, dE$ which will interpolate between the densities of the initial data. One approach is to take the free convolution of $z_i(0,\alpha)$ with the semicircle law. However, the resulting measures are not regular enough for our purposes. The following two lemmas assert that we can construct interpolating measures with better regularity properties. The first concerns the regularity of the $\rho(E,t,\alpha)$, and the second verifies these measures interpolate between the densities of the two initial ensembles. Here $m(z,t,\alpha)$ is the Stieltjes transform of $\rho(E,t,\alpha)\, dE$, and we let $\gamma_i(t,\alpha)$ denote the classical locations for $\rho(E,T,\alpha)\, dE$.

The set  $\mathcal G_\alpha$ used in the lemma is defined as follows. We fix $q^*  \in (0,1)$ and set $k_0$ to be the largest index such that
\beq |\gamma_{k_0}(t_0)| \le q^* G, \quad |\gamma_{-k_0}(t_0)| \le q^* G, \quad  |\gsc_{k_0}|  = |\gsc_{- k_0}| \le q^* G.\eeq
Then  $\mathcal G_\alpha$ is defined as 
\beq \mathcal G_\alpha = [\alpha \gamma_{-k_0}(t_0) + (1-\alpha)\gsc_{-k_0}, \alpha \gamma_{k_0}(t_0) + (1-\alpha)\gsc_{k_0}].\eeq

\bel A family of interpolating measures $\rho(E,t,\alpha)\, dE$ exists (in a sense made precise by the following lemma) such that the following holds. Let $\delta >0$. For $|E| \le N^{-\delta} t_0$, $t\le N^{-\delta}t_0$, and $N^{-1+\delta} \le \eta \le 10$ all of the following estimates are true with overwhelming probability.
\beq |\partial_E \rho(E,t,\alpha) | \le \frac{C}{t_0},\quad \rho(0,0,\alpha) = \rho(0,0,0) = \sc(0)\eeq
\beq |\rho(E,t,\alpha) - \sc(0)| \le C\left( \frac{t\log(N)}{t_0} + \frac{|E|}{t_0}  \right)\eeq
\beq |\rho(E,t,\alpha) - \rho(0,t,\alpha) - (\rho(E,t,0) - \rho(0,t,0))| \le C\left( \frac{t\log(N)}{t_0} + \frac{|E|}{t_0}  \right)\eeq
Further, for $q\in(0,1)$ and $E\in q \mathcal G_\alpha$, $N^{-1+\delta}\le \eta \le 10$, with overwhelming probability
\beq\label{eqn:stbounds} |m(z,t,\alpha)| \le C \log(N),\quad c \le \Im m(z,t,\alpha) \le C.\eeq
For $q\in (0,1), 0 \le t \le t_1$, and $E\in q\mathcal G_\alpha$, with overwhelming probability
\beq|\partial_z m(z,t,\alpha)| \le \frac{C}{t_0 + \eta}.\eeq
\eel
\bp
We indicate where such claims are shown in \cite{fixed}, where the proofs are identical. The estimates on $\rho(E, t , \alpha)$ are in Lemma A.6. The first bound in (\ref{eqn:stbounds}) comes from the construction of $\rho(E,t,\alpha)$ and the reasoning in (7.12) of \cite{LY}, and the second is Lemma A.4. The last claim is Lemma A.5. 
\ep

Define $d(i,j) = |i-j|$ if $ij>0$ and $d(i,j)=|i-j| -1$ if $ij <0$. This is just to define an appropriate distance for our indexing, since no element is indexed with $0$. The proof of the following lemma is the same as Lemma 3.4 in \cite{fixed}.

\bel\label{lem:gammarigid} The following estimates hold for the $\rho(E,t,\alpha)\, dE$ constructed in the previous lemma with overwhelming probability. We have for $\eps >0$ and $\omega_1 < \omega_0/2,$ 
\beq \sup_{0\le t \le 10t_1} |\gamma_1(t,1) - \gamma_1(t_0 +t)| \le N^{-1 -\omega_0/2 + \omega_1 + \eps},\eeq
\beq \sup_{0\le t \le 10t_1} | \gamma_1(t,0) - 0 | \le  N^{-1 -\omega_0/2 + \omega_1 + \eps}.\eeq
For $|j|, |k| \le N^{\omega_0/2}$, and any choice of $t\le 10 t_1$, $\omega_1 \le \omega_0/2$, $\alpha \in[0,1]$,
\beq \gamma_k(t,\alpha) - \gamma_j(t,\alpha) = \frac{d(k,j)}{\sc(0)} + O(N^{-1})   .\eeq
\eel

We further have an analogue of Lemma 3.5 in \cite{fixed}, giving rigidity and a local law. Here we choose a centered interval $\hat C_q$ of indices, asymptotically of size $qGN$, so particles with indices in $\hat C_q$ can be controlled uniformly in $\alpha$. Precisely, we define $k_1$ to be the largest integer such that 
\beq \bigcup_{0\le \alpha \le 1}  [\alpha \gamma_{-k_1}(t_0) + (1-\alpha)\gsc_{-k_1}, \alpha \gamma_{k_1}(t_0) + (1-\alpha)\gsc_{k_1}] \subset \bigcap_{0\le \alpha\le 1} \mathcal G_\alpha \cap  \{ - \mathcal G_\alpha \},\eeq
and set for $q\in (0,1)$
\beq \hat {\mathcal C}_q = \{ j \colon |j| \le q k_1 \}.\eeq
We let $m_N(z,t,\alpha)$ be the Stieltjes transform of the $z_i(t,\alpha)$.
\bel\label{lem:uniformlocal} Fix $\eps, \delta, \delta_1, D >0$ and $q\in (0,1)$. Then the following two estimates hold.
\beq \P\left[  \sup_{0\le t \le N^{-\delta_1} t_0} \sup_{i\in \hat C_q} \sup_{0\le \alpha \le 1} |z_i(t,\alpha) - \gamma_i(t,\alpha)| \ge \frac{N^\eps}{N} \right] \le N^{-D} \eeq
\beq \P\left[  \sup_{N^{-1+\delta} \le \eta \le 10} \sup_{0\le t \le N^{-\delta_1} t_0} \sup_{0\le \alpha \le 1} \sup_{E\in q \mathcal G_\alpha}    | m_N(z,t,\alpha) - m(z,t,\alpha) | \ge \frac{N^\eps}{N\eta} \right] \le N^{-D} \eeq
\eel
The key input to the proof of this lemma is a deformed local law for small times $t$. In the eigenvalue context this law was shown in \cite{LY}. We prove the needed singular value version of this result in Section \ref{sec:dqclaw}. Excepting this change, the proof is identical to the one described in Appendices A and B of \cite{fixed}.

\subsection{Short range equation}

Define a centered process
\beq \tilde z_i(t,\alpha) = z_i(t,\alpha) - \gamma_1(t,\alpha).\eeq
Note this differs from (3.35) of \cite{fixed} because we use the classical location $\gamma_1$ to center the particles, as we have no $\gamma_0$. The new classical locations are \beq\tilde \gamma_i(t,\alpha) = \gamma_i(t,\alpha) - \gamma_1(t,\alpha).\eeq 
We recall from \cite{schnelli2} that 
\beq \partial_t \gamma_i(t,\alpha) = - \Re[ m(\gamma_i(t,\alpha), t ,\alpha)].\eeq
Then the SDE that governs the $\tilde z_i$ is
\beq d \tilde z_i(t,\alpha) = \frac{dB_i}{\sqrt{N}} + \left( \frac{1}{2N} \sum_{j\neq \pm i}  \frac{1}{\tilde z_i(t,\alpha) - \tilde z_j(t,\alpha)} + \Re[m(\gamma_1(t,\alpha), t, \alpha) ]\right)\, dt .\eeq

Because we do not have good control over the extremal particles, and because we are interested only in particles near the origin, it is convenient to introduce a short range cutoff. We fix $q_*\in(0,1)$, and $\omega_l,\omega_A >0$ such that 
\beq 0 <  \omega_1 < \omega_\l < \omega_A < \omega_0/2.\eeq
We choose the parameters in this way so the error term in the forthcoming Lemma \ref{lem:hatapprox} is $o(1/N)$. Given $q\in(0,1)$, we define 
\beq A_q = \{ (i,j) \colon |i - j| \le N^{\omega_l} \text{ or } ij>0, i\notin \hat C_q, j\notin \hat C_q \},\eeq
and we let ${A_{{q_*},(i)}}$ be the indices $j$ such that $(i,j)\in A_{q_*}$ and $A^c_{{q_*}, (i)}$ be the indices $j$ such that $(i,j) \notin  A_{q_*}$.

Define $\hz_i(t,\alpha)$ as the solution to, for $|i| \le N^{\omega_A}$,
\beq d\hz_i(t,\alpha) =\frac{dB_i}{\sqrt{N}} + \frac{1}{2N} \sum_{j}^{A_{{q_*},(i)}} \frac{1}{\hz_i(t,\alpha) - \hz_j(t,\alpha)}\, dt.\eeq
For $|i| > N^{\omega_A}$,
\beq d\hz_i(t,\alpha) = \frac{dB_i}{\sqrt{N}}+ \frac{1}{2N} \sum_{j}^{A_{{q_*},(i)}} \frac{1}{\hz_i(t,\alpha) - \hz_j(t,\alpha)} dt +  \frac{1}{2N} \sum_j^{A^c_{q_*}} \frac{dt }{\tilde z_i - \tilde z_j}  + \Re[ m(\gamma_1(t,\alpha), t ,\alpha)] \, dt.\eeq
Here the initial condition is $\hz_i(0)=\tilde z_i(0)$. 

The $\hz_i$ are good approximations to the $\tilde z_i$. This is the content of the following lemma, whose proof is identical to the proof of Lemma 3.7 in \cite{fixed}.
\bel\label{lem:hatapprox} Fix $\eps, D>0$. Then, for large $N$,
\beq \P\left[  \sup_{0\le t \le t_1 } \sup_{i} \sup_{0\le \alpha \le 1} |\hz_i(t,\alpha) - \tilde z_i(t,\alpha) |  \ge N^{\eps}t_1\left( \frac{N^{\omega_A}}{N^{\omega_0}}  + \frac{1}{N^{\omega_\l}} + \frac{1}{\sqrt{NG}}\right) \right] \le N^{-D}.\eeq
\eel

\subsubsection{Short Range Kernel}

Here we introduce a coupled parabolic equation with no Brownian motion term. We define $u_i = \partial_\alpha \hz_i(t,\alpha)$, where the differentiation is justified in Appendix \ref{sec:interpolationdetails}.
% By fixing $\alpha_1, \alpha_2$, and using the maximum principle for the equation for $\hz_i(t,\alpha_1) - \hz_i(t,\alpha_2)$, we see as in (3.66) of \cite{fixed} that $\hz_i$ is Lipschitz in $\alpha$, so the derivative exists almost everywhere and the Fundamental Theorem of Calculus holds. 

We have 
\beq \partial_t u_i = \sum_j^{A_{q_*, (i)}} B_{ij}(u_j - u_i) +\xi_i = - (\B u)_i + \xi_i,\eeq
where \beq B_{ij} = \frac{\one_{i\neq \pm j}}{2N(\hz_i - \hz_j)^2},\eeq
the operator $\B$ is implicitly defined by the above equation, and $\xi_i$ is an error term that vanishes for $|i| \le N^{\omega_A}$. Because it vanishes for small values of $i$, we will show its effect on particles near the origin is negligible.

We let $\mathcal U$ be the semigroup associated to $\mathcal B$. Its elements $\mathcal U_{ij}(s,t)$ are defined so that, if $v(t)$ is any solution of the system $\partial_t v_i  = - (\B v)_i$, then for any times $t,s \ge 0,$
\beq  v_i(t) = \sum_{j = - N}^N \U_{ij}(s,t) v_j(s).\eeq
Note that the $\U_{ij}$ are random.

\subsubsection{Finite speed estimates}
We state two results on the decay of the semigroup elements $\U_{ij}$. Their proofs are straightforward adaptations of those given in Section 4 of \cite{fixed} for the kernel $\U^{(B)}$ in that reference.

\bel\label{lem:finprop}
Let $0\le s \le t \le t_1$. Fix $0< q_1 < q_2 < q_*$ and $D,\eps >0$. For every $\alpha$ there exists an event $\Fa$ such that $\P(\Fa)\ge 1 - N^{-D}$ and on which the following estimates hold. If $i\in \hc_{q_2}$ and $0\le s \le t \le 10t_1$, then
\beq | \U_{ji}(s,t,\alpha)| \le N^{-D},\quad |i-j| > N^{\omega_\l +\eps}.\eeq
If $i\notin \hc_{q_2}$, $j\in \hc_{q_1}$, and $0\le s \le t \le 10t_1$, then
\beq | \U_{ji}(s,t,\alpha)| \le N^{-D}.\eeq
\eel

\bel\label{lem:poisson}
Let $0< q_1 < q_*$ and $D, \eps >0$. For every $\alpha$ there exists an event $\Fa$ such that $\P(\Fa)\ge 1 - N^{-D}$ and on which he following bound holds. For $i, j\in \hc_{q_1}$ and $0\le s \le t \le 10t_1$,
\beq \U_{ij}(s,t) \le \frac{N^\eps}{N} \frac{(t-s) \vee N^{-1}}{\left( \frac{i-j}{N}\right)^2 + \left( (t-s) \vee N^{-1} \right)^2 }.\eeq
\eel

\subsection{\emph{A priori} estimates}\label{sec:homog}

We first make some definitions necessary for the homogenization argument. We fix a constant $\eps_B>0$ such that $\omega_A - \eps_B > \omega_l$, and fix an integer $a$ such that $0 < |a| \le N^{\omega_A  - \eps_B}$. We also define the deterministic particle locations  $\flat_j  = j(2N\sc(0))^{-1}$.

We consider a particular solution $w$ of the short range equation that will play a key role:
\beq \partial_t w_i = - (\B w)_i, \quad w_i(0) = 2N\delta_a(i).\eeq
Define the cutoff $\eta_l = N^{\omega_\l} (2N\sc(0))^{-1}$. We will need the kernel $p_t(x,y)$ of the following equation
\beq \partial_t f (x) = \int_{|x-y| \le \eta_\l} \frac{f(y) -f(x)}{(x-y)^2}\sc(0) \, dy.\eeq
The main result of this subsection is that this kernel is a good approximation to elements of the short range semigroup.
We choose parameters $s_0$ and $s_1$ such that
\beq N^{-1} \ll s_0 \ll s_1 \ll t_1 \ll t_0\eeq
and introduce
\beq f(x,t) =\frac{1}{2N} \sum_{j\neq 0 } p_{s_0 +t -s_1}(x,\flat_j)w_j(s_1), \quad f_i(t) = f(\hz_i(t,\alpha), t).\eeq
We now collect some bounds on these objects. We let $p^{(k)}_t(x,y)$ denote the derivative in $x$. The following lemma is Lemma 3.11 in \cite{fixed}.

\bel\label{lem:pbound} Fix $\eps_1, \eps_2, D_1>0$. For $t$ such that $N^{-D_1} \le t \le N^{-\eps_1}\eta_\l$ and any $D_2>0$,
\beq p_t(x,y) \le C \frac{t}{(x-y)^2 + t^2},\quad p^{(k)}_t(x,y) \le \frac{C}{t^k} p_t(x,y) + N^{-D_2},\quad |\partial_t p_t(x,y)| \le \frac{C}{x^2 + y^2} + N^{-D_2}.\eeq
If additionally $|x-y| > N^{\eps_2}\eta_l$,
\beq p_t(x,y) \le N^{-D_2},\quad p^{(k)}_t(x,y) \le N^{-D_2}.\eeq
\eel

\bel\label{lem:wbound} Fix $\eps>0$ and $t>0$. The following estimates hold with overwhelming probability.
\beq \| w(t) \|_1 \le 2N,\quad \|w(t)\|_\infty \le N^\eps t^{-1},\quad \|w(t)\|^2_2 \le N^\eps t^{-1}.\eeq
Additionally, if $|i|\ge |a| + N^{\omega_l + \eps_1}$ then for any $D>0$ we have with overwhelming probability
\beq |w_i(t)| \le N^{-D}.\eeq
 \eel
 
 \bp We see from the definition of $\B$ and examining $\min w_i$ that all $w_i(t)$ are non-negative for $t\ge 0$. Hence 
 \beq \| w(t) \|_1 = \|w(0)\|_1 = 2N.\eeq
 The second estimate is a consequence of Lemma \ref{lem:poisson} and the third estimate is (3.109) in \cite{fixed}. The last estimate follows from applying Lemma \ref{lem:finprop}.
 \ep
 
\bel\label{lem:fbound} We have 
\beq f(t,\hz_i)\le \frac{C}{t+s_1},\eeq
and for any $\eps_1>0$, if $|i| \ge |a| + N^{\omega_\l + \eps_1}$,
\beq |f_i| \le N^{-D}.\eeq
For $k\ge 1$,
\beq f^{(k)}(t,\hz_i) \le \frac{C}{t-s_1 + s_0} f(t,\hz_i) + N^{-D}.\eeq
\eel
\bp See the beginning of the proof of Lemma 3.13 in \cite{fixed}.\ep
\subsection{Semigroup estimate}

In the lemma below, we temporarily break with our convention of omitting the zero index in order to match the notation of \cite{fixed}. The strange indexing is due to the fact that we will need to apply it in the context of an even number of particles labeled $-N$ to $-1$ and $1$ to $N$. We think of shifting the positive indices of $w$ and $f$ down $1$, so we will define for $i\ge 0$,
\beq \ w^{\mathrm {new}}_i = w_{i+1}\eeq
and for $i \le -1$,
\beq \ w^{\mathrm {new}}_i = w_{i}\eeq
Then the largest positive index in the following sum will not have a matching negative index, but we will see this makes no difference to our application below. 
\bel\label{lem:initial} For any $\eps_1, \eps_2, D >0$ and $\alpha$, there exists an event $\Fa$ such that $\P(\Fa)\ge 1 - N^{-D}$ and on which the following estimate holds. 
\beq \| w^{\mathrm {new}}(s_1) - f^{\mathrm {new}}(s_1) \|^2_2  \le C s_0 \sum_{|i| \le N^{\omega_A - \eps_B} + N^{\omega_\l + \eps_2}} \sum_{|i-j|\le \l, j \neq  i, -(i+1)} \frac{(w_i^{\mathrm {new}}(s_1) - w_j^{\mathrm {new}}(s_1))^2}{(i-j)^2} \eeq \beq+ \frac{N^{\eps_1}}{s_1}\left( \frac{1}{(N s_0)^2} + \frac{(Ns_0)^2}{\l^2}  \right) \eeq
\eel
\bp We omit the superscript for this proof. Lemma 3.12 in \cite{fixed} gives 
\beq \| w(s_1) - f(s_1) \|^2_2  \le C s_0 \sum_{|i| \le N^{\omega_A - \eps_B} + N^{\omega_\l + \eps_2}} \sum_{|i-j|\le \l, j \neq  i} \frac{(w_i(s_1) - w_j(s_1))^2}{(i-j)^2} + \frac{N^{\eps_1}}{s_1}\left( \frac{1}{(N s_0)^2} + \frac{(Ns_0)^2}{\l^2}  \right) .\eeq
Note the additional terms in the sum where $j = - ( i + 1)$. We separate them out and bound them.
\beq s_0 \sum_{|i| \le N^{\omega_A - \eps_B} + N^{\omega_\l + \eps_2}} \sum_{j = -(i+1)} \frac{(w_i(s_1) - w_j(s_1))^2}{(i-j)^2} \le s_0 \sum_i  (w_i(s_1) - w_{-{i+1}}(s_1))^2 \le s_0 C \| w(s_1)\|^2_2 \le  N^\eps s_0 s_1^{-1}. \eeq
In the last inequality we used Lemma \ref{lem:wbound}, and $\eps>0$ is arbitrary. Now we are done, because
\beq \frac{s_0}{s_1} \ll \frac{1}{s_1}\frac{1}{(Ns_0)^2}.\eeq
\ep

The next lemma is the key estimate for the homogenization. It is here we deal with the missing repulsion term in the symmetrized equation \eqref{eqn:sym}.

\bel For $t\ge s_1$, the It\^{o} differential of $\|w(t) - f(t)\|_2^2$ takes the form
\beq d \frac{1}{N} \sum_i (w_i - f_i)^2 = - \langle (w(t) -f(t))  , \B (w(t) -f(t))\rangle \, dt + X_t\, dt + d M_t,\eeq
where $M_t$ is a martingale and $X_t$ is a process defined by the above equality.

Fix $\eps, D >0$. Then additionally for every $\alpha$ there is an event $\Fa$ such that $\P(F_\alpha)\ge 1  - N^{-D}$ and on which the following bounds hold. For $s_1 \le t \le 9 t_1$,
\beq |X_t| \le \frac{1}{5}\langle (w(t) -f(t))  , \B (w(t) -f(t))\rangle + \frac{C}{t+s_1}\frac{N^\eps}{t-s+s_0}\left( \frac{1}{\sqrt{N(t-s_1 +s_0)}} \right). \eeq
For any $u_1$ and $u_2$ with $9t_1 > u_2 > u_1 \ge s_1$,
\beq \left| \int_{u_1}^{u_2} d M_t   \right| \le \frac{N^\eps}{N} \frac{1}{(u_1 + s_1)^{3/2} (u_1 - s_1 + s_0)^{1/2}}.\eeq
\eel

\bp  This result is the analogue of Lemma 3.13 in \cite{fixed}, and most of the proof goes through with only notational changes. We comment on one important difference. The term (3.137) becomes, using the notation of that reference,
\beq \left( \frac{1}{2N} \sum_{j\neq \pm i}^{A_{q_*}\setminus A_{2, (i)}} \frac{f_j - f_i}{(\hz_i - \hz_j)^2}  - \int_{\eta_{\l,2} \le |y - \hz_i| \le \eta_\l}  \frac{f(y) - f(\hz_i)}{(\hz_i - y)^2} \sc(0) \, dy\right) \eeq
The difference is that if $-i\in{A_{q_*}\setminus A_{2, (i)}}$ we omit the $j= - i$ term. We will show in this case that the term is negligible, so that it can be reinserted and the proof can be completed in the same way.

By the mean value theorem, there exists $\xi\in \R$ such that
\beq \frac{1}{2N} \frac{f_{-i} -f_{i}}{(\hz_{-i} - \hz_{i})^2} = \frac{1}{2N} \left(  \frac{f'_{-i}}{(\hz_{-i} - \hz_{i})} + \frac{f''_{-i}(\xi)}{2} \right). \eeq
Then, using Lemma \ref{lem:uniformlocal}, Lemma \ref{lem:fbound}, and the hypothesis on $i$, we have the bound
\beq \left|  \frac{1}{2N} \frac{f_{-i} -f_{i}}{(\hz_{-i} - \hz_{i})^2} \right| \le  \frac{1}{2N} \left( N \frac{N^{\eps}}{N^{\omega_{\l,2}}}\frac{1}{t-s_1 +s_0}\frac{1}{t+s_1}  + \frac{1}{(t-s_1 +s_0)^2}\frac{1}{t+s_1} \right). \eeq 
This is the same size as the error obtained for (3.137) and (3.138) in the bound (3.142) of \cite{fixed}.\ep
  
  We obtain the following lemma by integrating. 
  \bec Fix $\eps, D>0$. For each $\alpha$ there exists an event $\Fa$ such that $\P(\Fa)\ge 1 - N^{-D}$ and on which the following estimate holds. 
 \beq \int_{s_1}^{2t_1} \langle (w-f), \B (w-f) \, dx\rangle \le C\| (w(s_1) - f(s_2) \|_2^2  + \frac{N^\eps}{s_1}\frac{1}{\sqrt{Ns_0}}\eeq
  \eec
  
 We now obtain a time-averaged version of our desired result. 
 
\bet Fix $a$ and $i$ such that $|a| \le N^{\omega_A - \eps_B}$ and $|i-a|\le \l/10$. Fix also $\eps, D>0$. For every $\alpha$ there is an event $\Fa$ such that $\P(\mathcal F_\alpha)\ge 1  - N^{-D}$ and on which 
\beq  \frac{1}{t_1}\int_0^{t_1} \left( U_{t_1+u}(i,a) - \frac{1}{N}p_{t_1 +u}(\flat_i, \flat_a)\right)^2 \, du \le \frac{N^\eps}{(Nt_1)^2} \left(  \frac{(Nt_1)^4}{\l^4}  + \frac{s_1^2}{t_1^2}  + \frac{t_1}{s_1} \left( \frac{1}{\sqrt{Ns_0}} + \frac{s_0}{s_1}  \right)  \right)
\eeq
for $0\le u \le t_1$.
\eet
\bp
We follow the proof of Theorem 3.15 in \cite{fixed}. There are two changes. First, we must modify the bound in display (3.165) to use our kernel $\B$. Recall that we omit the index $i=0$ in $w$ and $f$. As in our discussion of Lemma \ref{lem:initial}, we shift the positive indices down one in order to apply a Sobolev inequality in Appendix D of \cite{fixed}. As in the proof of Theorem 3.15 of that reference, we apply the Sobolev inequality to obtain
\beq \frac{1}{4N^2}(w_{t_1 + u}(i)  - f_{t_1 + u}(i))^2 \le  C\left( \frac{1}{N\l}\sum_{|j-i|\le \l} w_{t_1 +u}(j) -\frac{1}{N\l}\sum_{|j-i|\le \l} f_{t_1 +u}(j) \right)^2\eeq \beq + C\log(N)N^{-2} \left( \sum_{|i|, |j| \le l, j\neq i, -(i+1)} \frac{[(w-f)_i(t_1+u) - (w-f)_{j}(t_1+u)]^2}{(i-j)^2} \right)\eeq \beq  + C\log(N)N^{-2}\sum_i [(w-f)_i(t_1+u) - (w-f)_{-i}(t_1+u)]^2 .  \eeq 
The first two terms are bounded as in the proof in \cite{fixed}. In particular, using rigidity, we see the second is bounded by
\beq \frac{N^\eps}{N^2} \langle (w-f)(t_1 + u) , \B (w-f)(t_1 + u)\rangle. \eeq
The third term is, using Lemma \ref{lem:wbound} and the analogous bound for $\|f\|^2$, bounded by
\beq  \frac{C \log{N}}{N^2} \| (w -f)(t_1+u) \|^2  \le \frac{C \log{N}}{N^2} \left( \| w(t_1 +u)\|_2^2 +  \| f(t_1 +u)\|_2^2 \right) \le \frac{N^\eps}{N^2}\left( \frac{1}{t_1} + \frac{1}{s_1}  \right).\eeq
This error is bounded by $N^{-1+\delta}$ for some small $\delta>0$, which smaller than the errors obtained when bounding (3.165) in \cite{fixed}. Hence replacing the old kernel with the new one in the proof is permissible. 

Second, in (3.170) we apply our Lemma \ref{lem:initial} and change (3.171) to use our $\B$.
\ep

The removal of the time average as in Theorem 3.16 of \cite{fixed} goes through without change, since it only uses the abstract properties of $\U$ and the bounds in Lemma \ref{lem:pbound}. 

\bet\label{thm:noavg} Fix $a$ and $i$ such that 
\beq |a| \le \frac{N^{\omega_A - \eps_B}}{2},\quad |i-a| \le \frac{\l}{20}.\eeq
Fix also $\eps, D>0$. There exists an event $\Fa$ such that $\P(F_\alpha)\ge 1  - N^{-D}$ and on which the following bound holds.
\beq \left| \mathcal U_{t_1 + 2t_2}(i,a) - \frac{1}{N} p_{t_1}(\flat_i , \flat_a) \right| \le CN^{\eps + \eps_2 -1} t_1^{-1} \left(  \frac{s_1^2}{t_1^2}  + \frac{(Nt_1)^4}{\l^4} + \frac{t_1}{s_1}\left[ \frac{1}{\sqrt{Ns_0}} + \frac{s_0}{s_1} \right]\right)^{1/2} + N^{\eps + \eps_2/2 -1}t_1^{-1} \eeq

\eet

Then as in the proof of Theorem 3.10 in \cite{fixed}, choosing $s_0$ and $s_1$ in Theorem \ref{thm:noavg} such that $Ns_0 = (Ns_1)^{2/3}$ and $Ns_1 = (Nt_1)^{9/10}$ yields the final homogenization result.

\bet\label{thm:homog} Fix $a$ and $i$ such that $|a| \le N^{\omega_A - \eps_B}$ and $|i-a|\le \l/10$. With $t_1$ as above and $t_2 = N^{-\eps_2}t_1$ (for $\eps_2$ satisfying $\omega_1 - \eps_2 >0$), and fixed $\eps, D>0$, there exists for every $\alpha$ an event $\Fa$ with $\P(\Fa) \ge 1 - N^{-D}$ on which the following holds. 
\beq  \left|  \U_{ia}(0,t_1) - \frac{1}{N}p_{t_1}(\flat_i,\flat_j) \right| \le \frac{N^{\eps+\eps_2}}{Nt_1}\left(\frac{(Nt_1)^2}{\l^2} + \frac{1}{(Nt_1)^{1/10}}   \right) + \frac{N^{\eps - \eps_2/2}}{Nt_1}\eeq
\eet
\subsection{Conclusion}\label{sec:conclusion}
Here we will frequently need to use the above bounds, which hold for fixed $\alpha$ on a set of large probability, for all $\alpha\in[0,1]$ simultaneously when bounding an integrand. This is accomplished using Lemma E.1 and the Remark that follows it in Appendix E of \cite{fixed}.\\

\noindent {\bf Proof of Theorem \ref{thm:homomain}.}\ \  By Lemma \ref{lem:hatapprox},
\beq z_i(t_1, 1) - z_i(t_1,0) = \hz_i(t_1, 1) - \hz_i(t_1,0) + (\gamma_1(t_1,1) - \gamma_1(t_1,0)) + O\left( N^\eps t_1 \left( \frac{N^{\omega_A}}{N^{\omega_0}} + \frac{1}{N^{\omega_\l}}  + \frac{1}{\sqrt{NG}}\right)  \right). \eeq
The last term is $o(1)$, and by Lemma \ref{lem:gammarigid} the difference $(\gamma_1(t_1,1) - \gamma_1(t_1,0))$ is also $o(1)$ for $j,k \le N^{\omega_0/2}$. Hence it suffices to bound the first term. 

Recalling $u_i = \partial_\alpha \hz_i$, we have
\beq \hz_i(t_1, 1) - \hz_i(t_1,0)  = \int_0^1 u_i(t_1,\alpha)\, d\alpha.\eeq
We now use the finite speed of propagation estimates to show that the initial data far from zero makes a negligible contribution to $u_i(t_1,\alpha)$. We perform this estimate in two steps. First, we remark that $u$ that is a solution of
\beq\partial_t u  = \mathcal B u + \xi,\eeq
where $\xi$ satisfies the bound
\beq\label{eqn:xiest} | \xi_i(t)| \le \one_{\{|i| > N^{\omega_A}\}} N^C\eeq
with overwhelming probability for $0\le t \le 1$. We define $v_i$ as the solution
\beq \partial_t v = \mathcal B v,\quad v_i(0) = u_i(0).\eeq
Then by the Duhamel formula,
\beq u_i(t_1) - v_i(t_1) = \int_0^1 \sum_{|p| \le N} \mathcal U_{ip}(s,t_1)\xi_{p}(s)\, ds = \int_0^1 \sum_{N^{\omega_A} <|p| \le N} \mathcal U_{ip}(s,t_1)\xi_{p}(s)\, ds.\eeq
We fix $\delta_B>0$ and consider $i$ such that $|i| \le N^{\omega_A - \delta_B}$. By Lemma \ref{lem:finprop} and \eqref{eqn:xiest}, we obtain
\beq |  u_i(t_1) - v_i(t_1) | \le N^{-10}.\eeq

For the second step, we fix $\eps_a>0$ and consider $w$ defined as 
\beq\partial_t w = \mathcal B w,\quad w_i(0) = v_i(0)\one_{\{ |i| \le N^{1+\eps_a} t_1\}}.\eeq
We have
\beq v_i(t_1) - w_i(t_1) = \sum_{N^{\omega_1+\eps_a}}\mathcal U_{ij}(0,t_1)u_j(0).\eeq
As in the first step, Lemma \ref{lem:finprop} shows the terms with $|j| > N^{\omega_A}$ are negligible. Fix $\eps_b > 0$ such that $\eps_b < \eps_a$. Then by Lemma \ref{lem:poisson},
\beq |v_i(t_1) - w_i(t_1)| \le \left| \sum_{N^{\omega_1 + \eps_a} < |j| \le N^{\omega_A}}  \mathcal U_{ij}(0,t_1)u_j(0)  \right| \le N^{\eps + \omega_1} \sum_{|j| > N^{\omega_1 + \eps_a}} \frac{1}{(i-j)^2}\le N^{-1-\eps_a + \eps}.\eeq
We then have
\beq \left| \int_0^1 u_i(t_1,\alpha)\, d\alpha \right| \le \left|\int_0^1 \sum_{|j| \le Nt_1 N^{\eps_a}} \U(0, t_1,\alpha)_{ij}u_j(0)\, d\alpha\right|  + N^{-1-\eps_a + \eps} .\eeq
Therefore, it suffices to estimate 
\beq \int_0^1 \sum_{|j| \le Nt_1 N^{\eps_a}} \U(0, t_1,\alpha)_{ij}(z_j(0,1) - z_j(0,0)) \, d\alpha.\eeq

Recall that by hypothesis the initial data are symmetric, $z_j = - z_{-j}$. Hence we just need to estimate
\beq \int_0^1 \sum_{ 0 < j \le Nt_1 N^{\eps_a}} (\U_{ij}(0, t_1) - \U_{i,-j}(0,t_1)(z_j(0,1) - z_j(0,0)) \, d\alpha.\eeq
By Theorem \ref{thm:homog} we can replace $\U_{ij}(0, t_1)$ with $\frac{1}{N}p_{t_1}(\flat_i,\flat_j)$ and accrue an error of 
\beq E_0=\frac{N^{\eps+\eps_2}}{Nt_1}\left(\frac{(Nt_1)^2}{\l^2} + \frac{1}{(Nt_1)^{1/10}}   \right) + \frac{N^{\eps - \eps_2/2}}{Nt_1}. \eeq Then
\beq (\U_{ij}(0, t_1) - \U_{i,-j}(0,t_1) \le \frac{1}{N}(p_{t_1}(\flat_i,\flat_j) - p_{t_1}(\flat_i, -\flat_j) ) + E_0\eeq \beq=  \frac{1}{N}(p_{t_1}(\flat_i,\flat_j) - p_{t_1}(-\flat_i, \flat_j) ) + E_0 = \frac{1}{N}\int_{-\flat_1}^{\flat_1} p'_{t_1}(x, \flat_j) \, dx+ E_0 .\eeq
Summing over $j \le Nt_1 N^{\eps_a}$, then using Lemma \ref{lem:pbound} and the normalization of $p_t$ gives
\beq \sum_{ 0 < j \le Nt_1 N^{\eps_a}} (\U_{ij}(0, t_1) - \U_{i,-j}(0,t_1) \le  \sum_{ 0 < j \le Nt_1 N^{\eps_a}}  E_0 + \frac{1}{N}\int_{-\flat_1}^{\flat_1} p'_{t_1}(x, \flat_j) \eeq \beq\le    Nt_1 N^{\eps_a} E_0 +\frac{1}{ t_1}  \int_{-\flat_1}^{\flat_1}  \sum_{0<j \le Nt_1 N^{\eps_a}}  \frac{1}{N} p_{t_1}(x,\flat_j) \, dx \le Nt_1 N^{\eps_a} E_0 +  \frac{C}{t_1} \int_{-\flat_1}^{\flat_1} 1\, dx \le  Nt_1 N^{\eps_a} E_0 + \frac{N^\eps}{Nt_1}. \eeq 
Choosing $\eps, \eps_2, \eps_a >0$ small enough shows that for large enough $N$, 
\beq  \sum_{ 0 < j \le Nt_1 N^{\eps_a}} (\U_{ij}(0, t_1) - \U_{i,-j}(0,t_1) \le N^{-\delta} \eeq
for some small $\delta >0$. Finally, using Lemma \ref{lem:uniformlocal} we have with overwhelming probability
\beq \sup_{j \le Nt_1 N^{\eps_a}} |z_j(0,1) - z_j(0,0)| \le N^{-1+\eps}\eeq
for arbitrarily small $\eps >0$. Hence 
\beq \int_0^1 \sum_{|j| \le Nt_1 N^{\eps_a}} \U(0, t_1,\alpha)_{ij}(z_j(0,1) - z_j(0,0)) \, d\alpha  \le \int_0^1 \sum_{|j| \le Nt_1 N^{\eps_a}} |\U(0, t_1,\alpha)_{ij}| |(z_j(0,1) - z_j(0,0))| \, d\alpha \le N^{-1-\delta'}\eeq
for some $\delta'>0$. Undoing the time shift in the definition of $z_i(t,\alpha)$ completes the proof.
%\frac{1}{N} \int_{-\flat_j}^{\flat_j} p'_{t_1}(\flat_i, x)\, dx \le \frac{2\flat_j}{Nt_1}\left( \frac{1}{t_1}+ N^{-D}\right)\le \frac{2j N^\eps}{(N t_1)^2}  \eeq
%\beq \le \frac{Nt_1 N^{\eps_a}}{(Nt_t)^2} \le \frac{N^{\eps_a}}{Nt_1}\eeq
%Then
%\beq \hz_i(t_1, 1) - \hz_i(t_1,0)  \le N^\eps \| z(0,0) - z(0,1)\|_\infty \frac{1}{(Nt_1)^2}\sum_{j=1}%^{(Nt_1)N^{\eps_a}} j \le \| z(0,0) - z(0,1)\|_\infty \frac{1}{(Nt_1)^2}  \eeq 
\ep

\section{Deformed local law} \label{sec:dqclaw}

In this section we prove the deformed local law, Theorem \ref{thm:deformedlaw}, necessary for the proof of Lemma \ref{lem:uniformlocal}. A deformed local law for eigenvalues was previously shown in \cite{LY}.

The structure of this section is as follows. In Subsection \ref{sec:preliminaries} we define some notions necessary for the rest of the section, state the main result, and compute Green function elements. In Subsection \ref{sec:bootstrapestimates} we establish some preliminary estimates, then prove Lemma \ref{lem:weaklaw}, a weak version of the deformed local law. Finally, in Subsection \ref{sec:deformedlawproof}, we use a fluctuation averaging argument to upgrade this weak law and prove Theorem \ref{thm:deformedlaw}.

\subsection{Preliminaries}\label{sec:preliminaries}

\subsubsection{Deformed Semicircle Law}

It is well known that in the large $N$ limit, many symmetric matrix ensembles have a macroscopic eigenvalue density that obeys the semicircle law, whose density is
$$\sc(x) = \one_{\{|x|<2\}} \frac{1}{2\pi} \sqrt{4-x^2}.$$
The \emph{free convolution} of the semicircle law and a deterministic diagonal matrix $\overline V = (v_i)_{i=-N}^N$, with the $i=0$ index omitted, is defined by its Stieltjes transform,
$$\mfc^{(N)}(z) = \frac{1}{2N} \sum_{i=-N}^N \frac{1}{v_i - z -t\mfc^{(N)}(z)},$$
where again the $i=0$ term is omitted in the summation. There is a unique solution to this equation, and it is the Stieltjes transform of a measure absolutely continuous with respect to Lebesgue measure. We call the associated density $\dsc^{(N)}$. These facts and basic properties of the free convolution are proved in \cite{Bi97}. We often write $\mfc$ and $\dsc$ for these quantities, omitting dependence on $N$. 

We now state a result on the stability of $\mfc$, proved in Lemma 7.2 of \cite{LY}. Recall the notion of $(g,G)$-regularity introduced in Definition \ref{def:regular}.
\bel\label{lem:stability} Assume that $V$ is $(g,G)$-regular. For $q\in(0,1)$, $\sigma>0$, and $N$ large enough, the following statements hold for $E\in (-qG, qG)$, $\eta \in [N^{-5}, 10]$, and $t$ such that $g N^\sigma \le t \le N^{- \sigma}G^2$. The constants do not depend on $\sigma$ or $q$. 

\begin{enumerate}[label=(\roman*)]
\item We have 
\beq c\le \Im \mfc(z) \le C \eeq 
and hence 
\beq\label{eqn:stability} ct\le |v_i - z - t\mfc (z)| \eeq
for all $v_i$. Both statements hold uniformly in $N$.

\item We have 
$$ c \le \left|  1  - \frac{t}{2N} \sum_{-N}^N \frac{1}{( v_i - z -t\mfc (z))^2} \right| \le C$$
uniformly in $N$.
\end{enumerate}

\eel

\subsubsection{Stochastic Domination}

We recall the notion of stochastic domination introduced in \cite{EKY13}.
\bed\label{d:sd} Let 
\beq X = (X^{(N)}(u)\colon N\in \N, u\in U^{(N)}), \quad Y=(Y^{(N)}(u)\colon N \in \N, u\in U^{(N)}) \eeq
be two sets of nonnegative random variables, where $U^{(N)}$ is a possibly $N$-dependent parameter set. We say that $X$ is \emph{stochastically dominated by $Y$, uniformly in $u$}, if for all $\eps > 0$ and $D>0$ we have
\beq \sup_{u \in U^{(N)}} \P\left[  X^{(N)}(u) > N^\eps Y^{(N)}(u) \right] \le N^{-D} \eeq
for large enough $N\ge N_0(\eps, D)$. The stochastic domination is always uniform in all parameters that are not explicitly stated.

If $X$ is stochastically dominated by $Y$, uniformly in $u$, we write $X\prec Y$. If for some complex family $X$ we have $|X|\prec Y$, we also write $X = O_\prec (Y)$.
\eed

Observe that a sequence of events $E=(E^{(N)})$ holds with overwhelming probability if $1 - \one(E) \prec 0.$ We also recall the basic properties of $\prec$. 
\bel
\begin{enumerate}[label=(\roman*)]
\item Suppose that $X(u,v)\prec Y(u,v)$ uniformly in $u\in U$ and $v\in V$. If $|V|\le N^C$ for some constant $C$, then
$$\sum_{v\in V} X(u,v) \prec \sum_{v\in V} Y(u,v)$$ uniformly in $u$.
\item Suppose that $X_1(u)\prec Y_1(u)$ uniformly in $u$ and $X_2(u) \prec Y_2(u)$ uniformly in $u$. Then $X_1(u)X_2(u) \prec Y_1(u)Y_2(u)$ uniformly in $u$.
\item if $X\prec Y + N^{-\eps}X$ for some $\eps>0$, then $X\prec Y$.
\end{enumerate}
\eel

The following large deviations estimates will be important for our work. Proofs may be found in, for example, \cite{lectures}.

\bel\label{lem:largedev} Let $\left(  X^{(N)}_i \right)$, $\left(  Y^{(N)}_i \right)$, $\left(  a^{(N)}_{ij} \right)$, and $\left(  b^{(N)}_{ij} \right)$ be independent families of random variables, where $N\in\N$ and $i,j\in\{1,\dots, N\}$. Suppose that all entries $X_i^{(N)}$ and $Y_i^{(N)}$ are independent and satisfy
$$\E X = 0, \quad \| X \|_p \le \mu_p$$
for all $p\in \N$ with some constants $\mu_p$. In the following statements, $\Psi$ can be any random variable.

\benr

\item Suppose that $\left(\sum_i |b_i|^2\right)^{1/2} \prec \Psi$. Then $\sum_i b_i X_i \prec \Psi$.
\item Suppose that $\left( \sum_{i\neq j} |a_{ij}|^2 \right)^{1/2} \prec \Psi$. Then $\sum_{i\neq j} a_{ij} X_i X_j \prec \Psi$.

\item Suppose that $\left( \sum_{i,j} |a_{ij}|^{1/2} \right)^{1/2} \prec \Psi$. Then $\sum_{i,j} a_{ij} X_i Y_j \prec \Psi$. 
\eenr
If all of the above random variables depend on an index $u$ and the hypotheses of (i) -- (iii) are uniform in $u$, then so are the conclusions.

\eel

\subsubsection{Model and main result}

We consider as in Section \ref{sec:deterministic} a deterministic initial data matrix $V$. We let $W$ be a matrix of i.i.d. random variables with distribution $\mathcal N(0, N^{-1})$ and define
\beq H_t = V + \sqrt{t} W, \quad M_t = \bma 0 & H_t  \\   H_t^\trans & 0 \ema.\eeq
Recall the eigenvalues of $M_t$, which we will call $\lambda_i(t)$, are the singular values of of $H_t$ and their negatives. 
We also define 
\beq m_N(z)  = \frac{1}{N} \sum_{i=-N}^N \frac{1}{\lambda_i(t) - z}.\eeq
suppressing the dependence of $m_N$ on $t$ in our notation.

In what follows we fix $\delta_1>0$ and always suppose that $V$ is $(g,G)$-regular, as defined in Definition \ref{def:regular}. We also suppose $z\in \mathcal D$ and $t\in T_\sigma$, which are defined by fixing $\sigma >0$, $q\in (0,1)$, and setting
\beq T_\sigma = \{ t \colon g N^\sigma \le t \le N^{- \sigma}G^2 \}, \quad \mathcal D =\{ z = E + i\eta \colon E\in (-qG,qG), N^{\delta_1} \le N\eta \le 10N \}\eeq

The following theorem is the main result of this section. 
\bet\label{thm:deformedlaw} Let $H_t$ be defined as above for some $(g, G)$-regular $V$. Fix $q\in (0,1)$ and $\sigma> 0$. Uniformly for $t\in T_\sigma$ and $z\in \mathcal D$, we have
\beq\left|m_N(z) -\mfc(z) \right| \prec \frac{1}{N\eta}.\eeq
\eet

\subsubsection{Reduction to diagonal $V$}

In the proof below, we will assume that $V$ is diagonal with entries $\{v_i\}_{i = 1}^N$. If $V$ is a general deterministic matrix, it has a singular value decomposition $V = A D B^\trans$ with $D$ diagonal and $A, B$ orthogonal. Since the Gaussian ensemble $W$ defined above is invariant under multiplication on the left and right by unitary matrices, we have
\beq V + \sqrt{t} W = A( D + \sqrt{t} W' )B^\trans\eeq
for a Gaussian ensemble $W'$ with the same entry distribution as $W$. Hence $V+\sqrt{t} W$ and $D+ \sqrt{t} W'$ have the same distribution of eigenvalues, and we have shown that without loss of generality we may consider diagonal initial data.  

\subsubsection{Green Functions}\label{sec:greens}
We define the Green function matrix 
\beq G(z) = (M_t -zI)^{-1}.\eeq
For concreteness, we compute $G_{1,1}$. The other $G_{i,i}$ are analogous. We often write $G_{ij}$ for $G_{i,j}$. 

In the study of Hermitian random matrices, it is common to compute the diagonal Green function elements $G_{ii}$ using the Schur complement formula. The off-diagonal entries $G_{ij}$ where $i\neq j$ are asymptotically small and can be neglected. The Schur formula then produces an approximate fixed point equation for the diagonal elements $G_{ii}$. By analyzing the fixed point equations, one may bound the $G_{ii}$ and obtain useful information about the eigenvalues and eigenvectors.

 In the non-Hermitian setting, we would like to take the same approach. However, now the asymptotically non-trivial entries are not only the diagonal entries like $G_{ii}$, but also off-diagonal entries like $G_{i,i+N}$ and $G_{i+N, i}$, and the equations for the elements $(G_{ii}, G_{i,N+i}, G_{N+i,i}, G_{N+i,N+i})$ are closely coupled together. It turns out that considering these $2\times 2$ blocks instead of individual entries yields a similar formulation to the Hermitian case, though one that requires more careful bookkeeping. This motivates the decision to work with $2\times 2$ blocks in the following computation.

Applying the Schur complement formula with index set $\T = (1, N+1)$ yields 

\beq\bma
G_{1,1}& G_{1, N+1}\\
G_{N+1,1} & G_{N+1, N+1}
\ema
= (B_t - t v^\trans G^{(\T)} v)^{-1},\eeq
where
\beq B_t = \bma -z & v_1 + \sqrt{t} w_{11}\\  v_1 + \sqrt{t}w_{11}& -z   \ema, \quad v= \bma 
0 & w_{2,1}\\
0 & w_{3,1}\\
\vdots & \vdots \\
 w_{1,2} & 0 \\
 w_{1,3} & 0 \\
\vdots & \vdots \ema. \eeq
We compute
$$v^\trans G^{(\T)} v= \bma A_{11} & A_{12} \\ A_{21} & A_{22}\ema $$
\beq A_{11} = \sum_{i,j=1}^{N-1} w_{1, j +1} G^{(\T)}_{N-1+j, n -1 +i} w_{1, i+1},\quad A_{12}= \sum_{i,j=1}^{N-1} w_{1, j +1} G^{(\T)}_{N-1+j, i} w_{i+1, 1}, \eeq
\beq A_{21}=\sum_{i,j=1}^{N-1} w_{j+1, 1} G^{(\T)}_{j, n -1 +i} w_{1, i+1},\quad A_{22}= \sum_{i,j=1}^{N-1} w_{j +1,1} G^{(\T)}_{j, i} w_{ i+1,1}.\eeq

Notice that $A_{12}=A_{21}$, as $G^{(\T)}$ is symmetric. Further, we may apply the Schur complement formula to the $(N-1)\times (N-1)$ blocks of $G^{(\T)}$ and take the trace to see
\beq \sum_{i=1}^{N-1}  G^{(\T)}_{N-1+i, N -1 +i} = \sum_{i=1}^{N-1}  G^{(\T)}_{i, i} ,\eeq
so $\E A_{11}= \E A_{22}=m^{(\T)}(z)$.

From the formula for the inverse of a $2\times 2$ matrix, we find
\beq G_{11} = \frac{-z-tA_{22}}{(-z-tA_{11})(-z-tA_{22}) - (v_1+\sqrt{t}w_{11} -tA_{12})(v_1+\sqrt{t}w_{11} -tA_{21}) },\eeq

\beq =\frac{-z-tA_{22}}{(-z-(t/2)(A_{11}+A_{22}))^2- ((t/2)(A_{11}-A_{22}))^2- (v_1+\sqrt{t}w_{11} -tA_{12})(v_1+\sqrt{t}w_{11} -tA_{12}) }.
\eeq
Set \beq E_1=t^2 ((A_{11}-A_{22})/2)^2, \quad E_2=(A_{11}+A_{22})/2 - m^{(\T)}.\quad r=m-m^{(\T)}.\eeq Writing $m$ for $m_N(z)$, we have
\beq =\frac{-z-tA_{22}}{  ( -z-tm^{(\T)}+ v_1+\sqrt{t}w_{11} -tA_{12}-tE_2)(-z-tm^{(\T)}-v_1- \sqrt{t}w_{11} +tA_{12}-tE_2)- E_1},\eeq

\beq =\frac{-z-tA_{22}}{  ( -z-tm+ v_1+\sqrt{t}w_{11} -tA_{12}-tE_2 +t r )(-z-tm-v_1- \sqrt{t}w_{11} +tA_{12}-tE_2+ tr)- E_1}. \eeq

Set, for $i>0$, 
\beq g_i = \frac{1}{-z-tm +v_i}, \quad  g_{-i} = \frac{1}{-z-tm - v_i}. \eeq
We define 
\beq E_3 =  g_1(\sqrt{t}w_{11} -tA_{12}-tE_2+ tr) \eeq
and $E_4$ similarly. Then

\beq G_{11}=\frac{-z-tA_{22}}{  ( -z-tm+ v_1)(1+E_3)(-z-tm-v_1)(1+E_4)- E_1}, \eeq

\beq =g_1 g_{-1} \frac{-z-tA_{22}}{[(1+E_3)(1+E_4)-g_1  g_{-1}E_1]}. \eeq

Define
\beq K_1=E_3+E_4+E_3E_4-g_1 g_{-1} E_1,\quad B_1=\frac{-t([A_{22}-m^{(\T)}]-r)}{-z-tm}.\eeq
Our expression for $G_{ii}$ is, with $K_i$ and $B_i$ defined similarly,
\beq\label{grep} G_{ii} = \frac{1}{2}\left( g_i +  g_{-i}  \right)\frac{1+B_i}{1+K_i}.\eeq

\subsection{Weak law}\label{sec:bootstrapestimates}

We prove some bounds necessary for our bootstrap argument.

\bel\label{lem:gisum} For any $z\in \mathcal D$, $|\mfc(z) - m(z) | \le (\eta N)^{-1/2}$ implies
\beq \frac{1}{N} \sum_i |g_i| + |g_{-i}| \le C \log N.\eeq
\eel
\bp
This follows as in the proof of Lemma 7.5 in \cite{LY}, with the minor change that we are using $m$ instead of $\mfc$ in the definition of $g_i$. However, if $|\mfc - m | \le (\eta N)^{-1/2}$, then by Lemma \ref{lem:stability} we also have a lower bound $\Im m \ge c$, and this suffices to complete the proof.
\ep

\bel\label{lem:bootstrap}

Fix $z\in \mathcal D$. Let $\phi$ be the indicator function of some event, which may depend on $z$. If $\phi| m(z) - \mfc(z) | \prec N^{-c}$ for some $c>0$, then

\beq \phi \max_i |K_i| \, {\prec}\frac{1}{\sqrt{N\eta}} , \quad  \phi \max_i |B_i| \, {\prec}\frac{1}{\sqrt{N\eta}} ,\eeq
\beq \frac{\phi t^2(\mfc -m)}{2N} \sum_{-N}^N \frac{1}{(-z-t\mfc + v_i)^2(-z-tm +v_i)} = O_{\prec} (N^{-c}).\eeq
\eel 

\bp We first consider $B_1$. Note that by the assumption we have a stability bound for $m$ similar to Lemma \ref{lem:stability} and this shows 
\beq \frac{t}{-z-tm} =O(1), \eeq
so it suffices to bound $r$ and $A_{22} - m^{(\T)}$. We have 
\beq |r| \le \frac{C}{N\eta}, \quad \left|A_{22} - m^{(\T)}\right| \prec \frac{1}{\sqrt{N\eta}} \eeq 
The first inequality follows from the eigenvalue interlacing lemma. (See, for example, Lemma 7.5 in \cite{EYbook}.) For the second, we apply Lemma \ref{lem:largedev} and then Ward's identity to obtain
\beq \left| A_{22}-m^{(\T)} \right| \prec \frac{1}{N}\sqrt{\sum_{i,j} |G_{ij} ^{(\T)}|^2}  = \sqrt{\frac{\Im m^{(\T)}}{N\eta}} =\sqrt{\frac{\Im \mfc + r + (m-\mfc)}{N\eta}} \prec  \frac{1}{\sqrt{N\eta}}.\eeq
The final bound follows from the assumption on $m-\mfc$ and the upper bound on $\Im \mfc$ in Lemma \ref{lem:stability}. Finally, note these bounds are independent of $i$.

We now consider $K_1$. By the stability bound (\ref{eqn:stability}),
\beq |g_1 g_{-1} t^2| \le c,\eeq
and by the large deviations bound Lemma \ref{lem:largedev}, \beq \frac{A_{11} - A_{22} }{2} = O_\prec \left(\frac{1}{\sqrt{N\eta}} \right).\eeq
It remains to bound $E_3$ and $E_4$.  We just do $E_3$, as $E_4$ is similar. Using (\ref{eqn:stability}), and noting that  ${ \min (ct, \eta) \le |v_1 -z -t\mfc | }$, which implies ${ \sqrt{c\eta t}  \le |v_1 -z -t\mfc | }$,
\beq |E_3|\le  \frac{|\sqrt{t}w_{11} -tA_{12}-tE_2 +tr |}{|v_1-z-t\mfc|}\le \frac{|w_{11}|}{\sqrt{c\eta}} + |A_{12}| +|E_2|+ |r| .\eeq
The $r$ term was already bounded, and similar large deviations arguments apply to $A_{12}$ and $E_2$. Because $w_{11}$ has subexponential decay,
\beq w_{11} \prec \frac{1}{\sqrt N}.\eeq
Again, these bounds are independent of $i$.

For the final claim, we apply the Cauchy-Schwarz inequality and note
\begin{multline}\label{eqn:cs} \left| \frac{1}{2N}  \sum_{-N}^N  \frac{1}{(-z-t\mfc + v_i)^2(-z-tm +v_i)} \right| \\ \le \left(\frac{1}{2N}\sum_{-N}^N \frac{1}{|-z-t\mfc +v_i |^4} \right)^{1/2} \left(\frac{1}{2N} \sum_{-N}^N \frac{1}{|-z-tm+v_i|^2} \right)^{1/2}\end{multline}
We first bound the first factor in the right side of (\ref{eqn:cs}). By the stability bound Lemma \ref{lem:stability}, we get the bound
\beq\label{eqn:ff} \frac{1}{2N}\sum_{-N}^N \frac{1}{|-z-t\mfc +v_i |^4}  \le \frac{C}{2Nt^2}\sum_{-N}^N \frac{1}{|-z-t\mfc +v_i |^2}  \eeq
Taking imaginary parts in the equation that defines $\mfc$ gives
\beq \Im \mfc =  \frac{1}{2N}  \sum_1^{2N} \frac{ t\Im \mfc +\eta }{ |-z-t\mfc + v_i|^2},\eeq
so that \beq \frac{1}{2N}  \sum_1^{2N}\frac{ 1 }{ |-z-t\mfc + v_i|^2} = \frac{\Im \mfc}{t\Im \mfc+\eta} \le \frac{1}{t}. \eeq
We obtain
\beq \left( \frac{1}{2N}\sum_{-N}^N \frac{1}{|-z-t\mfc +v_i |^4} \right)^{1/2} \le  \frac{C}{t^{3/2}}. \eeq

We now consider the second factor in (\ref{eqn:cs}). We may sum the representation (\ref{grep}) for $G_{ii}$ to obtain
\beq\label{eqn:sc1} m=\frac{1}{2N} \sum_{-N}^{N} \frac{1}{-z-tm + v_i} + R,\eeq
with $R=O_\prec \left( \frac{1}{\sqrt{\eta N}}\right)$, where we have used the bounds for $K_i$ and $B_i$ to Taylor expand 
\beq \frac{1+B_i}{1+K_i} = 1 + R'\eeq
for large $N$, with $R'=O_\prec \left( \frac{1}{\sqrt{\eta N}}\right)$, and then Lemma \ref{lem:gisum} to control the sum of the errors. There is a factor of $\log(N)$ that is absorbed by the stochastic domination. 

Taking imaginary parts in (\ref{eqn:sc1}) gives 
\beq \frac{\Im m}{\eta + t\Im m} = \frac{1}{2N} \sum_{-N}^{N} \frac{1}{|-z-tm + v_i|^2} + \frac{\Im R}{\eta + t\Im m}.\eeq
Then, using Lemma \ref{lem:stability}  to write $2\Im m \ge \Im m - R$ for $z\in \mathcal D$ and absorb the error term, we have \beq   \left( \frac{1}{2N} \sum_1^{2N} \frac{1}{|-z-tm + v_i|^2}\right)^{1/2} \le  \frac{C}{t^{1/2}}.  \eeq
 
Putting these bounds on each factor together, the expression we want to control is bounded above by
\beq C \phi|m - \mfc| t^2   \frac{1}{t^{3/2}} \frac{1}{t^{1/2}} \prec  N^{-c} .\eeq
\ep

We also prove \emph{a priori} bounds.
\bel\label{lem:apriori}
If $\eta \ge 1$, then, 
\beq \max_i K_i = O_{\prec}\left(\frac{1}{\sqrt{N\eta}} \right), \quad  \max_i B_i = O_{\prec}\left(\frac{1}{\sqrt{N\eta}} \right),\eeq
\beq \frac{t^2}{2N} \sum_{-N}^N \frac{1}{(-z-t\mfc + v_i)(-z-tm +v_i)} \le 2t^2.\eeq
\eel 
\bp
The first two bounds are proved as before, except we use the trivial estimate
\beq \Im m \le \frac{1}{\eta}\eeq
instead of $|m-\mfc | \le N^{-c}$ to estimate $\Im m^{(\T)}$. 

For the last bound, we use
\beq |-z - tm + v_i|\ge \eta, \quad |-z - t\mfc + v_i| \ge \eta, \quad |m| \le \frac{1}{\eta}, \quad |\mfc| \le \frac{1}{\eta}.\eeq
Then \beq \frac{t^2(m - \mfc)}{2N} \left|\sum_{-N}^N \frac{1}{(z-tm + v_i)(-z -t\mfc v_i)} \right| \le \frac{2}{\eta}\frac{t^2}{2N} \frac{2N}{\eta^2} \le 2t^2. \eeq
\ep

We now prove the weak local deformed law at the optimal scale using a bootstrapping argument. Our presentation follows \cite{lectures}. 

\bel\label{lem:weaklaw} Suppose the initial values $V$ are $(g,G)$-regular as in Definition \ref{def:regular}. Then for $z\in \mathcal  D$, we have $|m - \mfc| \prec (N\eta)^{-1/2}$. 
\eel

\bp First, note that both $m$ and $\mfc$ are $N^2$-Lipschitz continuous on $\mathcal D$. For $m$ this is well known, and for $\mfc$ this is Lemma A.1 of \cite{fixed}. It then suffices to prove the statement for the lattice $\hat {\mathcal D} = \mathcal D \cap (N^{-3} \Z^2)$. We will verify at the end of the proof that for $z\in\hat {\mathcal D}$ with $\eta\ge 1$ the claim follows from Lemma \ref{lem:apriori}. We proceed assuming that the claim is true for such $z$.

For $E$ such that $z_0 = E + i\in \hat {\mathcal D}$, define $\eta_k=1-kN^{-3}$ and $z_k=E + i \eta_k$. Fix $\sigma_1 <  \delta_1/100$ and $D>0$. Define \beq \Omega_k = \left\{  |m(z_k) - \mfc(z_k) | \le \frac{N^{\sigma_1}}{\sqrt{N\eta}} \right\} .\eeq

Now recall the definition of $\mfc$ and the self-consistent equation (\ref{eqn:sc1}) for $m$ derived in the proof of Lemma \ref{lem:bootstrap}. Subtracting these yields
\beq m-\mfc= \frac{1}{2N}\sum_{i=-N}^N  \frac{1}{-z-tm + v_i} - \frac{1}{-z-t\mfc + v_i}  + R \eeq
\beq =\frac{1}{2N}\sum_{i=-N}^N  \frac{t(m-\mfc)}{(-z-tm + v_i)(-z-t\mfc + v_i)} + R. \eeq
We obtain
\beq\label{eqn:sc} (m-\mfc)\left( 1 - \frac{t}{2N} \sum_{i=-N}^N  \frac{1}{(-z-tm + v_i)(-z-t\mfc + v_i)}  \right) = R. \eeq
By Lemma \ref{lem:apriori}, $\P\left(\Omega_0^c\right)\le N^{-D}$, as the second factor in (\ref{eqn:sc}) is bounded below and ${|R|  \prec (N\eta)^{-1/2}}$. We now consider $\Omega_1$. Because $m-\mfc$ is $2N^2$-Lipschitz on $\mathcal D$, we have
\beq \one(\Omega_0)|m(z_1) - \mfc(z_1)|  \le \frac{N^{\sigma_1}}{\sqrt{N\eta}}  + \frac{2}{N}.\eeq
Hence the hypothesis of Lemma \ref{lem:bootstrap} is verified for some $c>0$ when $\phi=\one(\Omega_0)$ and $z=z_1$. Hence for $z=z_1$, in (\ref{eqn:sc}) we have $|R|\prec (N\eta)^{-1/2}$ and that the second factor on the left side is bounded below for large enough $N$. To see this, write
\beq1 - \frac{t}{2N} \sum_{i=-N}^N  \frac{1}{(-z-tm + v_i)(-z-t\mfc + v_i)}\eeq
\beq= \left( 1 - \frac{t}{2N} \sum_{i=-N}^N  \frac{1}{(-z-t\mfc + v_i)^2}  \right) + \left( \frac{t}{2N} \sum_{i=-N}^N \frac{1}{(-z-t\mfc + v_i)^2}  -   \frac{1}{(-z-tm + v_i)(-z-t\mfc + v_i)} \right).\eeq
The first term is bounded below by Lemma \ref{lem:stability}, and the error term is $o(N^{-c})$ by Lemma \ref{lem:bootstrap}, because it equals
\beq \frac{t}{2N} \sum \frac{1}{(-z-t\mfc + v_i)}\left(\frac{1}{(-z-t\mfc + v_i)} - \frac{1}{(-z-tm+ v_i)} \right) \eeq
\beq =\frac{t}{2N} \sum \frac{1}{(-z-t\mfc + v_i)}\left(\frac{-t(m-\mfc)}{(-z-t\mfc + v_i)(-z-tm+ v_i)} \right) \eeq
\beq =  \frac{(\mfc - m)t^2}{2N} \sum \frac{1}{(-z-t\mfc + v_i)^2(-z-tm+ v_i)}.\eeq
We conclude 
\beq \P( \Omega_0 \cap \Omega_1^c) \le N^{-D}. \eeq

Now we may apply this reasoning sequentially for all $k$ such that $z_k \in \mathcal D$. Note that the $c>0$ used to verify the hypothesis of Lemma \ref{lem:bootstrap} can be chosen to be the same for each step, so this lemma needs to be invoked only once. The conclusion follows by noting that 
$\P\left( \cap_k\, \Omega_k \right)$ can be made larger than $1-N^{D_1}$ for any $D_1$ by taking $D$ large enough.

For $z\in\hat {\mathcal D}$ with $\eta\ge 1$, we can use the same argument with the bounds in Lemma \ref{lem:apriori}, which hold unconditionally, so there is no need for a bootstrapping argument. In particular, we do not need to use Lemma \ref{lem:gisum}, since we have the trivial bound $|g_i| \le \eta^{-1}\le 1$. \ep

\subsection{Strong law}\label{sec:deformedlawproof}

We now improve the bound Lemma \ref{lem:weaklaw} using fluctuation averaging. For any random variable $X$, let $Q_i X$ denote the conditional expectation of $X$ with respect to the $i$th column and $i$th row of $M_t$. For $I=(i,i+N)$, $J=(j, j+N)$, and any index set $\T$ containing pairs $(k, k+N)$, set
\beq G_{IJ} = \bma   G_{i,j} & G_{i, N+j} \\   G_{i+n, j} & G_{N+i, N + j} \ema, \quad M^{(\T)}_l = \bma -z - t m^{(\T)}(z) & v_l \\ v_l & -z - t m^{(\T)}(z) \ema . \eeq

\bel\label{lem:taylor} For a matrix $A$ with $A=B+R$, \beq A^{-1} = B^{-1}- B^{-1}RB^{-1} + B^{-1}RB^{-1} R A^{-1} .\eeq\eel
\bp Iterate $A^{-1}=B^{-1}-B^{-1}RA^{-1}$.
\ep

With $i\le N$ and $\T=(i,N+i)$, Schur's complement formula gives
\beq G_{II}^{-1}  = M_i^{(I)} + Q_i\left( G_{II}^{-1}  \right),\eeq
and then Lemma \ref{lem:taylor} gives
\beq\label{eqn:greens}  G_{II} = \left[M_i^{(I)}\right]^{-1} - \left[M_i^{(I)}\right]^{-1} Q_1\left( G_{II}^{-1}\right)\left[M_i^{(I)}\right]^{-1} +  \left[M_i^{(I)}\right]^{-1} Q_1\left( G_{II}^{-1}\right) \left[M_i^{(I)}\right]^{-1}Q_i\left( G_{II}^{-1}\right) G_{II}  .\eeq

Define the deterministic quantities 
\beq S_i^{-1}=  \bma -z - t \mfc (z) & v_i \\ v_i & -z - t \mfc (z) \ema,\quad  f_i = \frac{1}{-z - t\mfc + v_i}, \quad  f_{-i} = \frac{1}{-z - t\mfc - v_i} . \eeq

\bel\label{lem:S} We have
\beq \left \|  S_i -  \left[M_i^{(I)}\right]^{-1}  \right \| \prec \frac{1}{\sqrt{N\eta}},\eeq
and \beq \|S_i\| \le \frac{|f_i| + |f_{-i}|}{2}.\eeq
\eel 
\bp
Note that 
\beq S_i = \frac{1}{2} \bma f_i + f_{-i} & f_i - f_{-i}  \\  f_i - f_{-i} & f_i +  f_{-i} \ema, \quad S_i^{-1} =  M_i^{(I)} + \bma t\eps_1 & 0  \\  0 &  t\eps_2 \ema,\eeq
with $|\eps_k| \prec (N\eta)^{-1/2}$ by the weak law Lemma \ref{lem:weaklaw} and eigenvalue interlacing. We conclude using the same algebraic manipulations as in Section \ref{sec:greens}.

\ep

\bel\label{lem:fisum} For $z\in \mathcal D$ we have
\beq \frac{1}{N} \sum_i |f_i| + |f_{-i}| \le C \log N.\eeq
\eel
\bp
This is Lemma 7.5 in \cite{LY}.\ep

\bel\label{lem:resolve} For $I=(i,i+N)$, $J=(j, j+N)$,
\beq (G_{II})^{-1} = \left( G_{II}^{(J)} \right)^{-1} - (G_{II})^{-1}G_{IJ}(G_{JJ})^{-1}G_{JI} \left( G^{(J)}_{II}\right)^{-1}. \eeq\eel
\bp This is a consequence of Schur's formula. \ep

\bel\label{lem:matrixest} The following claims hold for $z\in \mathcal D$ and $\T, \S$ with $|\T|, |\S| \le \log N$, and such that $\T$ and $\S$ are composed of pairs $(k,k+N)$.
\benr 
\item We have \beq \left[M^{(\T)}_i\right]^{-1} = \frac{1}{2} \bma g_i  + g_{-i} &   g_i  -  g_{-i}  \\  g_i  - g_{-i}   & g_i  + g_{-i}\ema   + R ,\eeq
where $E$ is a matrix such that $\|R \| \prec (|g_i| + |g_{-i}|) (N\eta)^{-1}$ and this bound is uniform in the index $i$. Hence \beq \left\| \left[M^{\T}_i\right]^{-1}\right\| \prec  (|g_i| + |g_{-i}| )\le \frac{C}{t}. \eeq

\item We have \beq \left\|  \left[G_{II}^{(\T)}\right]^{-1} - S_i^{-1}   \right\| \prec \frac{t}{\sqrt{N\eta}} \eeq
and therefore \beq\label{eqn:constantestimate} \left\| \left[G_{II}^{(\T)}\right]^{-1}S_i\right \|  \prec C.\eeq

\item We have \beq \left\| Q_i\left(\left[G^{(\T)}_{II}\right]^{-1}\right)\right\| \prec \frac{t}{\sqrt{N\eta}},\quad \left\|  G^{(\T)}_{II} \right\|  \prec \frac{C}{t}.\eeq
\item For $I\neq J$, and $\T$ not containing $I$ or $J$, \beq\label{eqn:pair} \left\|\left[G^{(\T)}_{II}\right]^{-1}G^{(\T)}_{IJ}\right\| \prec \sqrt{\frac{t}{N\eta}}, \eeq \beq \left\|G^{(\T)}_{IJ} \right\| \prec \frac{\min(|g_i| +|g_{-i}|,  |g_j| + |\overline g_j|)}{\sqrt{N\eta}},\eeq
   and therefore
\beq\label{eqn:triple}  \left\|G^{(\T)}_{IJ} \left[G^{(\T)}_{JJ}\right]^{-1}G^{(\T)}_{JI}\right\| \prec \frac{\sqrt{t}}{N\eta}(|g_i| +|g_{-i}|). \eeq

\item We have \beq \left\| [G_{II}^{(\T)}]^{-1} [G_{II}^{(\S)}] \right\|\prec C,\eeq and hence (\ref{eqn:pair}) holds for any superscripts $\T$ and $\S$, and hence (\ref{eqn:triple}) holds for any combination of superscripts on the Green function elements.
\eenr
\eel
\bp \benr \item We have \beq M^{(\T)}_l = \bma -z - t m(z) +tr & v_l \\ v_l & -z - t m(z) + tr \ema.\eeq

Here we have used the Cauchy interlacing theorem at most $\log N$ times to split off the error $|r| \prec (N\eta)^{-1}$, and the $\log N$ factor can be absorbed by the stochastic domination. The first claim now follows from the same algebraic manipulations as in our discussion of the Green functions. The second claim follows from the stability bound Lemma \ref{lem:stability} and the weak law Lemma \ref{lem:weaklaw}.

\item The first claim follows from using the representation of $G_{II}^{-1}$ developed in our discussion of the Green function, only now applied in the same way to $\left[G^{(\T)}_{II}\right]^{-1}$. The errors are bounded by Cauchy interlacing to control the removed rows and columns and the large deviations inequalities Lemma \ref{lem:largedev} are used to control the fluctuations. The second claim follows from the first and the analogue of the previous part for $S_i$.
\item The first claim is just a special case of (ii). The second follows from the explicit computation of the entries of $G_{II}$ in Section \ref{sec:greens} and the stability bound $|g_i| \le C t^{-1}$. 
\item We first establish a Green function identity. We the Schur complement formula on the index set $I=(i,  i+N)$ as in Section \ref{sec:greens}, except now we concentrate on the upper off-diagonal block. Write $G_{II^c}$ for the sub-matrix of $G$ whose rows are taken from the indices in $I$ and columns from the indices in $I^c$. We obtain 
\beq G_{II^c} = -G_{II} v^\trans G^{(I)}, \eeq
which implies
\beq G_{II} ^{-1} G_{IJ} = \sum_{K}^{(I)}  H_{IK} G_{KJ} ,\quad  H_{IK} = \bma  0 & \sqrt{t}w_{ik}\\ \sqrt{t}w_{ki} & 0  \ema .\eeq
We bound just the first entry of $G_{II} G_{IJ}$, as the rest are similar. This entry is 
\beq (G_{II}^{-1} G_{IJ})_{11}=  \sum_k \sqrt{t}w_{ik} G^{(I)}_{k+n, j}.\eeq 
By the large deviations estimates Lemma \ref{lem:largedev} and Ward's identity, we have 
\beq \left| \sum_k \sqrt{t}w_{ik} G^{(I)}_{k+n, j}\right| \prec \sqrt{\frac{t}{N}}\sqrt{\sum_k |G^{(I)}_{k+n, j}|^2 } 
\le\sqrt{\frac{t}{N}} \sqrt{\frac{\Im G_{kk}^{(I)}}{\eta}} \prec \sqrt{\frac{t}{N\eta}}.\eeq
In the last step, we used the bound $\Im G_{kk}^{(I)} \le \Im m^{(I)} \le \log(N)$ and absorbed the $\log N$ factor into the stochastic domination.
For the second claim, we have a similar identity, obtained in the same way using the other off-diagonal block in the Schur complement formula.
\beq G_{IJ} = \left(\sum^{(J)}_K  G_{IK}^{(J)}  H_{KJ}\right) G_{JJ}  \eeq 
We can expand this using the first identity to obtain
\beq G_{IJ} = G_{II} \left( \sum_{M,K}^{(I, J)}  H_{IM} G^{(IJ)}_{MK}  H_{KJ} - H_{IJ}\right) G_{JJ}^{(I)}  \eeq 
As shown in the work in Section \ref{sec:greens}, up to an $O(1)$ factor, we have 
\beq \| G_{II} \| \le |f_i| + | f_{-i}|, \quad \| G_{JJ} \| \le |f_j| + |f_{-j}|. \eeq
The middle factor has a norm that is stochastically dominated by
\beq t \sqrt{\frac{\log N}{N\eta}}  + \frac{\sqrt{t}}{\sqrt{N}} ,\eeq
where we used the large deviations bounds Lemma \ref{lem:largedev} and the Ward identity as before. The lemma follows from using ${|g_i|, |g_{-i}| \le \min(t^{-1}, \eta^{-1})}$. The third claim follows from the first two. Clearly we could repeat this argument for any $\T$ satisfying the given hypotheses.

\item First, we show this implies the modification of (\ref{eqn:pair}). We want to bound $[G^{(\T)}_{II}]^{-1} G^{(\S)}_{IJ}$. Write this as 
\beq [G^{(\T)}_{II}]^{-1}  [G^{(\S)}_{II}]  [G^{(\S)}_{II}]^{-1}  G^{(\S)}_{IJ} .\eeq
The final two terms are bounded in norm by the previous part, so we need to bound 
\beq\left\|[G^{(\T)}_{II}]^{-1}  [G^{(\S)}_{II}]\right\|\le C\eeq as claimed. This follows from the work in Section \ref{sec:greens} by estimating
\beq \left\| [G^{(\S)}_{II}]^{-1} - [G^{(\T)}_{II}]^{-1}\right\| \prec \sqrt{ \frac{1}{N\eta}}\eeq
using the explicit representation there.
\eenr
\ep
We record an elementary fact for later use.
\bel\label{lem:tracenorm}
For any symmetric matrix $M$ and integer $r>0$, \beq \| M \|^{2r} \le  \operatorname{Tr} M^{2r}.\eeq
\eel The following arguments are based on the proof of Lemma 7.15 in \cite{LY}. 
\bel\label{lem:moments} For $z\in \mathcal D$ and $p$ even,
\beq \E\left\| \frac{1}{N}\sum S_i Q_i\left( G_{II}^{-1}\right) S_i \right\|^{p} \prec \frac{1}{(N\eta)^{p}}.\eeq
\eel 
\bp
We proceed as in the proof of Lemma 4.7 in Appendix B of \cite{scgeneral}, by computing moments. We first consider the case $p=2$. Recall that if $Y$ is a random variable independent of the $i$th row and $i$th column, then $\E Q_i(X)Y=\E Q_i(XY)=0$. We have 

\beq \E \left\| \frac{1}{N}\sum_{i=1}^N S_i Q_i\left( G_{II}^{-1}\right) S_i\right\|^2 \le \operatorname{Tr} \frac{1}{N^2}\sum_{i=1}^N\E\left[ S_i Q_i\left( G_{II}^{-1}\right)S_i \right]^2\eeq \beq + \operatorname{Tr}  \frac{1}{N^2}\sum_{i\neq j}^N \E \left[S_i Q_i\left( G_{II}^{-1}\right)S_i S_j Q_j\left( G_{JJ}^{-1}\right)S_j \right]= A_1 + A_2. \eeq
By Lemma \ref{lem:matrixest}, and the fact that, up to a constant, the norm of a matrix bounds its trace, we have
\beq |A_1| \prec C \frac{1}{N\eta} \frac{1}{N^2} \sum_{i=1}^N (|f_i| + |f_{-i}|)^2  \le \left( \frac{1}{N\eta} \right)^2. \eeq
Here we used
\beq \frac{1}{N}\sum_{i=1}^N |f_{i}|^2 + |f_{-i}|^2 \le \frac{1}{\eta},\eeq
which comes from $|f_i|\le \eta^{-1}$. 

For $A_2$, we note that $S_i$ is deterministic, so we may move it inside the $Q_i$. 
\beq A_2 = \frac{1}{N^2}\sum_{i\neq j}^N \E Q_i\left(S_i G_{II}^{-1} S_i\right) Q_j\left(S_j G_{JJ}^{-1}S_j \right).\eeq
By Lemma \ref{lem:resolve}, we can write 
\beq Q_i\left( S_i \left[ \left( G_{II}^{(J)} \right)^{-1} - (G_{II})^{-1}G_{IJ}(G_{JJ})^{-1}G_{JI} \left( G^{(J)}_{II}\right)^{-1} \right] S_i \right)  \eeq 
Recall that if $Y$ is a random variable independent of the $j$th row and $j$th column, then $\E Q_j(X)Y=\E Q_j(XY)=0$. So the first term will cancel when multiplied against the $Q_j(\cdot)$ term. We are left with multiplying two terms of the following form.
\beq  Q_i\left( S_i (G_{II})^{-1}G_{IJ}(G_{JJ})^{-1}G_{JI} \left( G^{(J)}_{II}\right)^{-1} S_i \right) \eeq
Using Lemma \ref{lem:matrixest} again, we see  that this is stochastically dominated by $(N\eta)^{-1}(|g_i| + |g_{-i}|)$. Specifically, we use (\ref{eqn:triple}) on the middle three terms and (\ref{eqn:constantestimate}) on the outside pairs. Then
\beq A_2 \le\frac{1}{(N\eta)^2} \frac{1}{N^2} \sum_{i\neq j }  (|f_i| + |f_{-i}|) (|f_j| + |f_{-j}|)\eeq
Now we sum over $i$ and $j$ separately, picking up $\log$ terms by Lemma \ref{lem:fisum}, and obtain the bound $A_2 \prec (N\eta)^{-2}$.

We now discuss even $p>1$. Our strategy, as in \cite{LY}, is to adapt the proof of Theorem 4.7 in \cite{scgeneral} to the deformed case. 
Applying Lemma \ref{lem:tracenorm} to $M=\left\| \sum S_i Q_i\left( G_{II}^{-1}\right) S_i \right\|^{p}$ yields, using the notation in \cite{LY},
\beq \E\left\| \frac{1}{N}\sum_{i=1}^N S_i Q_i\left( G_{II}^{-1}\right) S_i \right\|^{p} \le \frac{1}{N^{p}} \E\left[ \operatorname{Tr} \sum_{k_1,\cdots,k_{2p}} \prod  S_{k_i} Q_{k_i}\left( G_{{K_i}{K_i}}^{-1}\right) S_{k_i}  \right]\eeq
\beq =\operatorname{Tr} \frac{1}{N^{p}} \sum_{\Gamma\in\mathcal{P}_{p} } \sum_{i_1,\dots, i_r} \one_{\{\Gamma = \Gamma(i) \}} \E\left[ S_{i_1} Q_{i_1}\left( G_{{I_1}{I_1}}^{-1}\right) S_{i_1}\dots S_{i_p} Q_{i_p}\left( G_{{I_p}{I_p}}^{-1}\right) S_{i_p} \right].\eeq 

We now apply the algorithm in \cite{scgeneral} to decompose each term, using the analogous substitutions (which also hold with any superscript $\T$ added to the Green function matrices)
\beq G_{IJ} \rightarrow G_{IJ}^{(K)} + G_{IK} G^{-1}_{KK} G_{KJ},\eeq
\beq (G_{II})^{-1} \rightarrow \left( G_{II}^{(J)} \right)^{-1} - (G_{II})^{-1}G_{IJ}(G_{JJ})^{-1}G_{JI} \left( G^{(J)}_{II}\right)^{-1}.\eeq
We see that these substitutions create non-commutative polynomials in the Green functions with the properties that inverted Green function matrices alternate with non-inverted ones, and that the lower indices pair across adjacent terms. The algorithm yields binary strings $\sigma_k$ and an expansion into monomials. 
\beq \E\left[ S_{i_1} Q_{i_1}\left( G_{{I_1}{I_1}}^{-1}\right) S_{i_1}\dots S_{i_p} Q_{i_p}\left( G_{{I_p}{I_p}}^{-1}\right) S_{i_p} \right] = \sum_{\sigma_1, \dots, \sigma_p} \E[S_{i_1} Q_{i_1}\left( F_{i_1}\right)_{\sigma_1} S_{i_1} \dots S_{i_p} Q_{i_p}\left( F_{i_p}\right)_{\sigma_p} S_{i_p} ]\eeq
Set $\Phi = (N\eta)^{-1/2}$. We now claim it suffices to establish the key bound 
\beq \left\|[\E S_{i_1} Q_{i_1}\left( F_{i_1}\right)_{\sigma_1} S_{i_1} \dots S_{i_p} Q_{i_p}\left( F_{i_p}\right)_{\sigma_p} S_{i_p} ]\right\| \prec C \Psi^{p+s} \prod_{i_k} (|f_{i_k}| + |f_{-i_k}|) \eeq
where $s$ is the number of lone labels. This is proved in the Lemma \ref{lem:key}, following this proof. Assuming this lemma, we will conclude the proof.

It follows from Lemma \ref{lem:key} that for a partition $\Gamma$ with $l=|\Gamma|$, 
\beq \left\| \frac{1}{N^{p}} \sum_{i_1, \dots, i_p} \one_{\Gamma = \Gamma(i) } \E[\dots ]\right\| \prec \frac{\Psi^{p+s}}{N^p} \sum_{k_1=1}^N \cdots \sum_{k_l =1}^N (|f_{k_1}| + | f_{-k_1}|)^{d_1}\cdots (|f_{k_l}| + |\overline f_{k_l}|)^{d_l},\eeq
where the $d_l$ are the sizes of the blocks in $\Gamma$. Then $|f_i|\le \eta^{-1}$ implies the bound
\beq \le \Psi^{p+s}\frac{1}{N^p} \frac{1}{\eta^{p-l}} \sum_{k_1=1}^N \cdots \sum_{k_l =1}^N (|f_{k_1}| + | f_{-k_1}|)\cdots (|f_{k_l}| + |f_{-k_l}|), \eeq
\beq \le\Psi^{p+s} \frac{1}{N^{p-l}} \frac{1}{\eta^{p-l}} \log(N)^l \prec \Psi^{2p+s}\Psi^{2p-2l}. \eeq
Using $p+s+(2p-2l)\ge 2p$ we have
\beq \E\left\| \frac{1}{N}\sum S_i Q_i\left( G_{II}^{-1}\right)S_i\right\|^{p} \prec \sum_{\Gamma \in \mathcal{P}_p} \Psi^{2p} \le (Cp)^{Cp} \Psi^{2p}\prec \Psi^{2p}. \eeq
This concludes the proof.
\ep

\bel\label{lem:key}
\beq\left\| \E [S_{i_1} Q_{i_1}\left( F_{i_1}\right)_{\sigma_1} S_{i_1} \dots S_{i_p} Q_{i_p}\left( F_{i_p}\right)_{\sigma_p} S_{i_p} ]\right \| \prec C \Psi^{p+s} \prod_{i_k} (|f_{i_k}| +| f_{-i_k}|) \eeq
\eel

\bp
We claim that if $b(\sigma_k)$ is the number of ones appearing in $\sigma_k$, then 
\beq\label{eqn:ones} \|S_{i_k} Q_{i_k}\left( F_{i_k}\right)_{\sigma_k} S_{i_k}\| \prec \Psi^{b(\sigma_k) + 1}(|f_{i_k}| +| f_{-i_k}| ).\eeq
If $\sigma=0$, then we have a single term of the form $S_{i_k}Q_{i_k}\left((G_{II}^{(T)})^{-1}\right)S_{i_k}$, and using Lemma \ref{lem:matrixest} we have
$\|S_{i_k}Q_{i_k}\left((G_{II}^{(T)})^{-1}\right)S_{i_k}\| \le \Psi(|f_{i_k}| +| f_{-i_k}|) . $
If $\sigma>1$, we can proceed similarly, except each interior pair of diagonal and off-diagonal terms gains a factor of $\Psi$ by (\ref{eqn:pair}), and there are $\sigma$ of these.

Now we have various cases.

\benr
\item One of the monomials $(F_{i_k})_{\sigma_k}$ is not maximally expanded. Then $(F_{i_k})_{\sigma_k}$ contains $\ge 2p$ off-diagonal Green function entries, and $b(\sigma_k)\ge 2p-1 \ge 2r +s$. The result follows from the previous discussion.
\item Every monomial is maximally expanded, and for every lone label $a$ there is a label $b\in \{1,\dots, p\}\setminus \{a\}$ such that the monomial $(F_{i_b})$ contains an off-diagonal resolvent entry with lower index $I_a$. Then $\sum b(\sigma_k)\ge s$ and again we are finished.
\item Every monomial is maximally expanded, and there is a lone label without a matching lower index in another monomial, say $i_1$. 
Then, since the $S_i$ are deterministic and the other terms are maximally expanded, the other ${S_{i_k} Q_{i_k}(F_{i_k})_{\sigma_k} S_{i_k}}$ are independent of the $i_1$ row and column. Then we can write 
\begin{multline}\E \left[S_{i_1} Q_{i_1}\left( F_{i_1}\right)_{\sigma_1} S_{i_1} \dots S_{i_p} Q_{i_p}\left( F_{i_p}\right)_{\sigma_p} S_{i_p} \right] \\= \E \left[Q_{i_1} \left(S_{i_1} \left( F_{i_1}\right)_{\sigma_1} S_{i_1} \dots S_{i_p} Q_{i_p}\left( F_{i_p}\right)_{\sigma_p} S_{i_p}\right) \right] = 0.\end{multline}
\eenr
\ep

\bel\label{lem:averaging} For $z\in \mathcal D$, we have entrywise 
\beq \frac{1}{N}\sum S_i Q_i\left( G_{II}^{-1}\right)S_i \prec \frac{1}{N\eta}.\eeq

\eel
\bp
This follows from Lemma \ref{lem:moments} and Markov's inequality.
\ep

\noindent {\bf Proof of Theorem \ref{thm:deformedlaw}.}\ \ 
We consider each term in (\ref{eqn:greens}), average over $N$, and take trace. For the first term we get
\beq g_i + g_{-i} + \mathrm{error}, \eeq
where \beq |\mathrm{error}| \le  \frac{C}{N\eta} \frac{1}{2N} \left(\sum_{i=1}^N |g_i| + |g_{-i}| \right) \le \frac{C\log N}{N\eta} .\eeq

It remains to show the contributions from averaging the last two terms are negligible. The second term is dealt with  by using Lemma \ref{lem:S} to make the replacement 
\beq \left[M_i^{(I)}\right]^{-1} Q_1\left( G_{II}^{-1}\right)\left[M_i^{(I)}\right]^{-1} \rightarrow S_i  Q_1\left( G_{II}^{-1}\right)S_i \eeq
and then Lemma \ref{lem:averaging} to bound the averaged error from this replacement. The largest terms in the error have the form
\beq \left\| \left(S_i  - \left[M_i^{(I)}\right]^{-1}\right)Q_1\left( G_{II}^{-1}\right)\left[M_i^{(I)}\right]^{-1} \right\| \prec \frac{1}{\sqrt{N\eta}} \frac{t}{\sqrt{N\eta}} \frac{1}{t} = \frac{1}{N\eta}.\eeq
By Lemma \ref{lem:matrixest}, the third term is 
\beq \left[ \frac{1}{2}\bma g_i  + g_{-i} &   g_i  -  g_{-i}  \\  g_i  - g_{-i}   & g_i  + g_{-i}\ema + O\left( \frac{|g_i| + |g_{-i}|}{N\eta} \right) \right]Q_1\left( G_{II}^{-1}\right) \left[M_i^{(I)}\right]^{-1}Q_i\left( G_{II}^{-1}\right) G_{II}  \eeq
Again by Lemma \ref{lem:matrixest}  \beq  \left \| Q_1\left( G_{II}^{-1}\right) \left[M_i^{(I)}\right]^{-1} Q_i\left( G_{II}^{-1}\right) G_{II} \right \| \le C  \frac{t}{\sqrt{N\eta}}  \frac{1}{t} \frac{t}{\sqrt{N\eta}} \frac{1}{t} \le \frac{C}{N \eta},\eeq
and \beq  \left\| \frac{1}{2}\bma g_i  + g_{-i} &   g_i  -  g_{-i}  \\  g_i  - g_{-i}   & g_i  + g_{-i}\ema + O\left( \frac{|g_i| + |g_{-i}|}{N\eta} \right) \right\| \le C (|g_i| + |g_{-i}|).\eeq 
Hence the average over the third term is negligible. Note the average over the $|g_i|\pm |g_{-i}|$ picks up a negligible $\log N$ factor, as above.

Finally, we have 
\beq m=\frac{1}{2N} \sum_{-N}^{N} \frac{1}{-z-tm + v_i} + R,\eeq
where $|R| \prec (N\eta)^{-1}$, and we can repeat the proof of Lemma \ref{lem:weaklaw} with the improved error $(N\eta)^{-1}$ in place of $(N\eta)^{-1/2}$ to conclude.
\ep

\section{Removal of time evolution}\label{sec:removal}

In this section we show how to complete the proof of universality given the main homogenization result Theorem \ref{thm:homomain}. In Subsection \ref{sll} we prove a local law for sparse matrices, which is necessary in the following subsections. In Subsection \ref{sec:newdynamics} we prove short time universality for sparse random matrices. Finally, in Subsection \ref{gfc} we show how to remove the time evolution through a Green function comparison argument. 

\subsection{Sparse local law}\label{sll}

For clarity, in this subsection we consider just the case of sparse ensembles where the variances are equal, $s_{ij} = N^{-1}$, but the case of a general doubly stochastic variance matrix satisfying the conditions in Definition \ref{def:sparse} can be handled with minor modifications. The key point is that in each of these cases the limiting spectral distribution is a semicircle. More general variance matrices (which give rise to new limit distributions), along with correlated entries, are considered in Section \ref{a:correlated}.

We prove a weak local law for the singular values of sparse matrices. The symmetric case was considered in \cite{EKYY13}.  We recall that our model is $M = B + f \ket{w}\bra{w}$, using the notation of Section \ref{sec:main}. The next lemma implies it is enough to prove such a law for $B$. Let the singular values of $M$ be $(\mu_i)_{i=1}^N$ and the singular values of $B$ be $(\lambda_i)_{i=1}^N$.

\bel The singular values of $M$ and $B$ are interlaced,
\beq \lambda_{j+ 1}  \ge \mu_j \ge \lambda_{j- 1}.\eeq
\eel 
\bp
Note that $M$ is a rank $1$ perturbation of $B$. The result follows from Weyl's inequality and Majorization for singular values.
\ep

\bel\label{lem:sparsedecay}
Letting $m_{ij}$ denote the entries of $M$, we have 
\beq |m_{ij}| \prec \frac{1}{q}.\eeq
\eel
\bp
This follows from Markov's inequality.
\ep

We also require a slight modification of Lemma 3.8 from \cite{EKYY13}. 

\bel\label{lem:sparsedev}
Let $a_1,\dots, a_N$ be centered and independent random variables satisfying
\beq \E|a_i|^p \le \frac{C^p}{N q^{p-2}}\eeq
for all $p$. Then for any $A_i\in \C$ and $B_{ij}\in \C$,
\beq \left| \sum_{i=1}^N A_i a_i \right| \prec \frac{\max_i |A_i|}{q} + \left(\frac{1}{N}\sum_{i=1}^N |A_i|^2 \right)^{1/2},\eeq
\beq \left| \sum_{i=1}^N \overline{a_i}B_{ii} a_i  - \sum_{i=1}^N \sigma^2_i B_{ii} \right| \prec \frac{B_d}{q},\eeq
\beq \left|  \sum_{1\le i \neq j \le N} \overline{a_i}B_{ij}a_j  \right| \prec  \frac{B_o}{q} + \left(\frac{1}{N^2} \sum_{i\neq j} |B_{ij}|^2 \right)^{1/2} ,\eeq
where $\sigma_i^2$ is the variance of $a_i$ and 
\beq B_d = \max_i |B_{ii}|, \quad B_o = \max_{i\neq j} |B_{ij}| .\eeq

Further, if $a_1,\dots, a_N$ and $b_1,\dots, b_N$ are independent random variables satisfying the above moment condition, then for $B_{ij}\in\C$ we have
\beq \left| \sum_{i,j=1}^N a_i B_{ij} b_j \right| \prec \left[ \frac{B_d}{q^2} + \frac{B_o}{q}  \left( \frac{1}{N^2} \sum_{i\neq j} |B_{ij}|^2 \right)^{1/2} \right] \eeq

\eel

Let $K = \bma 0 & B  \\ B^\trans & 0 \ema$ be the symmetric $2N\times 2N$ block matrix formed from $B$ and $G_{ij}(z)$ be entries of the Green function of $K$. Let $m(z)$ be the Stieltjes transform of $K$. Define
\beq \Lambda_o = \max_{i\neq j} |G_{ij}|,  \quad \Lambda_d = \max_i |G_{ii} - \msc|, \quad \Lambda  = | m - \msc|.\eeq

Repeating the Green function calculations for the deformed case, we have

\beq\label{eqn:sparseg} G_{ii} = g\left( \frac{1+B_i}{1+K_i} \right),\eeq
with
\beq g = \frac{1}{-z-m}, \quad K_i = E_3 + E_4 + E_3 E_4 -g^2 E_1,\quad B_i = \frac{-([A_{22} - m^{(\T)}] - r )}{-z-m},\eeq
\beq E_1 = \left(\frac{A_{11} - A_{22} }{2} \right)^2, \quad E_2 = \frac{A_{11} + A_{22}}{2} - m^{(\T)}, \quad E_1 = g(h_{ii} - A_{12} + E_2 +r), \quad r = m - m^{(\T)}. \eeq

For any $\delta>0$, set

\beq\mathcal D_\delta =\{ z = E + i\eta \colon E\in (-1,1), N^{\delta} \le N\eta \le 10N \}.\eeq
We will proceed largely as in the proof of the deformed law. The key difference is that we have better stability for $\msc$ (see Lemma 6.2 in \cite{EYbook}), so we can prove the local law for a larger spectral domain.

\bel\label{lem:sparsebootstrap}

Suppose $z\in \mathcal D_\delta$. Let $\phi$ be the indicator function of some event, which may depend on $z$. If  $\phi (\Lambda_o + \Lambda_d) \prec N^{-c}$ for some $c>0$, then

\beq \phi \max_i |K_i| {\prec}\left( q^{-1} + (N\eta)^{-1/2} \right), \quad  \phi \max_i |B_i|  {\prec}\left( q^{-1} + (N\eta)^{-1/2} \right), \eeq \beq \left(1 - \frac{1}{(-z-m)(-z-\msc)} \right) =O_\prec(N^{-c}) .\eeq
\eel 

\bp
For concreteness, we consider $B_1$, but our bounds will be uniform in $i$. By the stability bound for $\msc$ and the hypothesis on $m$, we have $\frac{1}{-z-m} \le C$. By Cauchy's interlacing lemma, $r\le C(N\eta)^{-1}$. Using Lemma \ref{lem:sparsedev}, reasoning as in the proof of the deformed law, we have
\beq |A_{22} - m^{(\T)}| \prec \frac{\Lambda_o}{q} + \frac{1}{\sqrt{N\eta}}. \eeq
Combining these completes the proof for the $B_i$.

For the $K_i$, it remains to bound $E_3$ and $E_4$. As these are similar to what was done before, we just sketch the proof for $E_3$. 
\beq |E_3| \le | h_{11} - A_{12} - E_2 +r|\eeq
The $r$ term was already bounded. Large deviations arguments using Lemma \ref{lem:sparsedev} suffice to bound $A_{12}$ and $E_2$, and $h_{11}$ is bounded using Lemma \ref{lem:sparsedecay}. Combining these completes the proof.

For the final bound, we write
\beq 1 - \frac{1}{(-z-m)(-z-\msc)}  = \left(1 - \frac{1}{(-z-\msc)^2}\right) + \left(\frac{1}{(-z-\msc)^2} - \frac{1}{(-z-m)(-z-\msc)}\right).\eeq
The first term equals $1-\msc^2$ and by Lemma 6.2 of \cite{EYbook} it is bounded above by a constant in $\mathcal D_\delta$. The second term is $O(N^{-c})$, which follows from the hypotheses and the bounds in the aforementioned lemma.
\ep

\bel\label{lem:lambdacontrol} For $z\in \mathcal D_\delta$, Let $\phi$ be the indicator function of some event, which may depend on $z$. If $\phi (\Lambda_o + \Lambda_d) \prec N^{-c}$ for some $c>0$, then
\beq \Lambda_d \prec \Lambda + \frac{1}{q} + \frac{1}{\sqrt{N\eta}}, \quad \Lambda_o \prec \frac{1}{q} + \frac{1}{\sqrt{N\eta}}. \eeq 
\eel
\bp
The first claim is proved the same way as in display (3.39) in \cite{EKYY13}. We use the explicit expression (\ref{eqn:sparseg}) for $G_{ii}$ above to compute
\beq G_{ii} - G_{jj} \le C\left( \frac{1}{\sqrt{N\eta}} + \frac{1}{q}  \right).\eeq
The claim follows by fixing $i$ and averaging over $j$. 

The second claim is proved as in Lemma 3.13 of \cite{EKYY13}.  We use \beq G_{IJ} = G_{II} \left( \sum_{M,K}^{(I, J)}  H_{IM} G^{(J)}_{MK}  H_{KJ} - H_{IJ}\right) G_{JJ}^{(I)}.\eeq 
By hypothesis we find
\beq \| G_{II} \| \le C,\quad  \| G_{JJ} \| \le C,\eeq
so using Lemma \ref{lem:sparsedev} on the individual entries of the matrix expression,
\beq \left\| G_{II} \left( \sum_{M,K}^{(I, J)}  H_{IM} G^{(J)}_{MK}  H_{KJ} - H_{IJ}\right) G_{JJ}\right\| \prec C\left( \frac{1}{q} + \frac{\Lambda_o}{q}  + \left( \frac{1}{N^2} \sum_{k,l}^{(IJ)} \left|G^{(IJ)}_{kl}\right|^2 \right)^{1/2} \right).\eeq  
By Ward's identity,
\beq \frac{1}{N^2} \sum_{k,l}^{(IJ)} \left|G^{(IJ)}_{kl}\right|^2 = \frac{1}{N^2 \eta} \sum^{(IJ)}_k \Im G_{kk}^{(ij)} \le \frac{\Im m}{N\eta} + \frac{C \Lambda^2_o}{N\eta}, \eeq
where the last inequality follows from using 
\beq |G_{ij}|\le C,  \quad c \le |G_{ii}|\le C, \quad G_{ij} = G_{ij}^{(k)} + \frac{G_{ik} G_{kj}}{G_{kk}}. \eeq
repeatedly. Taking the maximum over $i\neq j$ gives
\beq \Lambda_o \prec \frac{C}{q} + o(1) \Lambda_o + \sqrt{\frac{\Im m}{ N\eta}},\eeq
which implies the claim.
\ep
\bel\label{lem:sparseapriori}
If $\eta \ge 2$, then, 
\beq \max_i |K_i| {\prec}\left( q^{-1} + (N\eta)^{-1/2} \right), \quad  \max_i |B_i| {\prec}\left( q^{-1} + (N\eta)^{-1/2} \right),\eeq \beq \left(1 - \frac{1}{(-z-m)(-z-\msc)} \right)\ge c.\eeq
\eel 
\bp
The proof is similar to the above using the trivial estimates for $\eta \ge 2$ as in Lemma \ref{lem:apriori}.
\ep

\bel\label{lem:apriorisc} If $\eta \ge 2$,
\beq \Lambda_d(z) + \Lambda_o(z) \prec \frac{1}{\sqrt N} + \frac{1}{q}.\eeq
\eel
\bp
The bound on $\Lambda_o$ follows from the same calculation as in Lemma \ref{lem:lambdacontrol}, where we now use Lemma \ref{lem:sparseapriori} to bound the error terms and the trivial estimates $|G_{ij}| \le \eta^{-1}$ to bound the Green function entries. For $\Lambda_d$ we estimate using (\ref{eqn:sparseg})
\beq \Lambda_d = | \msc  - G_{ii} | =\left| \msc -  \frac{1}{-z-m} (1+ R)  \right| \le  \left|   \frac{1}{z+m} - \frac{1}{z+\msc} \right| + C|R| \eeq
\beq  \le \left|  \frac{m-\msc}{(m+z)(\msc + z) }\right| + C|R|  \le \frac{\Lambda_d}{3/2} + C|R|, \eeq
where we used 
\beq |z + \msc| = |\msc|^{-1} \ge 2, \quad |m - \msc | \le 1.\eeq
The conclusion follows by combining the $\Lambda_d$ terms on the left and using Lemma \ref{lem:sparseapriori} to bound $R$.
 \ep
Let \beq H = \bma 0 & M  \\ M^\trans & 0 \ema\eeq be the symmetric $2N\times 2N$ block matrix formed from $M$, and define $\tilde m$ to be the Stieltjes transform of $H$.
\bel\label{lem:sparseweaklaw}
Uniformly for $z\in \mathcal D_\delta$, we have $|\tilde m -\msc | \prec q^{-1/2} + (N\eta)^{-1/2}$.
\eel
\bp
As noted above, it is enough to establish the theorem for $m$, the Stieltjes transform of $K$. Define, following the proof of Lemma \ref{lem:weaklaw}, the lattice $\hat {\mathcal D_\delta} = \mathcal D_\delta \cap (N^{-3} \Z^2)$. We have already shown in Lemma \ref{lem:apriorisc} that the claim holds for $z\in\hat {\mathcal D_\delta}$ with $\eta\ge 2$. As in the proof of the deformed weak law, it suffices to prove the result holds uniformly for elements of the lattice $\hat {\mathcal D_\delta}$ with $\eta < 2$. Define $n_k=2 - k N^{-3}$ and $z_k = E + i\eta_k$. Fix $\sigma >0$ and $D>0$, and define
\beq\Omega_k = \left\{  \Lambda_o(z_k) + \Lambda_d(z_k) \le \frac{N^\sigma}{\sqrt{N\eta}} + \frac{1}{q} \right\}.\eeq
Note that $\Lambda_d \ge \Lambda$, so $\Lambda(z_k)  \le \frac{N^\sigma}{\sqrt{N\eta}} + \frac{1}{q}$ on $\Omega_k$.

It is well known that $\msc$ satisfies a self consistent equation
\beq \msc = \frac{1}{-z-\msc}.\eeq
Using Lemma \ref{lem:sparseapriori}, we may Taylor expand (\ref{eqn:sparseg}) to find 
\beq G_{ii} = \frac{1}{-z-m}(1 +R'),\eeq
where $|R'|\prec q^{-1} + (N\eta)^{-1/2}$. By the stability estimate $\Im \msc \ge c$ and working on the set $\Omega_0$ to control $\Lambda$, we have 
\beq m = \frac{1}{-z-m} + R\eeq
where $|R|\prec q^{-1} + (N\eta)^{-1/2}$.
Subtracting the two self consistent equations yields 
\beq (m-\msc)\left(1 - \frac{1}{(-z-m)(-z-\msc)} \right) = R.\eeq
By Lemma \ref{lem:apriorisc}, $\P(\Omega_0^c)\le N^{-D}$. We now consider $\Omega_1$. Because $m-G_{ii}$ is $2N^2$-Lipschitz on $\mathcal D$, we have \beq \one(\Omega_0)| m(z_1) - G_{ii}(z_1)| \le \frac{N^\sigma}{\sqrt{N\eta}} + \frac{2}{N} + \frac{1}{q}.\eeq
This shows that $\Lambda_d(z_1) \le N^{-c}$ for some $c>0$, and similar reasoning applies to $\Lambda_o$. Then, by Lemma \ref{lem:sparsebootstrap}, the coefficient of $(m-\msc)$ in the self-consistent equation is bounded below and $|R| \prec q^{-1} + (N\eta)^{-1/2}$. Hence $\Lambda(z_1) \prec (N\eta)^{-1/2} + q^{-1}$. We conclude by Lemma \ref{lem:lambdacontrol} that
\beq \P(\Omega_0 \cap \Omega_1^c) \le N^{-D}.\eeq

We may apply this reasoning sequentially for all $k$ such that $z_k\in \mathcal D$. The conclusion follows by noting that $\P(\cap_k \Omega_k)$ can be made larger than $1-N^{D_1}$ for any $D_1$ by taking $D$ large enough. \ep

\subsection{Short time universality for sparse matrices}\label{sec:newdynamics}

 Let $M_N$ be a sparse matrix ensemble and define $H$ by forming a $2N\times 2N$ symmetric block matrix as in (\ref{eqn:ht}):
 \beq H= \bma 0 & M_N \\ M_N^\trans & 0 \ema.\eeq
  
We now define the perturbed matrix $H_t$ that we show universality for. Because we must accommodate the possibly unequal variance structure, we cannot simply add a Gaussian matrix. Instead, we evolve the nonzero entries of $H$ according to the following Ornstein-Uhlenbeck dynamics from \cite{QUE}. 
\beq\label{eqn:OU} d(h_{ij}(t) -f ) = \frac{dB_{ij}(t)}{\sqrt{N}} - \frac{1}{2Ns_{ij}} (h_{ij}(t) -f)\,dt\eeq
Here the $B_{ij}$ are independent Brownian motions. The entries of the evolved matrix $H_t$ satisfy
\beq h_{ij}(t) = f + \exp\left( - \frac{t}{2Ns_{ij} }\right)(h_{ij}(0) - f)  + \frac{1}{\sqrt{N}} \int_0^t \exp\left( -\frac{s-t}{2Ns_{ij}} \right)\, dB_{ij}(s).\eeq
We choose this dynamics because it preserves the mean and variance of the initial entries, and because the resulting entries are Gaussian divisible. With $r = \min\left\{ N s_{ij} \right\}$ and $G$ a $2N\times 2 N$ symmetrized version of a $N\times N$ ensemble of independent standard Gaussian variables,
\beq H_t \overset{d}{=}  H^{(1)}_t +  \sqrt{\frac{r(1-\exp(-t/r))}{N}}G,\eeq
where
\beq \left( H^{(1)}_t \right)_{ij}  \overset{d}{=}  f+ e^{ - \frac{t}{2Ns_{ij} }} (h_{ij}(0) - f) + \sqrt{Ns_{ij} \left(1 - e^{ - \frac{t}{Ns_{ij} }} \right) - r (1 - e^{-t/r}) } \frac{\tilde B_{ij}(t)}{\sqrt{N}},\eeq
and the $\tilde B_{ij}$ are a family of symmetric, independent Brownian motions. Note that \beq \sqrt{t} \asymp \sqrt{r(1-\exp(-t/r))}\eeq because $r$ is bounded below.

Before invoking Theorem \ref{thm:homomain}, we prove a lemma that assists in showing that sparse matrices are $(g,G)$-regular.
\bel\label{lem:largesv} The largest singular value of $M$ satisfies $\mu_n \le N^C$ for some $C$ with overwhelming probability. \eel
\bp The largest singular value of $M$ is equal to the largest eigenvalue of the symmetrized sparse matrix $ \bma 0 & M \\
M^\trans &0 \ema$. The lemma then follows from the proof of Lemma 4.3 in \cite{EKYY13}.
\ep

We now obtain short time universality for the new dynamics. Given a matrix $M$, let $\lambda_1(t, M)$ denote the least singular value of $M$ evolved according to the dynamics $(\ref{eqn:OU})$.

\bel\label{lem:shorttime}
Let $(M_N)$ be a sparse matrix ensemble and let  $W$ be a Gaussian ensemble of i.i.d. $\mathcal N(0, N^{-1})$ variables. Given $\eps>0$, there exists $\delta>0$ and a coupling of the processes $\lambda_1(t, M_n)$ and $\lambda_1(t, W)$ such that
\beq \left| \lambda_1(t_a, M_N) -   \lambda_1(t_a, W) \right| \le  N^{-1-\delta}\eeq
for some $t_a\le N^\eps/ N$ with overwhelming probability.
\eel
\bp 
Recall that for any $t$ we have \beq H_t  \overset{d}{=}  H^{(1)}_t +  \sqrt{\frac{r(1-\exp(-t/r))}{N}}G \overset{d}{=} H^{(1)}_t +  \frac{1}{\sqrt N}B_{r(1-\exp(-t/r))}, \eeq
where $G$ is a matrix of i.i.d. standard Gaussians and $B_{r(1-\exp(-t/r))}$ is a matrix Brownian motion considered at the fixed time $r(1-\exp(-t/r))$.

Note that up to a factor of $1+O(t)$, $H_t^{(1)}$ is a sparse matrix and obeys the weak local semicircle law, Lemma \ref{lem:sparseweaklaw}. Then Lemma \ref{lem:sparseweaklaw} holds for $(1+O(t))^{-1}H_t^{(1)}$ on the optimal scale $g=N^{-1+\nu}$ for any $\nu >0$. Further, by invoking Lemma \ref{lem:largesv}, this matrix is $(g,G)$-regular for any such $g$. The additional factor of $1+O(t)$ does not affect the $(g,G)$-regularity of the singular values if $t\le N^{-1-\eps}$ by the argument at the end of the proof of Lemma 6.3 in \cite{HLY15}. In summary, we find that $H_t^{(1)}$ is $(g,G)$-regular for any $g=N^{-1+\nu}$ with overwhelming probability.

By making $\nu$ small enough and choosing constants appropriately in the statement of \mbox{Theorem \ref{thm:homomain}}, we may take $t_a \le N^{-1+\eps}$ in the statement of that theorem and apply it to complete the proof. More precisely, we condition on the entries of $H^{(1)}_{t_a}$ and apply Theorem \ref{thm:homomain} to obtain a conditional coupling of the singular value processes. Then, because the hypotheses of Theorem \ref{thm:homomain} hold for the singular values of $H_{t_a}^{(1)}$ with overwhelming probability, by the weak law Lemma \ref{lem:sparseweaklaw} and Lemma \ref{lem:largesv}, we obtain the desired coupling with overwhelming probability after removing the conditioning. \ep

\subsection{Green function comparison}\label{gfc}

We now control the distribution of the least singular value of a stable matrix ensemble $H_N$ in terms of Green functions. Fix a matrix $H$ from this ensemble. We retain the notation $H_t$ for the dynamics in the previous subsection.

For any $r>0$, define $\chi_r = \one_{(-r,r)}$. For $\eta>0$, we set 
\beq\theta_\eta = \frac{\eta}{\pi(x^2 + \eta^2)} = \frac{1}{\pi} \Im \frac{1}{x - i\eta}.\eeq
In particular, for any $r>0$ we have
\beq \tr \chi_r \star \theta_\eta(H_s) = \frac{N}{\pi} \int_{-r}^r \Im m_s(y+ i\eta)\, dy.\eeq
We fix $\eps,r>0$ and set \beq \eta_1 = N^{-1 - 99\eps}, \quad l = N^{-1 - 3\eps}, \quad l_1= l N^{2\eps}, \quad E = \frac{r}{N}. 
\eeq

We also consider a time parameter $t$. We are interested in the case $0\le t \le N^{\eps_0}/N$ for some small $\eps_0>0$. For $t\le N^{\eps_0}/N$, note that $H_t$ still satisfies the weak local law at the optimal scale. As described in the proof of Theorem 6.3 of \cite{HLY15}, it is a consequence of the weak local law (in particular the fact that $\Im m(z)$ is bounded down to the optimal scale) that there exists $C$ such that, for any interval $I$ with length $| I | \ge N^{-1 + \delta}$, \beq\label{scstates} {|\{ \lambda_i(t) \in I \}| \le C|I|N}\eeq
holds with overwhelming probability.
\bel\label{lem:smear}
Fix $t$ such that $0\le t \le N^{\eps_0}/N$ and $\eps>0$. With overwhelming probability,
\beq |\tr \chi_E(H_t) - \tr \chi_E \star \theta_{\eta_1}(H_t) | \le C\bigg(N^{-2\eps} + n(-E-l, -E+l) + n(E-l, E+l)\bigg).\eeq
\eel
\bp
By the argument in the proof of Lemma 6.1 in \cite{rigidity},
\beq | \chi_E(x) - \chi_E \star \theta_{\eta_1}(x) | \le C \eta_1 \left( \frac{2E}{d_1(x)d_2(x)} + \frac{\chi_E(x)}{d_1(x) + d_2(x)} \right)\eeq
where $d_1=|E-x| + \eta_1$ and $d_2=|-E-x| + \eta_1$, and the right side is bounded by a constant if $\min\{d_i\} \le l$ and is $O(\eta_1/l)$ if $\min\{d_i\}\ge l$. We obtain by (\ref{scstates}), on a set of probability greater than $1 - N^{-D}$,
\begin{align}  | \tr \chi_E(H_t) - \tr \chi_E \star \theta_{\eta_1}(H_t) | &\le C\bigg( \tr f_1(H_t) + \tr f_2(H_t) + \frac{\eta_1}{l} n(-E+l, E-l) \bigg) \\&  + C \bigg(n(-E-l, -E+l) + n(E-l, E+l)\bigg)
 \\ &\le C\bigg( \tr f_1(H_t) + \tr f_2(H_t) + \frac{\eta_1 N^\eps}{l} \bigg) \\ &+ C\bigg(n(-E-l, -E+l) + n(E-l, E+l) \bigg),\end{align}
where \beq f_1(x) = \frac{2\eta_1 E}{d_1(x) d_2(x)}\one(x\le - E  -l ), \quad f_2(x) =  \frac{2\eta_1 E}{d_1(x) d_2(x)}\one(x\ge  E + l ).\eeq
We now describe how to bound $\tr f_2(H_t)$. The term $\tr f_1(H_t)$ is similar. For $x\ge E + l$, we have 
\beq f_2(x) = \frac{2 E \eta_1 }{(|E-x| + \eta_1)(2E + \eta_1 + |E-x| )} \le  \frac{ N^{-1-99\eps} }{|E-x| }. \eeq
Set $\alpha = 3 - E$. We consider the $N$ intervals
\beq \left[E + l, E + \frac{\alpha}{N} \right], \left[ E+ \frac{\alpha}{N} , E+ \frac{2\alpha}{N}\right], \left[ E+  \frac{2\alpha}{N} , E+  \frac{3\alpha}{N}\right], \dots, \left[3 - \frac{\alpha}{N} , 3  \right],\eeq
where the first is of a different size than the rest. By (\ref{scstates}), each interval contains at most $N^\epsilon$ eigenvalues. We also consider the interval $[3,\infty]$. Using this decomposition, we obtain
\beq   \tr f_1(H_t)  \le N^{-1-99\eps} N^\eps \left(\frac{1}{l}+ \frac{N}{\alpha} + \frac{N}{2\alpha} + \frac{N}{3\alpha} + \dots + \frac{N}{N\alpha}  \right) + N N^{-1-99\eps} \le \frac{N^{-98\eps}}{\alpha}(N^{3\eps} + \log(N)).\eeq
This completes the proof.
 \ep

\bel\label{lem:gcount} Fix $t$ such that $0\le t \le N^{\eps_0}/N$ and $\eps>0$. There is a constant $C$ such that, with overwhelming probability,
\beq\tr \chi_{E-l_1}\star \theta_{\eta_1}(H_t)  - CN^{-\eps}\le  \tr \chi_{E}(H_t) \le   \tr \chi_{E+l_1}\star \theta_{\eta_1}(H_t)  + CN^{-\eps}.\eeq\eel
\bp
We see Lemma \ref{lem:smear} holds with $E$ replaced by $y\in [ E-l ,E + l]$. Recall $l_1= l N^{2\eps}$. Hence, with overwhelming probability,
\beq \tr \chi_{E}(H_t) \le \frac{1}{l_1} \int_{E}^{E+l_1} \tr \chi_y(H_t) \, dy \eeq \beq \le  \frac{1}{l_1} \left(\int_{E}^{E+l_1} \tr \chi_y(H_t)\star \theta_{\eta_1}(H_t) \, dy  + CN^{-2\eps} + C n(y-l, y+l) + C n(-y-l, -y+l)\right)\, dy\eeq
\beq \le \tr \chi_{E+l_1}\star \theta_{\eta_1}(H_t) + CN^{-2\eps} + \frac{C l}{l_1}\left( n(E-2l_1, E+2l_1) + n(-E-2l_1, -E+2l_1)\right).\eeq
By (\ref{scstates}), the two counting functions in the above expression are at most $CN^\eps$, so we obtain an error of $N^{-\eps}$ and
\beq \tr \chi_{E}(H_t) \le   \tr \chi_{E+l_1}\star \theta_{\eta_1}(H_t)  + N^{-\eps}.\eeq
A matching lower bound is proved similarly. \ep

As in \cite{PY}, we fix a smooth function $q:\mathbb R \rightarrow \mathbb R_+$ such that $q(x)$ is decreasing for $x\ge 0$,
$q(x)=1$ for $|x|\le 1/9$, and $q(x)=0$ for $|x| \ge 2/9$. 
\bel\label{lem:qbound} Fix $t$ such that $0\le t \le N^{\eps_0}/N$ and $\eps>0$. For any $D>0$, we have
\beq\E q(\tr \chi_{E+l_1}\star \theta_{\eta_1}(H_t))  - N^{-D} \le\P(n(-E,E)=0)\le \E q({\tr \chi_{E-l_1}\star \theta_{\eta_1}(H_t)}) + N^{-D}.\eeq
\eel
\bp When Lemma \ref{lem:gcount} holds, $n(-E,E)=0$ implies ${\tr \chi_{E-l_1}\star \theta_{\eta_1}(H_t)\le 1/9}$ with overwhelming probability. Hence 
\beq\P(n(-E,E)=0) \le \P({\tr \chi_{E-l_1}\star \theta_{\eta_1}(H_t)\le 1/9}) + N^{-D} ,\eeq
and Markov's inequality applied to $q({\tr \chi_{E-l_1}\star \theta_{\eta_1}(H_t)})$ yields
\beq \le \P(q({\tr \chi_{E-l_1}\star \theta_{\eta_1}(H_t)}) \ge 1  ) +N^{-D} \le \E q({\tr \chi_{E-l_1}\star \theta_{\eta_1}(H_t)}) + N^{-D}.\eeq
Also, again using Lemma \ref{lem:gcount},
\beq \E q(\tr \chi_{E+l_1}\star \theta_{\eta_1}(H_t)) \le \P(\tr \chi_{E+l_1}\star \theta_{\eta_1}(H_t) \le 2/9 ) \eeq \beq \le \P( n(-E,E)\le 2/9 + CN^{-\eps} ) +  N^{-D} = \P(n(-E,E)=0) + N^{-D}.\eeq\ep

In the work \cite{QUE}, which analyzed the eigenvector moment flow for generalized Wigner matrices and covariance matrices, the authors developed a purely dynamical approach to Green function comparison. We implement it here in the current context. We require the following modification of Lemma A.1 in \cite{QUE}, which asserts the continuity of the above dynamics. The proof is essentially the same. 

The deformed matrix $\theta^{ab} H_t$ is defined as
$$(\theta^{ab} H_t)_{kl} = f + \theta^{ab}_{kl}( h_{kl}(t)  - f),$$
where $\theta^{ab}_{kl}=1$ if ${k,l} \neq {a,b}$ and some number $0 \le \theta^{ab}_{kl} \le 1$ otherwise, where we impose the symmetry condition $\theta^{ab}_{ab} = \theta^{ab}_{ba}$.
We define the index set $\mathcal I$ to be the entries in the off-diagonal blocks of the $2N\times 2N$ symmetrized matrix:
\beq \mathcal I = \{ (i,j) \colon 1\le i,j,\le 2N, i\le N < j \text{ or } j\le N < i  \}.\eeq

\bel\label{lem:continuity}
Let $H$ be a $2N\times 2N$ symmetric matrix, with entries independent up to the symmetry constraint $h_{ij}=h_{ji}$ and the $N\times N$ blocks on the main diagonal all zero. Suppose the other entries satisfy $\E[h_{ij}] =f$ and $\E[(h_{ij} -f)^2]=s_{ij}$ with $cN^{-1} \le s_{ij} \le C N^{-1}$. Denote $\partial_{ij} = \partial_{h_{ij}}$. Suppose $F$ is a smooth function of the matrix elements $(h_{ij})$ in the upper off-diagonal block of $H$ satisfying
\beq \sup_{0\le s \le t} \sup_{\theta^{ab}} \E [(N^2|h_{ab}(s)-f|^3 +  N |h_{ab} -f| )\partial^{(3)}_{ab} F(\theta^{ab}H_s ) ] \le B \eeq
where the supremum is taken over deformations in the off-diagonal block indices $(i,j)\in\mathcal I$. Then
\beq |\E F(H_t)  - \E F(H_0) | \le C t B.\eeq
\eel
We will use Lemma \ref{lem:continuity} to study the expressions appearing above:
\beq \E q(\tr \chi_{E+l}\star \theta_{\eta_1}(H_s)) =\E q\left( \frac{N}{\pi} \int_{-E-l}^{E+l} \Im m(y+ i \eta_1)\, dy\right). \eeq

\bel\label{lem:gcompare}
There exists $\eps >0$ and $c>0$ such that, for all $t\le N^\eps/N$,
\beq |\E q(\tr \chi_{E+l}\star \theta_{\eta_1}(H_0)) - \E q(\tr \chi_{E+l}\star \theta_{\eta_1}(H_t))| = O(N^{-c}),\eeq
and a similar statement holds with $\chi_{E+l}$ replaced by $\chi_{E-l}$.
\eel
\bp It suffices to control the derivatives of \beq \frac{N}{\pi} \int_{-E-l}^{E+l} \Im m(y+ i \eta_1)\, dy\eeq in order to apply the Lemma \ref{lem:continuity}, since the derivatives of $q$ are bounded independent of $N$. Recall here that $\eta_1= N^{-1-99\eps}$ is below the natural scale. By Lemma \ref{lem:tiny}, which is proved below, we have
\beq \P( |\partial^k_{ab} m_s(y+ i \eta_1)| \le CN^{3(k+1)100\eps} ) \ge 1 - N^{-D}\eeq
for any $D>0$ and $s\in[0,t]$, with a deterministic upper bound of $CN^{3(1+100\eps)}$. Then, moving the derivatives inside the integral, bounding the resulting integrand, and using the fact that the factor of $N$ is canceled by the fact the integral is over an interval of length $O(N^{-1})$, we obtain
\beq B \le C(N N^{-\alpha} N^{\sigma_1}  + N^{-D} N^{\sigma_2}), \eeq
where $\sigma_1$ can be made as small as desired by adjusting $\eps$. Hence $B = O(N^{1-\sigma})$ for some $\sigma>0$. We conclude by choosing $\eps$ small enough in order to make $Bt=o(1)$.\ep

\bel\label{lem:tiny}
With the notation above,
\beq \P( |\partial^k_{ab} m_s(y+ i \eta_1)| \le CN^{3(k+1)100\eps} ) \ge 1 - N^{-D},\eeq
\beq |\partial^k_{ab} m_s(y+ i \eta_1)| \le CN^{3(1+100\eps)}.\eeq
 \eel
 
 \bp
We first suppose that each deformed matrix $\theta^{ab}H_s$ has Green functions elements $G_{ij}(z)$ bounded by a constant $C$ independent of all other parameters for energies $E\in [-1,1]$ and $\eta \ge N^{-1+\eps}$, with overwhelming probability. We fix $s\in [0,t]$ and indices $a,b$. Let $G(z)$ be the Green function for this matrix. Defining 
\beq \Gamma(z) = \max_{i,j} |G_{ij}(z)| \vee 1,\eeq
we have by Lemma 2.1 of \cite{BKY} that
\beq \Gamma(E+ i \eta_1) \le N^{100\eps} \Gamma(E + i \eta),\eeq
where $\eta = N^{-1+\eps}$. We have $\Gamma(E + i \eta)\le C$ for some constant $C$ with overwhelming probability. Given this control over the individual Green functions elements, the argument used in the proof of Lemma 5.2 in \cite{HLY15} proves the first claim. The second claim follows from the trivial deterministic bound
\beq G(E + i \eta) \le \frac{1}{\eta}.\eeq

We now show that we have a uniform constant bound on the deformed $G_{ij}$. This follows from the proof of Theorem 6.3 in \cite{HLY15}. Note that in that reference the bound is obtained with overwhelming probability for fixed $s$, and implicitly a standard stochastic continuity argument gives the bound uniformly in $s$.
 \ep
 
  \section{Random matrices with correlated entries} \label{a:correlated}
  
 In this section, we prove the universality of the least singular value for a class of non-symmetric square random matrices with correlated entries. Our main result is Theorem \ref{correlatedthm}. We first sketch the proof of a local law for the smallest singular values.  Our method extends the one in \cite{Ch16}, and here we consider the more general case of exponentially decaying correlations. We then show how to adapt the method in the previous section on the removal of the time evolution. As a consequence, universality holds for the smallest singular value. The same result can be proved for non-Hermitian random matrices with complex entries in the same way.

 \subsection{Model}

Consider a family of centered real random variables $(x_{ij})_{1\leq i,j\leq N}$ that satisfy
\[
	\E[x_{ij}x_{kl}] =\frac{1}{N} \xi_{ijkl},\quad 1\leq i,j,k,l\leq N.
\]
We introduce a sparsity parameter $q=N^\tau$ for some $\tau\in(0,1]$.  Assume that for any $p\geq 2$, there is a constant $\mu_p<+\infty$ such that
\[
	\sup_{i,j} \E[\abs{x_{ij}}^p]\leq \frac{\mu_p^p}{Nq^{p/2-1}},\quad \forall N\in\N.
\]
Furthermore, we assume that $\xi$ has a profile, in the sense that there is a function $\phi:[0,1]^2\times \Z^2\to \R$ such that
\[
	\xi_{ijkl}=\phi (i/N,j/N, k-i,l-j),
\]
and that $\phi$ is piecewise H{\"o}lder-continuous with respect to the first two variables. We impose exponential decay on the correlation. For any index set $\mathcal{A}\in \{(i,j):1\leq i,j\leq N\}$ we define $\mathcal{F}_\mathcal{A}$ to be the $\sigma$-algebra generated by $(X_{ij})_{(i,j)\in\mathcal{A}}$. For any two index sets $\mathcal{A},\mathcal{B}$ we define their distance by
\[
	\dd(\mathcal{A},\mathcal{B})= \max_{(i,j)\in\mathcal{A},(i',j')\in\mathcal{B}} \abs{i-i'}\vee\abs{j-j'}.
\]
We assume that there are universal constants $c_1,c_2>0$ such that for any random variables $Z_1\in\mathcal{F}_\mathcal{A},Z_2\in\mathcal{F}_\mathcal{B}$ with $\mathrm{Var}[Z_1]=\mathrm{Var}[Z_2]=1$, the following inequality holds,
\be\label{decay}
	\mathrm{Cov}[Z_1,Z_2] \leq c_1\exp(-c_2\dd(\mathcal{A},\mathcal{B})).
\ee
As usual we write the symmetrized version of $X$, namely a $2N$ by $2N$ matrix $H$ defined by
\[
	H = \bma 0 & X^* \\X & 0\ema.
\]
It is easy to see that $H$ is a special case of the model in \cite{Ch16} without the positive definite condition (see Definition 2.2 in \cite{Ch16}), since $H$ has many $0$ entries.  However,  we can still consider an alternative positive definiteness condition in the current case.
\bed\label{pd}
	Let $\Sigma^{(N)}\in \R^{N^2\times N^2}$ be the covariance matrix of the family of real random variables $(x_{ij})_{1\leq i,j\leq N}$ with $\E[X_{ij}X_{kl}] = \xi^{(N)}_{ijkl}$.  We say that $\xi$ is positive definite with lower bound $c_0>0$ if $\Sigma^{(N)}\geq c_0$ for all $N$.
\eed

We now state the main result of this section.

\bet 
\label{correlatedthm}
For the class of correlated sparse matrices whose correlation comes from a positive definite profile function, as defined above, the conclusion of Theorem \ref{thm:mainresult} holds.
\eet 

\subsection{Local law}
\subsubsection{Concentration}

Condition \ref{decay} is weaker than the finite-ranged correlation enforced in \cite{Ch16}, because every pair of entries could be correlated, although exponentially weakly.  Nevertheless the same concentration estimates hold for linear combinations and quadratic forms (see Lemma 3.6 in \cite{Ch16}). 

\bel\label{ldp}
Let $k\in\N$, $\mathcal{A}\subset\N$ be such that $\dd(\{k\},\mathcal{A})\geq \log^2N$. Let $(A_i),(B_{ij})$ be families of random variables that are $\mathcal{F}_\mathcal{A}$-measurable with upper bounds $\max_i\abs{A_i} \vee\max_{i\neq j} \abs{B_{ij}}\leq N^{100}$ . Then, 
\[
	\absa{\sum_i A_ix_{ki} - \E\left[\sum_i A_ix_{ki} \right]} \prec \frac{\max_i\abs{A_i}}{\sqrt{q}} + \sqrt{\frac{1}{N} \sum_i \abs{A_i}^2}+\exp(-c\log^2N).
\]
\[
	\absa{\sum_{i,j} B_{ij}x_{ki}x_{kj} -\E\left[ \sum_{i,j} B_{ij}x_{ki}x_{kj}\right]} \prec \frac{\max_{i\neq j} \abs{B_{ij}}}{\sqrt{q}}+\sqrt{\frac{1}{N} \sum_{i,j} \abs{B_{ij}}^2}+\exp(-c\log^2N).
\]
\eel
\begin{proof}
 Note that $x_{ki}$ is centered, therefore
\[
	\sum_i\E[A_ix_{ki}] = \sum_i \mathrm{Cov}[A_i,x_{ki}] \leq c_1\exp(-c_2 \log^2N)\sum_i\left( \mathrm{Var}[A_i] \mathrm{Var}[x_{ki}]\right),
\]
which is bounded by $\exp(-c\log^2N)$ for some constant $c>0$. Therefore the first inequality is reduced to proving
\[
	\absa{\sum_i A_ix_{ki} }\prec \frac{\max_i\abs{A_i}}{\sqrt{q}} + \sqrt{\frac{1}{N} \sum_i \abs{A_i}^2}.
\]
We split the sum into $\lfloor\log^2N\rfloor$ parts, each part being the sum of weakly correlated random variables.  Specifically, write $S_l= \sum_i A_{i\lfloor\log^2N\rfloor+l}$ so that
\[
	\sum_i A_ix_{ki} = \sum_{ 1\leq l\leq \lfloor\log^2N\rfloor} S_l.
\]
Heuristically, each $S_l$ can be viewed as the sum of independent random variables, because the summands have very weak correlation with each other.  Let  $\tilde A_i$, $\tilde x_{ki}$ be independent copies of $A_i$ and $x_{ki}$, and let $\tilde S_l$ be defined likewise by replacing $A_i,x_{ki}$ with their copies.  It is easy to see that for any $p\geq 2$,
\[
	\E\abs{S_l}^p =\E\abs{\tilde S_l}^p + O(c_p\exp(-c\log^2N))
\]
by expanding out $\abs{S_l}^p$ and collecting all the cross terms, which are exponentially small.  This implies that
\[
	\abs{S_l} \prec\abs{ \tilde{S}_l }+ O(\exp(-c\log^2N)).
\]
For $\tilde{S}_l$ one can apply Lemma A.1 in \cite{EKYY13} to get 
\[
	\abs{ \tilde{S}_l }\prec \frac{\max_i\abs{A_i}}{\sqrt{q}} + \sqrt{\frac{1}{N} \sum_i \abs{A_i}^2}.
\]
Therefore we summarize the estimates above and see
\[
	\sum_i A_ix_{ki} \prec \log^2N  \left(\frac{\max_i\abs{A_i}}{\sqrt{q}} + \sqrt{\frac{1}{N} \sum_i \abs{A_i}^2}+O(\exp(-c\log^2N))\right).
\]
The factor $\log^2N$ can be absorbed into $\prec$ by definition.  This proves the first inequality.

The second inequality follows from a very similar argument.  One only needs to split the sum into $O(\log^2N)$ parts, each of which is the sum of weakly correlated random variables.  The weak correlation will not worsen the estimate, as we saw in the proof of the first inequality.\qed
\end{proof}

\subsubsection{Self-consistent equation}
As usual we define the Green function $G(z)$ by
\[
	 G(z)= (H-z)^{-1}.
\]
Define a map $\Xi: \C^{2N\times 2N} \to  \C^{2N\times 2N}$ through
\[
	\Xi(M)= \E[HMH],
\]
which can be explicitly defined entrywise:
\[
	\Xi(M)_{ik}= \sum_{1\leq j,\, l\leq 2N} \E[h_{ij}h_{kl} M_{jl}].
\]
We introduce three control parameters:
\beq
	\Gamma= 1\vee \max_{i,j}\abs{G_{ij}},\quad \gamma= 1\vee \max_i \max_{\mathbb{I},\mathbb{J}}\norma{\left(G_{\mathbb{I},\mathbb{I}}^{(\mathbb{J})}\right)^{-1}},\quad \Phi= \frac{1}{\sqrt{N\eta}}+\frac{1}{\sqrt{q}}.
\eeq
With the help of Lemma \ref{ldp} and repeating the argument of Lemma 3.9 in \cite{Ch16}, one can prove an estimate that is the same as Lemma 3.9 in \cite{Ch16}:
\be\label{schwinger1}
	G(-z-\Xi(G)) = I + O_\prec(\Gamma^5\gamma^3 \Phi).
\ee
Here the $O_\prec$ notation is in the entrywise sense. We omit the details here, but point out that the only difference in the proof is that the error term here is bigger by a factor of $O(\log^2N)$, which is negligible in the context of stochastic domination.  It remains to show that $G$ is close to the solution $M$ to the following equation.
\be\label{schwinger2}
	M(-z-\Xi(M)) = I
\ee

The whole argument in \cite{Ch16} goes through except for the part where the positive-definiteness condition on the tensor $\xi$ is used to prove the  stability of \eqref{schwinger2} (that is, the solution is stable under small perturbation of the equation). We now discuss the necessary changes. 

In \cite{Ch16} the equation was transformed into a continuous version.  It was shown by a discretization argument that \eqref{schwinger2} is stable in the bulk if and only if the following continuous equation is stable in the bulk. 
\be\label{qve}
	u(\theta,s) = \frac{1}{-z-S u(\theta,s)},\quad  \theta,s\in[0,1]
\ee
Here the operator $S$ is given by 
\[
	Su(\theta,s) = \iint \hat\phi(\theta,\vartheta,s,t)\, d \vartheta\, d t,
\]
where $\hat\phi$ is the Fourier transform of the symmetrized version of $\phi$ in the latter two variables (since $H$ is the symmetrized version of $X$). In particular,
\[
	\hat\phi(\theta,\vartheta,s,t)= \sum_{1\leq k,l\leq N}\phi(\theta,\vartheta,k,l)\exp(\mathrm{i}2\pi(sk-tl)),\quad \text{ for } (\theta,\vartheta)\in[0,1/2]\times [1/2,1],
\]
and $\hat\phi(\theta,\vartheta,\cdot,\cdot) = \hat\phi(\vartheta,\theta,\cdot,\cdot)$ for $ (\theta,\vartheta)\in[1/2,1]\times [0,1/2]$ and $\hat\phi=0$ for other $(\theta,\vartheta)$.  
In \cite{Ch16} it was shown that $0<c_0\leq \hat\phi(\theta,\vartheta,s,t) \leq C_0$ for some universal constants  $c_0$ and $C_0$ (Lemma 4.15) under the positive definite condition in that paper.  Here we use the new positive definite condition (Definition \ref{pd}) to prove upper and lower bounds on $\hat\phi$. 
\bel
Suppose that $\xi$ is positive definite  with lower bound $c_0>0$ in the sense of Definition \ref{pd}.  For $\theta,\vartheta\in[0,1/2]\times [1/2,1]$ or $\theta,\vartheta \in [1/2,1]\times[0,1/2]$, we have
\[
	\hat\phi(\theta,\vartheta,s,t)\in[c_0,C_0],\quad \forall s,t\in[0,1].
\]
Here $C_0$ is an universal constant depending on $c_1,c_2>0$ in \eqref{decay}.
\eel
\begin{proof}
	For any $(s,t)\in[0,1]^2$, take an arbitrary real continuous function $g\in C([0,1]^2)$ with $\iint \absa{g}^2 =1$. For each $N\in\N$ define a random variable
	\[
		Y_N= \frac{1}{N}\sum_{1\leq i,j\leq N}h_{ij} g\left(\frac{i}{N},\frac{j}{N}\right) \exp(\mathrm{i}2\pi(si-tj)).
	\]
	By definition of positive definiteness and the decay of correlation, we have $c_0\leq \mathrm{Var}Y_N\leq C_0$.  One can explicitly compute the variance of $Y_N$:
	\[
		\mathrm{Var}Y_N = \frac{1}{N^2} \sum_{i,j,k,l} \phi\left( \frac{i}{N},\frac{j}{N},k,l\right) g\left(\frac{i}{N},\frac{j}{N}\right)g\left(\frac{i+k}{N},\frac{j+l}{N}\right) \exp(\mathrm{i}2\pi(sk-tl)).
	\]
	Let $N\to\infty$ and use the fact that $c_0\leq \mathrm{Var}Y_N\leq C_0$. We have
	\[
		 \iint g(\theta,\vartheta)\hat\phi(\theta,\vartheta,s,t) \, d \theta\, d \vartheta \in[c_0,C_0].
	\]
	Since $g$ was arbitrary, we conclude that $\hat\phi(\theta,\vartheta,s,t)\in[c_0,C_0]$. \qed
\end{proof}

The stability of equation \eqref{qve} is very similar to the case in \cite{Ch16} and was analyzed in \cite{gram}. It follows from Proposition 3.10 (ii) in \cite{gram} that the following estimate holds.
\bep
Let $u$ solve \eqref{qve} and $u'$ solve a perturbed version of \eqref{qve}, namely
\[
	u' = \frac{1}{-z-Su'} + r.
\]
There exist universal constants $\eps>0$, and $C>0$ such that if $\abs{\Re z}\leq\eps, \Im z\in(0,10], {\norma{u-u'}_\infty\leq \eps}$, then
\[
	\norma{u'-u}_\infty \leq C \norma{r}_\infty.
\]	
\eep\label{stability}
Via the discretization method in \cite{Ch16}, we can prove the stability for equation \eqref{schwinger2}.  Below, $\norma{A}_\infty$ means $\max_{i,j} \abs{A_{ij}}$.

\bep
Let $M$ be the solution to \eqref{schwinger2} and let $M'$ solve a perturbed version of \eqref{schwinger2}, namely
\[
	M'(-z-\Xi(M'))=I +R.
\]
There exist universal constants $\eps>0$, and $C>0$ such that if $\abs{\Re z}\leq\eps, \Im z\in(0,10], {\norma{M-M'}_\infty<\eps}$, then,
\[
	\norma{M-M'}_\infty<C\norma{R}_\infty.
\]
\eep
Proposition \ref{stability} and the estimate (\ref{schwinger1}) allow us to prove the following theorem.
\bet
Let $M$ the the solution of \eqref{schwinger2}.  There exists a universal constant $\eps>0$ such that for any $\kappa>0$,
\[
	\max_{ij} \abs{G_{ij}-M_{ij}}\prec\Phi,
\]	
uniformly for all $z\in\{E+\mathrm{i}\eta: \eta\in[N^{-1+\kappa}, 10), E\in[-\eps,\eps]\}$.
\eet
\bec\label{cor:corweaklaw}
There exist universal constants $\eps>0$ and $c>0$ such that with overwhelming probability,
\beq
	c\leq \Im\left( \frac{1}{2N}\tr G\right) \leq c^{-1}
\eeq
uniformly for all $z\in\{E+\mathrm{i}\eta: \eta\in[N^{-1+\kappa}, 10), E\in[-\eps,\eps]\}
$.
\eec

\subsection{Universality}
Let $(B_{ij}(t))_{1\leq i,j\leq N}$ be a family of Brownian motions that has the same correlation structure as $(x_{ij})$:
\beq\label{eqn:corOU}
	\E[B_{ij}(t)B_{kl}(t)] =t \E[x_{ij}x_{kl}] = t\xi_{ijkl}/N.
\eeq
Define $x_{ij}(t)$ by the SDE
\[
	\dd x_{ij} =  d B_{ij} - \frac{x_{ij}}{2}\, d t.
\]
We show that when $t\ll N^{-1}\sqrt{q}$, the evolution does not affect the local statistics of the smallest singular values.  Following the argument in \cite{Ch16}, we prove the same result as in Lemma 6.1 in \cite{Ch16}.
\bel
	Let $x=(x_k)_{1\leq k\leq m}$ be an array of real centered random variables such that $\sup_k\E[\abs{x_k}^3]\leq \kappa_3^3$ and $\mathrm{Corr}[Z_1,Z_2]\leq c_1\exp(-c_2d)$ for all nontrivial random variables ${Z_1\in\sigma(x_1,\cdots,x_k)}$, ${Z_2\in\sigma(x_{k+d},\cdots,x_m)}$ and any $1\leq k\leq k+d\leq m$. Let $f$ be a $C^2$ function on $\R^m$ with $\norma{D^2f}\vee \kappa_3\vee m\leq N^{100}$. Then, 
\[
	\E[f(x)x_i] = \sum_k \E[\partial_kf(x)]\E[x_ix_k] + O\left(\log^2 N\norm{D^2f}_\infty \kappa_3^3+ \exp(-c\log^2N)\right).
\]
\eel
\bp
If $f$ is a linear function in $x$, then the equality is exact without error terms.  In general, define
\[
	\mathbb{T} = \{j: \abs{i-j}\leq \log^2N\},\quad \mathbb{U}= \{ j: \abs{i-j} \leq 2\log^2N\}.
\]
Denote $x^{(\mathbb{T})}= (x_k\one_{k\notin \mathbb{T}})$, $x^{(\mathbb{U})}= (x_k\one_{k\notin\mathbb{U}})$. By Taylor's expansion,
\[
	f(x) = f(x^{(\mathbb{T})}) + \sum_{k\in\mathbb{T}} \partial_k f(x^{(\mathbb{T})})x_k + \frac{1}{2} \sum_{k,l\in\mathbb{T}} \int_0^1 (1-t)\partial_{kl}f(x^{(\mathbb{T})}+t(x-x^{(\mathbb{T})}) ) x_kx_l\, d t.
\]
We expand the second term further,
\[
	\sum_{k\in\mathbb{T}} \partial_k f(x^{(\mathbb{T})})x_k =\sum_{k\in\mathbb{U}} \partial_k f(x^{(\mathbb{T})})x_k +\sum_{k\in\mathbb{T},l\in\mathbb{U}} \int_0^1 (1-t) \partial_{kl} f(x^{(\mathbb{U})}+t(x^{(\mathbb{T})}-x^{(\mathbb{U})}))x_kx_l\, d t.
\]
Therefore, $f(x)$ can be written as
\beq\label{eqn:defexpand}
	f(x) = f(x^{(\mathbb{T})})+ \sum_{k\in\mathbb{U}}\partial_k f(x^{(\mathbb{U})})x_k + O\left( \sum_{k,l\in\mathbb{U}} \sup_{\theta\in[0,1]^\mathbb{U}} \abs{f(\theta x)}\abs{x_kx_l}\right).
\eeq
Here $\theta x= (x_k\one_{k\in\mathbb{U}}\theta_k + x_k\one_{k\notin\mathbb{U}})$.
Therefore, we can compute
\[
	\E[f(x)x_i] = \E[f(x^{(\mathbb{T})})x_i] + \sum_{k\in\mathbb{T}}\E[\partial_k f(x^{(\mathbb{U})})x_kx_i] + O\left( \norma{D^2 f}_\infty \kappa_3^3\right).
\]
The first and second term are expectations of products of weakly correlated random variables. So,
\[
	\E[f(x)x_i] =\sum_{k,l} \E[\partial_{kl}f(x^{(\mathbb{U})})] \E[x_kx_l]+O\left( \norma{D^2 f}_\infty \kappa_3^3+\exp(-c\log^2N)\right).
\]
Using Taylor expansion again, we can replace $\partial_{kl}f(x^{(\mathbb{U})})$ by $\partial_{kl}f(x)$ with the cost of a small error term. Therefore,
\[
	\E[f(x)x_i] =\sum_{k,l} \E[\partial_{kl}f(x)] \E[x_kx_l]+O\left( \norma{D^2 f}_\infty \kappa_3^3+\exp(-c\log^2N)\right).\]\ep

We can similarly imitate the proof of Lemma 6.2 of \cite{Ch16} to prove the following.
\bel
Suppose $f$ is a $C^{3}$ function on $\mathbb R^{N\times N}$. Then
\beq \E[f(H_t) - f(H_0) ]  =O\left(\exp(-c\log^2N) + tNq^{-1/2}  \E  \sup_{\theta} \partial^{(k)} f(\theta H) \right) .\eeq
\eel

The rest of the argument is essentially the same as the one in Section \ref{sec:removal}. Similarly to Section \ref{sec:newdynamics}, we may decompose the correlated OU dynamics (\ref{eqn:corOU}) as
\beq x_{ij}(t) = \tilde x_{ij}(t) + w_{ij},\eeq
where the $(w_{ij})$ are i.i.d. Gaussian and the $\tilde x_{ij}(t)$ have a positive definite correlation structure. Then we see the local law holds for $\tilde X(t) = (\tilde x_{ij}(t))$, and the rest of the arguments in Section \ref{sec:removal} go through, since they do not rely on the structure of the matrix. The only necessary change is that in the proof of Lemma \ref{lem:tiny} we need to use a different method to show the (non-deformed) Green function entries are uniformly bounded by a constant down to the scale $N^{-1+\delta}$. Using the local law, this reduces to showing the entries of $M$ are bounded, and this is a consequence of the definition of $M$ and the analogue of Lemma 4.20 in \cite{Ch16}. The regularity necessary for Theorem \ref{thm:homomain} is provided by Corollary \ref{cor:corweaklaw} and the analogue of Lemma \ref{lem:largesv} for the correlated model in this Section. The analogue is proved by splitting into $\log^4(N)$ weakly correlated matrices and applying the argument in the proof of Lemma \ref{lem:largesv}. This completes the proof of Theorem \ref{correlatedthm}.

\appendix

\section{Singular value dynamics}  \label{a:sde}

This appendix collects information on the SDE

\beq\label{eqn:appendixSDE} d\lambda_k = \frac{1}{\sqrt{N}} dB_k + \frac{1}{2N}\sum_{ j\neq k}\left( \frac{1}{\lambda_k - \lambda_j } + \frac{1}{\lambda_k + \lambda_j}  \right) dt. \eeq
\subsection{Existence and uniqueness of solutions}
Let $\Delta$  be the region where $\lambda_1 < \lambda_2 <\dots< \lambda_N$ and $|\lambda_1| < |\lambda_2| <\dots< |\lambda_N|$. We show that given initial data in $\Delta$, there is a unique strong (continuous) solution that stays in $\Delta$ for all time. We follow the arguments in Section 4.3 of \cite{AGZ}, explaining the necessary changes. We also show the solutions of this equation are the singular values of a matrix Brownian motion process. 

Throughout, $N$ will be fixed, and $\lambda(t) = (\lambda_1(t),\dots,\lambda_N(t) )$.

\bel Fix an initial condition $\lambda(0)\in \Delta$. There exists a unique strong solution ${(\lambda(t)_{t\ge 0}) \in C(\R^+, \Delta)}$.

\eel

\proof In the proof of Lemma 4.3.3 in \cite{AGZ}, replace the given definition of $f$ with 

\beq f(x)= \frac{1}{N}\sum_i x_i^2  - \frac{1}{2 N^2} \sum_{i\neq j}  \log|x_i - x_j| + \log | x _i + x_j|.\eeq
Then the estimates (4.3.6) given still hold, and
\beq df(\lambda^R(t)) = \sum_{i=1}^N  \partial_i f(\lambda^R(t))\, d\lambda^R_i + \frac{1}{2} \sum_{i,j} \partial_i\partial_j f(\lambda^R(t))\, d\langle \lambda_i^R , \lambda_j^R \rangle.\eeq

As in \cite{AGZ}, we define 
\beq u_{i,1} = \sum_{k\neq i} \frac{1}{x_i - x_k}, \quad u_{i,2} = \sum_{k\neq i} \frac{1}{(x_i - x_k)^2},\eeq
\beq \overline u_{i,1} = \sum_{k\neq i} \frac{1}{x_i + x_k}, \quad \overline u_{i,2} = \sum_{k\neq i} \frac{1}{(x_i + x_k)^2}.\eeq
We have the identities 
\beq \sum x_i(u_{i,1}(x) + \overline  u_{i,1}(x)) = N(N-1)\eeq
\beq\sum (u_{i,1} + \overline u_{i,1})^2 - u_{i,2} - \overline u_{i,2} = 0.\eeq

The equation becomes (suppressing the $\lambda^R_i$ in the arguments of the $u$ functions):

\beq df(\lambda^R(t)) = 1 \, dt + \frac{1}{N^2} \sum_{i}\left( \lambda^R_i  -\frac{1}{2N}(u_{i,1} +\overline u_{i,1} ) \right) u_{i,1} dt + \frac{1}{N^3} \sum u_{i,2} \, dt + dM(t)\eeq
\beq=   1 +  \frac{1}{2} - \frac{1}{N} + dM(t).\eeq
Since $M$ is a martingale with zero expectation,
\beq E[f(\lambda^R(t\wedge T_M) ] \le 1.5 E[t\wedge T_M]  + f(\lambda^R(0)).\eeq
Now we are finished, as in \cite{AGZ}, by a Borel-Cantelli argument. 
\qed \\

 We now explain why solutions to \eqref{eqn:appendixSDE} have the same distribution as the singular values of $M_t$ from Section \ref{sec:deterministic}. We recall the SDE for the eigenvalues of ${M_t^\trans M_t}$ given in Appendix C of \cite{QUE}:

\beq\label{eqn:covarianceSDE} d\lambda_k = 2\sqrt{\lambda_k} \frac{dB_k}{\sqrt{N}} + \left(1 + \sum_{j\neq k} \frac{\lambda_k + \lambda_l}{\lambda_k - \lambda_l} \right)\, dt.\eeq
Existence and uniqueness of solutions to \eqref{eqn:covarianceSDE} was shown in \cite{Bru}.

An application of It\^{o}'s lemma to \eqref{eqn:appendixSDE} yields

\beq d(\lambda_k^2) = \frac{2\lambda_k}{\sqrt N} \, dB_k + 1\, dt + \sum_{j\neq k} \frac{\lambda^2_k + \lambda^2_j}{\lambda^2_k - \lambda^2_j} \,dt.\eeq
Then $\lambda_k^2$ almost solves the SDE given in \cite{QUE}, except here we are not always choosing the positive square root of $\lambda_k^2$. However, we obtain a weak solution by noting that $\lambda^2$ solves \eqref{eqn:covarianceSDE} with the Brownian motions chosen as $\overline B_k = \operatorname{sgn}(\lambda_k) B_k$. (Note that by the L\'{e}vy criterion, the $\overline B_k$'s are indeed independent Brownian motions.) Hence the solutions of \eqref{eqn:appendixSDE} have the desired distribution.

\begin{comment}
It remains show that SDE \eqref{eqn:appendixSDE} governs the evolution of the singular values of the matrix process $M_t = V + B_t$, and that $dB_i$ relates to $B_t$ as claimed in Section \ref{sec:deterministic}.  This can be done by imitating the proof of Lemma 4.3.4 in \cite{AGZ}.
\end{comment}

\subsection{Interpolation}\label{sec:interpolationdetails}
Here we provide the details of the construction of the interpolated solutions (\ref{eqn:interp}) and show they are differentiable with respect to $\alpha$. 

We first construct solutions $z_i(t,\alpha)$ for $\alpha\in \mathbb Q \cap [0,1]$ using the argument in the previous subsection. Because there are a countable number of solutions, each of which exists individually except possibly on some set of measure zero in the probability space $\Omega$, they all exist and satisfy the SDE simultaneously on a set of full measure $E_1\subset \Omega$. For $\alpha_1,\alpha_2 \in \mathbb Q \cap [0,1]$, define
\beq \tilde u_i(t,\alpha_1,\alpha_2)=z_i(t,\alpha_1)-z_i(t,\alpha_2).\eeq
Then 
\beq \partial_t \tilde u_i(t)=\sum_j B_{ij}(\tilde u_j - \tilde u_i), \quad  B_{ij} = \frac{\one_{i\neq \pm j}}{2N(z_i(\alpha_1)-z_j(\alpha_1))(z_i(\alpha_2)-z_j(\alpha_2))},\eeq
\beq \tilde u_i(0) = (\alpha_1 - \alpha_2)(z_i(0,1) - z_i(0,0)).\eeq
Suppose $\tilde u_k(t)$ is a particle where $\tilde u_k(t)=\max_{i} \tilde u_i(t)$. Because the particles are ordered, the coefficients $B_{ij}$ are positive, so $\partial_t \tilde u_k (t) \le 0$. We conclude that $\| \tilde u(t) \|_\infty$ is non-increasing, and \beq\label{eqn:lip} \| \tilde u(t) \|_\infty \le (\alpha_1 - \alpha_2) (\|z(0,1)\|_\infty + \|z(0,0)\|_\infty)\eeq
holds uniformly for all $t$ on $E_1$.

Since $z_i(t, \alpha)$ is Lipschitz in $\alpha$ (with Lipschitz constant depending on the random initial data), it extends uniquely to a (random) function $z(t,\alpha)$ continuous in $\alpha\in[0,1]$. Further, since the uniform limit of continuous functions is continuous, the Lipschitz estimate shows the paths in the variable $t$ are continuous for all $\alpha$. 

Fix $\alpha_0\in [0,1]$. If $\tilde z_i(t,\alpha_0)$ is a solution a.s., then the same reasoning that led to (\ref{eqn:lip}) shows that $\tilde z_i(t,\alpha_0,\omega)=z_i(t,\alpha_0,\omega)$ for a set of full measure in $\Omega$. (Note, however, that this set of full measure may vary with the choice of $\alpha_0$). By Fubini's theorem, $z_i(t,\alpha,\omega)$ is solution for a set of $(\omega,\alpha)$ of full measure in the product space $\Omega \times [0,1]$. This completes the construction of the interpolated solutions.

Fix $\omega \in E_1$. Since Lipschitz functions are differentiable almost everywhere and their derivatives satisfy the fundamental theorem of calculus, we see that $\partial_\alpha z_i(t,\alpha,\omega)$ exists for almost every $\alpha \in [0,1]$ (with the exceptional set depending on $\omega$) and 
$$z_i(t,1,\omega) - z_i(t,0,\omega) = \int_0^1 \partial_\alpha z(t,\alpha,\omega)\, d\alpha.$$
Hence this relation holds for every $\omega \in E_1$ and therefore almost surely.

%%%%%%%%%%%

\bibliography{mybib2}{}

\begin{thebibliography}{10}

\bibitem{sparsecovar}
B.~Adlam.
\newblock The local {M}archenko-{P}astur law for sparse covariance matrices.
\newblock {\em \url{http://www.people.fas.harvard.edu/~adlam/thesis.pdf}},
  2013.

\bibitem{gram}
J.~Alt, L.~Erd\H{o}s, and T.~Kr{\"u}ger.
\newblock Local law for random {G}ram matrices.
\newblock {\em Electronic Journal of Probability}, 22(25), 2017.

\bibitem{AGZ}
G.~W. Anderson, A.~Guionnet, and O.~Zeitouni.
\newblock {\em An Introduction to Random Matrices}.
\newblock Cambridge University Press, 2010.

\bibitem{BSbook}
Z.~Bai and J.~Silverstein.
\newblock {\em Spectral Analysis of Large Dimensional Random Matrices}.
\newblock Science Press, 2006.

\bibitem{BR15}
A.~Basak and M.~Rudelson.
\newblock Invertibility of sparse non-hermitian matrices.
\newblock {\em Advances in Mathematics}, 310:426--483, 2017.

\bibitem{BKY}
R.~Bauerschmidt, A.~Knowles, and H.-T. Yau.
\newblock Local semicircle law for random regular graphs.
\newblock {\em Comm. Pure Appl. Math.}, 70:1898--1960, Oct. 2017.

\bibitem{lectures}
F.~Benaych-Georges and A.~Knowles.
\newblock Lectures on the local semicircle law for {W}igner matrices.
\newblock {\em Preprint arXiv:1601.04055}, 2016.

\bibitem{Bi97}
P.~Biane.
\newblock On the free convolution with a semi-circular distribution.
\newblock {\em Indiana Univ. Math. J.}, 46(3):705--718, 1997.

\bibitem{wigfixed}
P.~Bourgade, L.~Erd{\H{o}}s, H.-T. Yau, and J.~Yin.
\newblock Fixed energy universality for generalized {W}igner matrices.
\newblock {\em Comm. Pure Appl. Math.}, Dec. 2015.

\bibitem{QUE}
P.~Bourgade and H.-T. Yau.
\newblock The eigenvector moment flow and local quantum unique ergodicity.
\newblock {\em Comm. Math. Phys.}, 2013.

\bibitem{Bru}
M.-F. Bru.
\newblock Diffusions of perturbed principal component analysis.
\newblock {\em Journal of Multivariate Analysis}, 29(1):127--136, 1989.

\bibitem{Ch16}
Z.~Che.
\newblock Universality of random matrices with correlated entries.
\newblock {\em Electronic Journal of Probability}, 22(30):1--38, 2017.

\bibitem{claeys2017boundaries}
T.~Claeys, T.~Neuschel, and M.~Venker.
\newblock Boundaries of sine kernel universality for {G}aussian perturbations
  of {H}ermitian matrices.
\newblock {\em arXiv preprint arXiv:1712.08432}, 2017.

\bibitem{C16}
N.~Cook.
\newblock Lower bounds for the smallest singular value of structured random
  matrices.
\newblock {\em The Annals of Probability}, 46(6):3442--3500, 2018.

\bibitem{dumitriu2018sparse}
I.~Dumitriu and Y.~Zhu.
\newblock Sparse general wigner-type matrices: Local law and eigenvector
  delocalization.
\newblock {\em arXiv preprint arXiv:1808.07611}, 2018.

\bibitem{EPR10}
Erd\H{o}s, S.~P{\'e}ch{\'e}, J.~A. Ramirez, and B.~Schlein.
\newblock Bulk universality for {W}igner matrices.
\newblock {\em Comm. Pure Appl. Math.}, 63(7):895--925, 2010.

\bibitem{EKY13}
L.~Erd\H{o}s, A.~Knowles, and H.-T. Yau.
\newblock Averaging fluctuations in resolvents of random band matrices.
\newblock {\em Ann. Inst. Henri Poincar{\'e}}, 14(8):1837--1926, 2013.

\bibitem{scgeneral}
L.~Erd\H{o}s, A.~Knowles, H.-T. Yau, and J.~Yin.
\newblock The local semicircle law for a general class of random matrices.
\newblock {\em Electronic Journal of Probability}, 18(59), 2013.

\bibitem{EKYY13}
L.~Erd\H{o}s, A.~Knowles, H.-T. Yau, and J.~Yin.
\newblock Spectral statistics of {E}rd{\H{o}}s-{R}{\'e}nyi graphs {I}: local
  semicircle law.
\newblock {\em Ann. Probab.}, 41(3B):2279--2375, 2013.

\bibitem{EYbook}
L.~Erd\H{o}s and H.-T. Yau.
\newblock Dynamical approach to random matrix theory.
\newblock {\em Courant Lecture Notes in Mathematics}, 28, 2017.

\bibitem{localrelax2}
L.~Erd{\H{o}}s, B.~Schlein, and H.-T. Yau.
\newblock Universality of random matrices and local relaxation flow.
\newblock {\em Invent. Math.}, 185(1):75--119, 2011.

\bibitem{localrelax1}
L.~Erd{\H{o}}s, B.~Schlein, H.-T. Yau, and J.~Yin.
\newblock The local relaxation flow approach to universality of the local
  statistics for random matrices.
\newblock {\em Annales de l'I.H.P. Probabilit{\'e}s et statistiques},
  48(1):1--46, 2012.

\bibitem{survey}
L.~Erd{\H{o}}s and H.-T. Yau.
\newblock Universality of local spectral statistics of random matrices.
\newblock {\em Bull. Amer. Math. Soc.}, 49(3):377--414, 2012.

\bibitem{G99}
D.~Grabiner.
\newblock Brownian motion in a {W}eyl chamber, non-colliding particles, and
  random matrices.
\newblock {\em Ann. Inst. Henri Pointcar{\`e} Probab. Stat.}, 35(2):177--204,
  1999.

\bibitem{he2018local}
Y.~He, A.~Knowles, and M.~Marcozzi.
\newblock Local law and complete eigenvector delocalization for supercritical
  {E}rd{\H{o}}s--{R}{\'{e}}nyi graphs.
\newblock {\em arXiv preprint arXiv:1808.09437}, 2018.

\bibitem{HLY15}
J.~Huang, B.~Landon, and H.-T. Yau.
\newblock Bulk universality of sparse random matrices.
\newblock {\em J. Math. Phys.}, 56(12):123301, 2015.

\bibitem{fixed}
B.~Landon, P.~Sosoe, and H.-T. Yau.
\newblock Fixed energy universality for {D}yson {B}rownian motion.
\newblock {\em Preprint arXiv:1609.09011}, 2016.

\bibitem{LY}
B.~Landon and H.-T. Yau.
\newblock Convergence of local statistics of dyson brownian motion.
\newblock {\em Comm. Math. Phys.}, 355:949–1000, Nov. 2017.

\bibitem{schnelli2}
J.~O. Lee, K.~Schnelli, B.~Stetler, and H.-T. Yau.
\newblock Bulk universality for deformed {W}igner matrices.
\newblock {\em Ann. Probab.}, 44(3):2349--2425, 2016.

\bibitem{lieb2001analysis}
E.~H. Lieb and M.~Loss.
\newblock Analysis, volume 14 of graduate studies in mathematics.
\newblock {\em American Mathematical Society, Providence, RI,}, 4, 2001.

\bibitem{LR12}
A.~Litvak and O.~Rivasplata.
\newblock Smallest singular value of sparse random matrices.
\newblock {\em Studia Mathematica}, 212(3), Jun. 2011.

\bibitem{MP}
V.~A. Marchenko and L.~A. Pastur.
\newblock The distribution of eigenvalues in certain sets of random matrices.
\newblock {\em Mat. Sb.}, 72:507--536, 1967.

\bibitem{PY}
N.~S. Pillai and J.~Yin.
\newblock Universality of covariance matrices.
\newblock {\em The Annals of Applied Probability}, 24(3):935--1001, 2014.

\bibitem{RV}
M.~Rudelson and R.~Vershynin.
\newblock Non-asymptotic theory of random matrices: extreme singular values.
\newblock {\em arXiv preprint arXiv:1003.2990}, 2010.

\bibitem{ST02}
D.~Spielman and S.~Teng.
\newblock Smoothed analysis of algorithms.
\newblock In {\em Proceedings of the International Congress of Mathematicians},
  volume~I, pages 597--606, Beijing, 2002. Higher Ed. Press.

\bibitem{TV09a}
T.~Tao and V.~Vu.
\newblock Inverse {L}ittlewood-{O}fford theorems and the condition number of
  random discrete matrices.
\newblock {\em Annals of Mathematics}, 169:595--632, 2009.

\bibitem{leastsv}
T.~Tao and V.~Vu.
\newblock Random matrices: the distribution of the smallest singular value.
\newblock {\em Geometric and Functional Analysis}, 20(1):260--297, 2010.

\bibitem{TV10a}
T.~Tao and V.~Vu.
\newblock Smooth analysis of the condition number and the least singular value.
\newblock {\em Mathematics of computation}, 79(272):2333--2352, 2010.

\bibitem{rigidity}
J.~Yin, L.~Erd{\H{o}}s, and H.-T. Yau.
\newblock Rigidity of eigenvalues of generalized {W}igner matrices.
\newblock {\em Adv. Math}, 229(3):1435--1515, 2012.

\end{thebibliography}
\bibliographystyle{abbrv}

\end{document}